\documentclass[a4paper,11pt]{article}
\usepackage{amsmath}
\usepackage{amssymb}
\usepackage{amsthm}
\usepackage{MnSymbol}

\usepackage{mathrsfs}
\usepackage{graphicx,subfigure}
\usepackage{paralist} %enumerate/itemize in paragraphs
\usepackage{enumerate}

\usepackage{hyperref}
\usepackage{bm} % for bold greek letters
\usepackage{placeins}
\usepackage{color}
\usepackage[margin=30mm]{geometry}

\newcommand{\R}{\mathbb{R}}

\newcommand{\dHaus}{\, \mathrm{d}\mathcal{H}^{d-1}}
\newcommand{\dx}{\, \mathrm{dx}}

\newcommand{\dt}{\, \mathrm{dt}}
\newcommand{\dd}{\, \mathrm{d}}
\newcommand{\dz}{\, \mathrm{dz}}

\newcommand{\pd}{\partial}
\newcommand{\abs}[1]{\left| #1 \right|}

\newcommand{\eps}{\varepsilon}
\newcommand{\Laplace}{\Delta}

\newcommand{\md}{\pd^{\bullet}_{t}}

\newcommand{\sech}{\mathrm{sech}}

\newcommand{\jump}[1]{\left [#1 \right ]_{H}^{T}}

\newcommand{\para}{g} %for the parameterisation of the interface 

\renewcommand{\div}{\, \mathrm{div} \, }

\allowdisplaybreaks

\theoremstyle{plain}
\newtheorem{rmk}{Remark}[section]

\newtheorem{assump}{Assumption}[section]
\numberwithin{equation}{section}
\title{A Cahn--Hilliard--Darcy model for tumour growth with chemotaxis and active transport}

\author{Harald Garcke \footnotemark[1] \and Kei Fong Lam\footnotemark[1] \and Emanuel Sitka \footnotemark[2] \and Vanessa Styles \footnotemark[3]}

\date{\today}

\begin{document}

\maketitle

\renewcommand{\thefootnote}{\fnsymbol{footnote}}
\footnotetext[1]{Fakult\"at f\"ur Mathematik, Universit\"at Regensburg, 93040 Regensburg, Germany
({\tt \{Harald.Garcke, Kei-Fong.Lam\}@mathematik.uni-regensburg.de}).}
\footnotetext[2]{Fakult\"at f\"ur Medizin, Universit\"at Regensburg, 93040 Regensburg, Germany ({\tt Emanuel.Sitka@stud.uni-regensburg.de}).}
\footnotetext[3]{University of Sussex, Sussex House, Falmer, Brighton, BN1 9RH, United Kingdom
({\tt V.Styles@sussex.ac.uk}).}

\begin{abstract}
Using basic thermodynamic principles we derive a Cahn--Hilliard--Darcy model for tumour growth including nutrient diffusion, chemotaxis, active transport, adhesion, apoptosis and proliferation.  The model generalises earlier models and in particular includes active transport mechanisms which ensure thermodynamic consistency.  We perform a formally matched asymptotic expansion and develop several sharp interface models.  Some of them are classical and some are new which for example include a jump in the nutrient density at the interface.  A linear stability analysis for a growing nucleus is performed and in particular the role of the new active transport term is analysed.  Numerical computations are performed to study the influence of the active transport term for specific growth scenarios.
\end{abstract}

\noindent \textbf{Key words. } Tumour growth, diffuse interface model, Cahn--Hilliard equation, chemotaxis, Darcy's flow, matched asymptotic expansions, stability analysis, finite element computations.\\

\noindent \textbf{AMS subject classification. } 92B05, 35K57, 35R35, 65M60

\section{Introduction}\label{sec:Intro}
In the last decades the understanding of tumour related illnesses has undergone a swift development.  Nowadays tumour therapy can be adapted to the genetic fingerprint of the tumour, resulting in a ``targeted therapy" that has dramatically improved the prognosis of many illnesses.  While some important mutations in tumour genomes have been identified and exploited by modern tumour drugs, basic growth behaviours of tumours are still far from being understood, e.g. angiogenesis and the formation of metastases.  The complexity of oncology has also attracted increasing interest of mathematicians, who are trying to find the appropriate equations to  provide additional insights in certain aspects of tumour growth, see for example \cite{article:BellomoLiMaini08} and \cite{book:CristiniLowengrub}.  In this paper we want to introduce a new diffuse interface model for tumour growth, and compare the resulting system of partial differential equations to some other recent contributions
\cite{article:ColliGilardiHilhorst15,article:CristiniLiLowengrubWise09,article:FrigeriGrasselliRocca15,article:HawkinsPrudhommevanderZeeOden13,article:HawkinsZeeOden12,article:HilhorstKaampmannNguyenZee15,preprint:JiangWuZheng14,article:LowengrubTitiZhao13,article:OdenHawkinsPrudhomme10,thesis:Sitka,article:WiseLowengrubFrieboesCristini08}.

In order to obtain a tractable system of partial differential equations, we will in this paper neglect some effects which could be addressed in further research and which then would lead to more complete theories.  From a medical point of view we will hence make the following assumptions as foundations for our further considerations:
\begin{enumerate}
\item Tumour cells only die by apoptosis. Hence we neglect the possibility of tumour necrosis, where we would have to take account of the negative effects of chemical species from the former intracellular space on the surrounding tumour cells.
\item The tissue around the tumour does not react to the tumour cells in any active way. In particular, we neglect any response of the immune system to the tumour tissue. 
\item Larger tumour entities are actually enforcing blood vessel growth towards themselves by secreting vessel growth factors.  This is a phenomenon that could be addressed in future in a generalised model.
\item We postulate the existence of an unspecified chemical species acting as a nutrient for the tumour cells.  This nutrient is not consumed by the healthy tissue.  We will also introduce terms which will reflect chemotaxis, which is the active movement of the tumour colony towards nutrient sources.  Additionally, the introduction of chemotaxis will also lead to the opposite process, meaning that the nutrient is moving towards the nearby tumour cells.  As we will point out later, this could be seen as a correlate of a nutrient uptake mechanism.
\end{enumerate}

Here we state a slightly simplified version of the general system, which will be derived in Section \ref{sec:derivation} from thermodynamic principles.  We will derive and analyse a two-component mixture model of tumour and healthy cells, whose behaviour is governed by the system
\begin{subequations}\label{Intro:CHDarcy:zeroexcesstotalmass}
\allowdisplaybreaks[3]
\begin{align}
\div \bm{v} & =  \alpha \Gamma, \label{Intro:div} \\
\bm{v} &= -K(\nabla p - \mu \nabla \varphi - \chi_{\varphi} \sigma \nabla \varphi), \label{Intro:Darcy} \\
\pd_{t}\varphi + \div (\bm{v} \varphi) & = \nabla \cdot (m(\varphi)\nabla \mu) + \bar{\rho}_{S} \Gamma, \label{Intro:varphi} \\
\mu & = \tfrac{\beta}{\eps} \Psi'(\varphi) - \beta\eps \Laplace \varphi - \chi_{\varphi}\sigma, \label{Intro:mu} \\
\pd_{t}\sigma + \div ( \sigma \bm{v} )& = \div( n(\varphi)(\chi_{\sigma} \nabla \sigma-\chi_{\varphi} \nabla \varphi)) - \mathcal{C} \sigma h(\varphi), \label{Intro:sigma} \\
\Gamma & = (\mathcal{P}\sigma - \mathcal{A}) h(\varphi). \label{Intro:MassTrans}
\end{align}
\end{subequations}
Here, $\bm{v}$ denotes the volume-averaged velocity of the mixture, $p$ denotes the pressure, $\sigma$ denotes the concentration of an unspecified chemical species that serves as a nutrient for the tumour, $\varphi \in [-1,1]$ denotes the difference in volume fractions, with $\{\varphi = 1\}$ representing unmixed tumour tissue, and $\{\varphi = -1\}$ representing the surrounding healthy tissue, and $\mu$ denotes the chemical potential for $\varphi$.  The particular simple form of (\ref{Intro:div}) is different to earlier modelling attempts and is based on the fact that we use volume-averaged velocities.

The positive constants $K$, $\beta$, $\mathcal{P}$, $\mathcal{A}$, and $\mathcal{C}$ denote the permeability, surface tension, proliferation rate, apoptosis rate, and consumption rate, respectively.  The constants $\overline{\rho}_{S}$ and $\alpha$ are related to the densities of the two components (see (\ref{defn:alpharhoS}) below), in particular, for the case of matched densities we have $\alpha = 0$.  Meanwhile $m(\varphi)$ and $n(\varphi)$ are non-negative mobilities for $\varphi$ and $\sigma$, respectively, and $\Psi(\cdot)$ is a potential with two equal minima at $\pm 1$.  In addition, we choose $h$ as an interpolation function with $h(-1) = 0$ and $h(1) = 1$.  The simplest choice is given as $h(\varphi) = \frac{1}{2}(1+\varphi)$.

We denote $\chi_{\sigma} \geq 0$ as the diffusivity of the nutrient, and $\chi_{\varphi} \geq 0$ can be seen as a parameter for transport mechanisms such as chemotaxis and active uptake (see below for more details).  Finally, the parameter $\eps$ is related to the thickness of the interfacial layers present in phase field systems.  The system (\ref{Intro:CHDarcy:zeroexcesstotalmass}) is a Cahn--Hilliard--Darcy system coupled to a convection-diffusion-reaction equation for the nutrient.

(\ref{Intro:div}) and (\ref{Intro:Darcy}) model the mass balance using a Darcy-type system, and in the situation of unmatched densities ($\alpha \neq 0$), the gain and loss of volume resulting from the mass transition $\Gamma$ leads to sources and sinks in the mass balance.  In (\ref{Intro:varphi}) and (\ref{Intro:mu}), $\varphi$ is governed by a Cahn--Hilliard type equation with additional source terms.  The mass transition from the the healthy cells to the tumour component and vice versa is described in (\ref{Intro:MassTrans}), where tumour growth/proliferation is represented by the term $\mathcal{P}\sigma h(\varphi)$, and the process of apoptosis is modelled by the term $\mathcal{A} h(\varphi)$.  In (\ref{Intro:sigma}), the nutrient is subjected to an equation of convection-reaction-diffusion type, and the term $\mathcal{C}\sigma h(\varphi)$ represents consumption of the nutrient only in the presence of the tumour cells.  As in \cite{article:ByrneChaplain95}, we could also consider the situation where the tumour possesses its own vasculature and the nutrient may be supplied to the tumour via a capillary network at a rate $\mathcal{B} (\sigma - \sigma_{B})$, where $\sigma_{B}$ is the constant nutrient concentration in the vasculature and $\mathcal{B}$ is the blood-tissue transfer rate which might depend on $\varphi$ and $x$.  This leads to the following nutrient balance equation instead of (\ref{Intro:sigma})
\begin{align*}
\pd_{t} \sigma + \div (\sigma \bm{v}) = \div (n(\varphi) (\chi_{\sigma} \nabla \sigma - \chi_{\varphi} \nabla \varphi)) - \mathcal{C} \sigma h(\varphi) + \mathcal{B} (\sigma_{B} - \sigma).
\end{align*}
Under appropriate boundary conditions the system (\ref{Intro:CHDarcy:zeroexcesstotalmass}) allows for an energy inequality (see (\ref{Model:CHDarcy:energyeq}) below) and we believe that this inequality will allow the well-posedness of the above system to be rigorously shown.

We now motivate the particular choices for the modelling of proliferation, apoptosis, chemotaxis, and mass transition in (\ref{Intro:CHDarcy:zeroexcesstotalmass}).
\begin{itemize}
\item In (\ref{Intro:MassTrans}), we obtain that $\Gamma = \mathcal{P} \sigma - \mathcal{A}$ holds in the tumour region $\{ \varphi = 1 \}$.  The implicit assumption that the tumour growth is proportional to the nutrient supply can be justified by the fact that malign tumours have the common genetic feature that certain growth inhibiting proteins have been switched off by mutations.  Hence, we can assume that while in healthy cells the mitotic cycle is rather strictly inhibited, tumour cells often show unregulated growth behaviour which is only limited by the supply of nutrients. 

Moreover, implicit in the choice of zero mass transition $\Gamma = 0$ in the healthy region $\{\varphi = -1\}$ is the assumption that the tumour proliferation rate is more significant than that of the healthy tissue.
\item In (\ref{Intro:varphi}) and (\ref{Intro:sigma}), the fluxes for $\varphi$ and $\sigma$ are given by
\begin{align*}
\bm{q}_{\varphi} & := -m(\varphi) \nabla \mu = -m(\varphi) \nabla \left ( \tfrac{\beta}{\eps} \Psi'(\varphi) - \beta \eps \Laplace \varphi - \chi_{\varphi} \sigma \right ), \\
\bm{q}_{\sigma} & := -n(\varphi)\nabla (\chi_{\sigma} \sigma - \chi_{\varphi}  \varphi),
\end{align*}
respectively.  It has been pointed out by Roussos, Condeelis and Patsialou in \cite{article:RoussosCondeelisPatsialou11}  that the undersupply of nutrient induces chemotaxis in certain tumour entities.  This is reflected in the  term $m(\varphi) \nabla (\chi_{\varphi} \sigma)$ of $\bm{q}_{\varphi}$, which drives the cells towards regions of high nutrient.  

On the other hand, we note that the term $n(\varphi) \nabla (\chi_{\varphi} \varphi)$ in $\bm{q}_{\sigma}$ drives the nutrient to regions of high $\varphi$, i.e., to the tumour cells, which indicates that the nutrient is actively moving towards the tumour cells.  This may seem to be counter-intuitive at first glance.  However, this term will only contribute to the equation significantly in the vicinity of the interface between the tumour and healthy cells.  This allows the interpretation that the term $n(\varphi) \nabla (\chi_{\varphi} \varphi)$ reflects active transport mechanisms which move the nutrient into the tumour colony.  Here we use the term ``active transport'' in the biological sense in order to indicate that some kind of mechanism is needed to maintain the transport (in contrast to passive transporters which are driven only by the concentration gradient of the substance).  The additional mechanism allows cells to establish persisting concentration differences between different compartments.  In particular we can expect that tumours, which have these active transporters on their cell membrane, are not dependent on diffusion but can establish high concentration of the vital nutrient even against the nutrient concentration gradient. 
\end{itemize}

Here we briefly give an example, where mechanisms like this have already been observed: Malign tumour cells often have a significantly increased need for glucose, a fact that is sometimes referred to as the Warburg effect.  As a consequence of several mutations in the tumour genome, these cells can adapt to their high rate of glucose consumption in several ways.  Apart from angiogenesis, which leads to a well perfused tumour environment providing large amount of glucose, the tumour cells can also express (i.e., build) more glucose transporters, which provide an improved glucose transport through the cell membrane.  Recently, both passive glucose transporters, so-called GLUT proteins, and active glucose transporters called SGLTs, have been observed on the cell membrane of several tumour entities.  For a more detailed description regarding the GLUT transporters we refer to \cite{article:Calvoetal10}, whereas SGLT expression of tumours has been described by \cite{article:Ishikawaetal01} and \cite{article:Scafoglioetal15}.  In the system we will derive in this paper, the existence of passive nutrient transporters like the GLUTs is implicitly assumed by including nutrient diffusion.  Apart from that, it will become more obvious in the corresponding sharp interface system that the term $n(\varphi) \nabla (\chi_{\varphi} \varphi)$ represents an active nutrient transport towards the tumour.

We note that in (\ref{Intro:CHDarcy:zeroexcesstotalmass}), the mechanism of chemotaxis and active transport are connected via the parameter $\chi_{\varphi}$.  In principle, it is possible to decouple the two mechanisms.  In order to do so, we introduce the following choice for the mobility $n(\varphi)$ and diffusion coefficient $\chi_{\sigma}$ (see also Section \ref{sec:compare:Cristini} below):  For $\lambda > 0$ and a non-negative mobility $\mathcal{D}(\varphi)$, we set
\begin{align}\label{activetransportchoice}
n(\varphi) = \lambda \mathcal{D}(\varphi) \chi_{\varphi}^{-1}, \quad \chi_{\sigma} = \lambda^{-1} \chi_{\varphi}.
\end{align}
Then, the corresponding fluxes for $\varphi$ and $\sigma$ are now given as
\begin{subequations}\label{activetransportflux}
\begin{align}
\bm{q}_{\varphi} & = -m(\varphi) \nabla \left ( \tfrac{\beta}{\eps} \Psi'(\varphi) - \beta \eps \Laplace \varphi - \chi_{\varphi} \sigma \right ), \\
\bm{q}_{\sigma} & = -\mathcal{D}(\varphi)\nabla (\sigma - \lambda \varphi).
\end{align}
\end{subequations}
For this choice, we can switch off the effects of active transport by sending $\lambda \to 0$, while preserving the effects of chemotaxis.

We now compare the new model (\ref{Intro:CHDarcy:zeroexcesstotalmass}) and some of the previous diffuse interface models in the literature:
\begin{itemize}
\item An important difference to other models is that we take a volume-averaged velocity which leads to the simple form $\div \bm{v} = \alpha \Gamma$, as a consequence of mass balance.  In other models, this equation has to be replaced by a more complicated transport equation.
\item Another significant difference is the presence of the term $-\div(n(\varphi)\chi_{\varphi}\nabla \varphi)$ in (\ref{Intro:sigma}), which only appears in cases where chemotaxis and active transport are taken into account.  As we have pointed out before, it represents active nutrient transport towards the tumour.  The corresponding nutrient equations in \cite{article:ColliGilardiHilhorst15,article:CristiniLiLowengrubWise09,article:FrigeriGrasselliRocca15,article:HawkinsPrudhommevanderZeeOden13,article:HilhorstKaampmannNguyenZee15,article:OdenHawkinsPrudhomme10,article:WiseLowengrubFrieboesCristini08} do not include an equivalent term.  However, we point out that this active transport mechanism is present in the nutrient equation of \cite{article:HawkinsZeeOden12}, who however used different source terms and no Darcy-flow contributions.

\item Our choice of the mass transition term $\Gamma$ in (\ref{Intro:MassTrans}) can also be found in \cite{article:CristiniLiLowengrubWise09,article:HawkinsPrudhommevanderZeeOden13,article:OdenHawkinsPrudhomme10,article:WiseLowengrubFrieboesCristini08}. 
Alternatively, one may consider equations of the form
\begin{align*}
\pd_{t} \varphi & = \div(m(\varphi) \nabla \mu) + P(\varphi) (\sigma - \chi \varphi - \mu), \\
\pd_{t} \sigma & = \div (n(\varphi) (\nabla \sigma - \chi \nabla \varphi)) - P(\varphi) (\sigma - \chi \varphi - \mu),  
\end{align*}
where the chemical potential $\mu$ enters as a source term for the equations of $\varphi$ and $\sigma$.  Here $\chi \geq 0$ is a constant, and $P(\cdot)$ denotes a non-negative proliferation function.  This type of mass transition term appears in \cite{article:HawkinsZeeOden12} and in \cite{article:ColliGilardiHilhorst15,article:FrigeriGrasselliRocca15,article:HilhorstKaampmannNguyenZee15} with $\chi = 0$.

\item The presence of chemotaxis, represented by the term $-\chi_{\varphi} \sigma$ in (\ref{Intro:mu}) can also be found in the models of \cite{article:CristiniLiLowengrubWise09,article:HawkinsPrudhommevanderZeeOden13,article:HawkinsZeeOden12,article:OdenHawkinsPrudhomme10}, while the corresponding Cahn--Hilliard systems in \cite{article:ColliGilardiHilhorst15,article:FrigeriGrasselliRocca15,article:HilhorstKaampmannNguyenZee15,preprint:JiangWuZheng14,article:LowengrubTitiZhao13,article:WiseLowengrubFrieboesCristini08} do not include an equivalent term.

\item In \cite{preprint:JiangWuZheng14,article:LowengrubTitiZhao13,article:WiseLowengrubFrieboesCristini08}, the nutrient does not enter the Darcy law for $\bm{v}$ like in (\ref{Intro:Darcy}).
\end{itemize}

In the diffuse interface model (\ref{Intro:CHDarcy:zeroexcesstotalmass}), the parameter $\eps$ is related to the thickness of the interfacial layer, which separates the tumour cell regions $\{ \varphi = 1 \}$ and the healthy cell regions $\{ \varphi = - 1\}$.  Hence, it is natural to ask if a sharp interface description of the problem will emerge in the limit $\eps \to 0$.  This means in the limit the interface between the tumour cells and the healthy cells is represented by a hypersurface of zero thickness.  

For convenience, suppose we take the mobilities $m(\varphi) = m_{0}$, $n(\varphi) = n_{0}$ to be constant.  A formally matched asymptotic analysis will yield the following sharp interface limit from (\ref{Intro:CHDarcy:zeroexcesstotalmass}) (see Section \ref{sec:Asym} for more details):  Let $\Omega_{T}$ and $\Omega_{H}$ denote the tumour cell region and the healthy cell region, respectively, which are separated by an interface $\Sigma$.  Then it holds that
\begin{subequations}\label{Intro:SharpInterface}\allowdisplaybreaks[3]
\begin{align}
\bm{v} & = -K \nabla p \text{ in } \Omega_{T} \cup \Omega_{H}, \\
\div \bm{v} & = \begin{cases}
\alpha(\mathcal{P} \sigma - \mathcal{A}) & \text{ in } \Omega_{T}, \\
0 & \text{ in } \Omega_{H}, 
\end{cases} \\
- m_{0} \Laplace \mu & = \begin{cases}(\overline{\rho}_{S} - \alpha)(\mathcal{P} \sigma_{0} - \mathcal{A}) & \text{ in } \Omega_{T}, \\
0 & \text{ in } \Omega_{H}, 
\end{cases} \\
\pd_{t} \sigma + \div ( \sigma \bm{v}) - n_{0}\chi_{\sigma} \Laplace \sigma & = \begin{cases}
-\mathcal{C} \sigma & \text{ in } \Omega_{T}, \\
0 & \text{ in } \Omega_{H}, 
\end{cases} \\
\jump{\bm{v}} \cdot \bm{\nu} = 0, & \quad \jump{\sigma} = 2 \frac{\chi_{\varphi}}{\chi_{\sigma}}, \quad \jump{p} = \beta \gamma \kappa  \text{ on } \Sigma, \label{Intro:SharpInterface:jump}   \\
\jump{\mu} = 0, & \quad 2\mu + \frac{\chi_{\sigma}}{2} \jump{\abs{\sigma}^{2}} = \beta \gamma \kappa  \text{ on } \Sigma, \\
2(-\mathcal{V} + \bm{v} \cdot \bm{\nu}) & = m_{0} \jump{\nabla \mu_{0}} \cdot \bm{\nu} \text{ on } \Sigma, \\
2 \frac{\chi_{\varphi}}{\chi_{\sigma}} (-\mathcal{V} + \bm{v} \cdot \bm{\nu})& = n_{0} \jump{\nabla \sigma} \cdot \bm{\nu}  \text{ on } \Sigma. 
\end{align}
\end{subequations}
Here, $\gamma$ is a constant related to the potential $\Psi$ (see (\ref{defn:gammaconst}) below), $\mathcal{V}$ denotes the normal velocity of $\Sigma$, $\kappa$ is the mean curvature of $\Sigma$, $\jump{f}$ denotes the jump of $f$ from $\Omega_{T}$ to $\Omega_{H}$ across $\Sigma$ (see (\ref{defn:jump})), and $\bm{\nu}$ is the outward unit normal of $\Sigma$, pointing towards $\Omega_{T}$.

In comparison to the formal sharp interface limits of \cite{article:CristiniLiLowengrubWise09,article:HilhorstKaampmannNguyenZee15,article:WiseLowengrubFrieboesCristini08}, the most significant difference is the jump condition $\eqref{Intro:SharpInterface:jump}_{2}$.  Let us remark on its physical meaning.  Let $\sigma_{T}$ and $\sigma_{H}$ denote the limiting values of the nutrient on the interface $\Sigma$ from tumour cell regions and from the healthy cell regions, respectively.  Then, $\eqref{Intro:SharpInterface:jump}_{2}$ implies that
\begin{align*}
\sigma_{T} = \sigma_{H} + 2 \frac{\chi_{\varphi}}{\chi_{\sigma}}.
\end{align*}
Thus, if $\chi_{\varphi}$ is positive, then $\eqref{Intro:SharpInterface:jump}_{2}$ tells us that the tumour cells will experience a higher level of nutrient concentration than the healthy cells on the interface, which reflects the effect of the active transport mechanism in (\ref{Intro:sigma}), attracting nutrients from the healthy cell regions into the tumour.

If we consider the fluxes (\ref{activetransportflux}) in (\ref{Intro:CHDarcy:zeroexcesstotalmass}), then one obtains the sharp interface model (\ref{Intro:SharpInterface}) with the following modification (see Section \ref{sec:compare:Cristini} for more details):  Instead of $(\ref{Intro:SharpInterface:jump})_{2}$, we now have
\begin{align*}
\jump{\sigma} = 2 \lambda.
\end{align*}
In particular, the parameter $\lambda$ only enters explicitly in the jump condition for $\sigma$, which relates to the above discussion regarding the physical interpretation of $(\ref{Intro:SharpInterface:jump})_{2}$.  Note that $\lambda \neq 0$ is a consequence of the active transport term we have discussed in the phase field model.  In the matched asymptotics expansion, this term directly leads to a jump of the nutrient concentration at the tumour interface.  Therefore, we obtain exactly the situation one would expect from active transport mechanisms: Close to the tumour surface, we observe a higher nutrient concentration inside the tumour than on the outside of the tumour, a situation that is only possible due to the transporter molecules.  Hence it should be considered if $\lambda$ could be referred to as a density parameter for the active transport proteins.

The plan of this paper is as follows:  In Section \ref{sec:derivation} we derive the new phase field model from thermodynamic principles and compare with previous phase field models of tumour growth in the literature.  In Section \ref{sec:Asym} we perform a formal asymptotic analysis to derive certain sharp interface models of tumour growth.  In Section \ref{sec:stabAna} we investigate the stability of radial solutions to a particular sharp interface model via a linear stability analysis, and highlight the effect of the active transport parameter on the stability.  In Section \ref{sec:Numerics} we present quantitative simulations for radially symmetric solutions and qualitative simulations for more general scenarios.

\section{Model Derivation}\label{sec:derivation}
Let us consider a two component mixture consisting of tumour and healthy cells in a bounded domain $\Omega \subset \R^{d}$, $d = 1,2,3$.  

We denote the first component as the component of healthy tissues, and the second component as the tumour tissues.  Let $\rho_{i}$, $i = 1,2$, denote the actual mass of the component matter per volume in the mixture, and let $\bar{\rho}_{i}$, $i = 1,2,$ be the mass density of a pure component $i$.  Then, $\rho := \rho_{1} + \rho_{2}$ denotes the mixture density (which is not necessarily constant), and we define the volume fraction of component $i$ as
\begin{align}\label{defn:volumefrac}
u_{i} = \frac{\rho_{i}}{\bar{\rho}_{i}}.
\end{align}
We expect that physically, $\rho_{i} \in [0, \bar{\rho}_{i}]$ and thus $u_{i} \in [0,1]$.  In addition to the considerations stated in Section \ref{sec:Intro}, we make the following modelling assumptions:
\begin{itemize}
\item There is no external volume compartment besides the two components, i.e.,
\begin{align}
\label{uisum1}
u_{1} + u_{2} = 1.
\end{align}
\item We allow for mass exchange between the two components.  Growth of the tumour is represented by mass transfer from component 1 (healthy tissues) to component 2 (tumour tissues), while tumour cells are converted back into the surrounding healthy tissues when they die.
\item We choose the mixture velocity to be the volume-averaged velocity:
\begin{align}
\label{defn:volavervelo}
\bm{v} := u_{1} \bm{v}_{1} + u_{2} \bm{v}_{2},
\end{align}
where $\bm{v}_{i}$ is the individual velocity of component $i$.
\item We model a general chemical species which is treated as a nutrient for the tumour tissues.  Its concentration is denoted by $\sigma$ and it is transported by the volume-averaged mixture velocity and a flux $\bm{J}_{\sigma}$.
\end{itemize}

\subsection{Balance laws}
The balance law for mass of each component reads as
\begin{subequations}\label{proto:individualmass}
\begin{align}
\pd_{t}\rho_{1} + \div (\rho_{1} \bm{v}_{1}) & = \Gamma_{1}, \\
\pd_{t}\rho_{2} + \div (\rho_{2} \bm{v}_{2}) & = \Gamma_{2}.
\end{align}
\end{subequations}
Observe that by (\ref{defn:volumefrac}), we can write (\ref{proto:individualmass}) in the following way: For $i = 1,2$,
\begin{align}
\label{proto:individualmass2}
\pd_{t} u_{i} + \div (u_{i} \bm{v}_{i}) = \frac{\Gamma_{i}}{\bar{\rho}_{i}}.
\end{align} 
We see that by (\ref{uisum1}),  (\ref{defn:volavervelo}), and (\ref{proto:individualmass2}),
\begin{align}
\label{divv}
\div \bm{v} = \div (u_{1} \bm{v}_{1}) + \div (u_{2} \bm{v}_{2}) = \sum_{i=1}^{2} \left (\frac{\Gamma_{i}}{\bar{\rho}_{i}} -\pd_{t}u_{i} \right ) = \frac{\Gamma_{2}}{\bar{\rho}_{2}} + \frac{\Gamma_{1}}{\bar{\rho}_{1}} =: \Gamma_{\bm{v}}.
\end{align}
We introduce the fluxes:
\begin{align}\label{defn:fluxJi}
\bm{J}_{i} := \rho_{i} (\bm{v}_{i} - \bm{v}), \quad \bm{\mathcal{J}} := \bm{J}_{1} + \bm{J}_{2}, \quad  \bm{J} := -\frac{1}{\overline{\rho}_{1}} \bm{J}_{1} + \frac{1}{\overline{\rho}_{2}} \bm{J}_{2}.
\end{align}
Then, we see that 
\begin{align*}
\bm{\mathcal{J}} + \rho \bm{v} = \bm{J}_{1} + \bm{J}_{2} + \rho \bm{v} = \rho_{1} \bm{v}_{1} + \rho_{2} \bm{v}_{2},
\end{align*}
and so, upon adding the equations in (\ref{proto:individualmass}) we obtain the equation for the mixture density:
\begin{align}
\label{proto:mixturemass}
\pd_{t}\rho + \div (\rho_{1} \bm{v}_{1} + \rho_{2} \bm{v}_{2}) = \pd_{t}\rho + \div (\rho \bm{v} + \bm{\mathcal{J}}) = \Gamma_{1} + \Gamma_{2}.
\end{align}
We now want to derive an equation for the phase field variable $\varphi$.  Recalling $\rho_{i} = \bar{\rho}_{i} u_{i}$, we obtain from (\ref{proto:individualmass2}) that
\begin{align}
\label{proto:individualmass3}
\pd_{t}u_{i} + \frac{1}{\bar{\rho}_{i}}\div \bm{J}_{i} + \div (u_{i} \bm{v}) = \frac{\Gamma_{i}}{\bar{\rho}_{i}}.
\end{align}
We define the order parameter $\varphi$ as the difference in volume fractions:
\begin{align}
\label{defn:orderparameter}
\varphi := u_{2} - u_{1},
\end{align}
then, subtracting the equation for $u_{1}$ from the equation for $u_{2}$, and using (\ref{defn:fluxJi}), we obtain the equation for $\varphi$:
\begin{align}
\label{orderparameter:equ}
\pd_{t}\varphi + \div (\varphi \bm{v}) + \div \bm{J} =  \frac{\Gamma_{2}}{\overline{\rho}_{2}} - \frac{\Gamma_{1}}{\overline{\rho}_{1}} =: \Gamma_{\varphi}.
\end{align}
We point out that from the constraint (\ref{uisum1}), we obtain
\begin{align*}
u_{2} = \frac{1 + \varphi}{2}, \quad u_{1} = \frac{1 - \varphi}{2}.
\end{align*}
Thus, the region of the tumour tissues is represented by $\{ x \in \Omega : \varphi = 1\}$ and the region of healthy tissues is represented by $\{ x \in \Omega : \varphi = -1\}$.  In particular, the mixture density $\rho$ can be expressed as a linear function of $\varphi$:
\begin{align}
\label{mixturedensity}
\rho = \rho_{1} + \rho_{2} = \bar{\rho}_{1} \frac{1-\varphi}{2} + \bar{\rho}_{2} \frac{1+\varphi}{2} = \frac{\bar{\rho}_{1} + \bar{\rho}_{2}}{2} + \varphi \frac{\bar{\rho}_{2} - \bar{\rho}_{1}}{2}.
\end{align}
For the nutrient, we postulate the following balance law:
\begin{align}
\label{proto:nutrient}
\pd_{t}\sigma + \div (\sigma \bm{v}) + \div \bm{J}_{\sigma} = -\mathcal{S},
\end{align}
where $\mathcal{S}$ denotes a source/sink term for the nutrient.  In addition, $\sigma \bm{v}$ models the transport by the volume-averaged velocity and $\bm{J}_{\sigma}$ models other transport mechanisms like diffusion and chemotaxis.

\subsection{Energy inequality}\label{sec:energyineq}
We postulate a general energy density of the form:
\begin{align}
\label{defn:energy}
e(\varphi, \nabla \varphi, \sigma) = f(\varphi, \nabla \varphi) + N(\varphi, \sigma).
\end{align}
Here, we neglected inertia effects, and so the kinetic energy does not appear in $e$.  Instead we refer the reader to \cite{thesis:Sitka} for the derivation of a model that includes inertia effects, leading to a Navier--Stokes--Cahn--Hilliard version of (\ref{Intro:CHDarcy:zeroexcesstotalmass}).  The first term $f$ in (\ref{defn:energy}) accounts for interfacial energy and unmixing tendencies, while the second term $N$ describes the chemical energy of the nutrient and energy contributions resulting from the interactions between the tumour tissues and the nutrient.  The latter will, for example, lead to chemotatic effects which are of particular interest as they result in the tumour tissue growing towards regions with high nutrient concentration.

In the following, we will consider $f$ to be of Ginzburg-Landau type:  For constants $A, B > 0$, we choose
\begin{align}\label{GinzburgLandau}
f(\varphi, \nabla \varphi) := A \Psi(\varphi) + \frac{B}{2} \abs{\nabla \varphi}^{2},
\end{align}
where $\Psi(s)$ is a potential with equal minima at $s = \pm 1$.

We will now derive the diffuse interface model based on a dissipation inequality for balance laws with source terms which has been used similarly by Gurtin \cite{article:Gurtin89,article:Gurtin96} and Podio-Guidugli \cite{article:Podio06} to derive phase field and Cahn--Hilliard type equations.  These authors used the second law of thermodynamics which in an isothermal situation is formulated as a free energy inequality.  We also refer to Chapter 62 of Gurtin, Fried, and Anand \cite{book:GurtinFriedAnand} for a detailed discussion of situations with source terms.  The second law of thermodynamics in the isothermal situation requires that for all volumes $V(t) \subset \Omega$ which are transported with the fluid velocity the following inequality has to hold (see \cite{article:Gurtin89,article:Gurtin96,book:GurtinFriedAnand,article:Podio06}):
\begin{align*}
\frac{\dd}{\dt} \int_{V(t)} e \dx \leq - \int_{\pd V(t)} \bm{J}_{e} \cdot \bm{\nu} \dHaus + \int_{V(t)} c_{\varphi} \Gamma_{\varphi} + c_{\bm{v}} \Gamma_{\bm{v}} + c_{S} (-\mathcal{S})  \dx,
\end{align*}
where $\bm{\nu}$ is the outer unit normal to $\pd V(t)$ and $\bm{J}_{e}$ is an energy flux yet to be specified.  Following \cite{book:GurtinFriedAnand}, we postulate that the source terms $\Gamma_{\bm{v}}$, $\Gamma_{\varphi}$ and the nutrient supply $(-\mathcal{S})$ carry with them a supply of energy described by
\begin{align}
\label{GurtinEnergySupply}
\int_{V(t)} c_{\bm{v}} \Gamma_{\bm{v}} + c_{\varphi} \Gamma_{\varphi} + c_{S} (-\mathcal{S}) \dx,
\end{align}
for some $c_{\bm{v}}, c_{\varphi}$ and $c_{S}$ yet to be determined.

Using the transport theorem and the divergence theorem, we obtain the following local form
\begin{align}
\label{energyineq}
\pd_{t} e + \div (e\bm{v}) + \div \bm{J}_{e} - c_{\bm{v}} \Gamma_{\bm{v}} - c_{\varphi} \Gamma_{\varphi} + c_{S} \mathcal{S} \leq 0.
\end{align}
We now use the Lagrange multiplier method of Liu and M\"{u}ller, see for example Section 2.2 of \cite{article:AbelsGarckeGrun12} and Chapter 7 of  \cite{book:Liu}.  Let $\lambda_{\bm{v}}$, $\lambda_{\sigma}$ and $\lambda_{\varphi}$ denote Lagrange multipliers for the divergence equation (\ref{divv}), the nutrient equation (\ref{proto:nutrient}) and the order parameter equation (\ref{orderparameter:equ}).  We require that the following inequality holds for arbitrary $(\varphi, \sigma, \bm{v}, \Gamma_{\bm{v}}, \Gamma_{\varphi}, \mathcal{S}, \md \varphi, \md \sigma)$:
\begin{align}
\notag -\mathcal{D} & := \pd_{t} e + \div (e \bm{v}) + \div \bm{J}_{e} - c_{\bm{v}} \Gamma_{\bm{v}} - c_{\varphi} \Gamma_{\varphi} + c_{S} \mathcal{S} \\
\notag & - \lambda_{\bm{v}} (\div \bm{v} - \Gamma_{\bm{v}}) \\
\notag & - \lambda_{\sigma} (\md \sigma + \sigma \div \bm{v} + \div \bm{J}_{\sigma} + \mathcal{S}) \\
& - \lambda_{\varphi} (\md \varphi +  \varphi \div \bm{v} + \div \bm{J} - \Gamma_{\varphi}) \leq 0, \label{dissipation}
\end{align} 
where we used the notation
\begin{align*}
\md \varphi := \pd_{t}\varphi + \nabla \varphi \cdot \bm{v},
\end{align*}
as the material derivative of $\varphi$ with respect to $\bm{v}$.  Using the identity
\begin{align*}
\nabla \varphi \cdot \md (\nabla \varphi) = \div \left ( \md \varphi \nabla \varphi  \right ) - \md \varphi \div \left ( \nabla \varphi \right ) - \left ( \nabla \varphi \otimes \nabla \varphi \right ) : \nabla \bm{v},
\end{align*}
we compute that 
\begin{align}
\notag -\mathcal{D} &  = \div \left ( \bm{J}_{e} - \lambda_{\varphi} \bm{J}   - \lambda_{\sigma} \bm{J}_{\sigma} + B  \md \varphi \nabla \varphi \right ) \\
\notag & + \left (A \Psi'(\varphi) + \frac{\pd N}{\pd \varphi} - B \Laplace \varphi - \lambda_{\varphi} \right ) \md \varphi - \nabla \bm{v} : B \left ( \nabla \varphi \otimes \nabla \varphi \right ) \\
\notag & + \left ( \frac{\pd N}{\pd \sigma} - \lambda_{\sigma} \right ) \md \sigma +  \mathcal{S} (c_{S} - \lambda_{\sigma} ) + \nabla \lambda_{\varphi} \cdot \bm{J} + \nabla \lambda_{\sigma} \cdot \bm{J}_{\sigma} \\
& + (e - \lambda_{\varphi} \varphi - \lambda_{\sigma} \sigma - \lambda_{\bm{v}} ) \div \bm{v} + \Gamma_{\bm{v}} (\lambda_{\bm{v}} - c_{\bm{v}}) + \Gamma_{\varphi} (\lambda_{\varphi} - c_{\varphi}). \label{dissipation2}
\end{align}
We use the following notation:
\begin{align*}
N_{,\sigma}  := \frac{\pd N}{\pd \sigma}, \quad N_{,\varphi} := \frac{\pd N}{\pd \varphi}, \quad \mu :=  A\Psi'(\varphi) + N_{,\varphi} - B \Laplace \varphi.
\end{align*}
Applying the product rule to the divergence term in (\ref{dissipation2}), we then obtain
\begin{align}
\notag -\mathcal{D} & = \div \left ( \bm{J}_{e} - \lambda_{\varphi} \bm{J}   - \lambda_{\sigma} \bm{J}_{\sigma} + B  \md \varphi \nabla \varphi+ (e - \lambda_{\varphi} \varphi - \lambda_{\sigma} \sigma - \lambda_{\bm{v}} ) \bm{v} \right ) \\
\notag & + \left ( \mu - \lambda_{\varphi} \right ) \md \varphi +  \mathcal{S} (c_{S} - \lambda_{\sigma} ) + \Gamma_{\bm{v}} (\lambda_{\bm{v}} - c_{\bm{v}}) + \Gamma_{\varphi} (\lambda_{\varphi} - c_{\varphi}) + \left ( N_{,\sigma} - \lambda_{\sigma} \right ) \md \sigma  \\
& - \nabla \bm{v} : B \left ( \nabla \varphi \otimes \nabla \varphi \right ) - \bm{v} \cdot \nabla (e - \lambda_{\varphi} \varphi - \lambda_{\sigma} \sigma - \lambda_{\bm{v}} ) + \nabla \lambda_{\varphi} \cdot \bm{J} + \nabla \lambda_{\sigma} \cdot \bm{J}_{\sigma}. \label{dissipation3}
\end{align}
Employing the following identities
\begin{align*}
\nabla \bm{v} : (\nabla \varphi \otimes \nabla \varphi) & = \div ( (\nabla \varphi \otimes \nabla \varphi) \bm{v}) - \bm{v} \cdot \div (\nabla \varphi \otimes \nabla \varphi), \\
\frac{1}{2} \nabla \left ( \abs{\nabla \varphi}^{2} \right ) & = \div (\nabla \varphi \otimes \nabla \varphi) - \Laplace \varphi \nabla \varphi, \\
\pd_{t} \varphi \nabla \varphi & = \md \varphi \nabla \varphi - (\nabla \varphi \cdot \bm{v}) \nabla \varphi = \md \varphi \nabla \varphi - (\nabla \varphi \otimes \nabla \varphi) \bm{v},
\end{align*}
in (\ref{dissipation3}), we arrive at
\begin{align}
\notag -\mathcal{D} & = \div \left ( \bm{J}_{e} - \lambda_{\varphi} \bm{J}  - \lambda_{\sigma} \bm{J}_{\sigma} + B \pd_{t} \varphi \nabla \varphi  + (e - \lambda_{\varphi} \varphi - \lambda_{\sigma} \sigma - \lambda_{\bm{v}} ) \bm{v}  \right ) \\
\notag & + \left ( \mu - \lambda_{\varphi} \right ) \md \varphi +  \mathcal{S} (c_{S} - \lambda_{\sigma} ) + \Gamma_{\bm{v}} (\lambda_{\bm{v}} - c_{\bm{v}}) + \Gamma_{\varphi} (\lambda_{\varphi} - c_{\varphi}) + \left ( N_{,\sigma} - \lambda_{\sigma} \right ) \md \sigma  \\
& - \bm{v} \cdot \left ( \nabla (e - \lambda_{\varphi} \varphi - \lambda_{\sigma} \sigma - \lambda_{\bm{v}} - \tfrac{B}{2} \abs{\nabla \varphi}^{2} ) - B\Laplace \varphi \nabla \varphi \right ) + \nabla \lambda_{\varphi} \cdot \bm{J} + \nabla \lambda_{\sigma} \cdot \bm{J}_{\sigma}. \label{dissipation4}
\end{align}

\subsection{Constitutive assumptions and the general model}
We are now seeking for a model fulfilling the second law of thermodynamics in the version of a dissipation inequality stated in Section \ref{sec:energyineq}.  We do not aim for the most general model but will state certain constitutive assumptions which take the most relevant effects into account.  We hence make the following constitutive assumptions:
\begin{subequations}\label{assump:constitutive}
\begin{align}
\bm{J}_{e} & =  \lambda_{\varphi} \bm{J}   + \lambda_{\sigma} \bm{J}_{\sigma} - B  \pd_{t} \varphi \nabla \varphi - (e - \lambda_{\varphi} \varphi - \lambda_{\sigma} \sigma - \lambda_{\bm{v}} ) \bm{v} , \label{constitutive:Je} \\
c_{\mathcal{S}} & = \lambda_{\sigma} = N_{,\sigma}, \quad c_{\varphi} = \lambda_{\varphi} = \mu, \quad c_{\bm{v}} = \lambda_{\bm{v}}, \label{constitutive:c} \\
\bm{J}_{\sigma} & = -n(\varphi) \nabla N_{,\sigma}, \quad \bm{J} = -m(\varphi) \nabla \mu,\label{constitutive:fluxes} 
\end{align}
\end{subequations}
where $n(\varphi)$ and $m(\varphi)$ are non-negative mobilities.  We introduce a pressure-like function $p$ and choose
\begin{align}\label{pressurechoice}
\lambda_{\bm{v}} = p - A \Psi(\varphi) - \frac{B}{2} \abs{\nabla \varphi}^{2} + e - \mu \varphi - N_{,\sigma} \sigma, 
\end{align}
and for a positive constant $K$,
\begin{align}
\notag \bm{v} & = K \left ( \nabla (e - \mu \varphi - N_{,\sigma} \sigma - \lambda_{\bm{v}} - \tfrac{B}{2} \abs{\nabla \varphi}^{2} ) - B\Laplace \varphi \nabla \varphi \right ) \\
\notag & =   K \left ( \nabla (-p + A \Psi(\varphi)) - B \Laplace \varphi \nabla \varphi \right ) \\
& = -K (\nabla p - (\mu - N_{,\varphi}) \nabla \varphi). \label{constutitive:velo}
\end{align}
Eq. (\ref{constitutive:Je}) makes a constitutive assumption for the energy flux $\bm{J}_{e}$ which guarantees that the divergence term in (\ref{dissipation4}) vanishes.  It contains classical terms like $\mu \bm{J}$ and $N_{,\sigma} \bm{J}_{\sigma}$ which describe energy flux due to mass diffusion and the non-classical term $B \pd_{t} \varphi \nabla \varphi$ which is due to moving phase boundaries, see also \cite{article:AltPawlow92,article:AltPawlow96} where this term is discussed.  The last term in (\ref{constitutive:Je}) will result in energy changes due to work by macroscopic stress, compare \cite{article:AbelsGarckeGrun12}.  Meanwhile, (\ref{constitutive:c}), (\ref{constitutive:fluxes}), (\ref{pressurechoice}) and (\ref{constutitive:velo}) are considered in order for the right hand side of (\ref{dissipation4}) to be non-positive for arbitrary values of $(\varphi, \sigma, \bm{v}, \Gamma_{\bm{v}}, \Gamma_{\varphi}, \mathcal{S}, \md \varphi, \md \sigma)$.  We mention that (\ref{constutitive:velo}) is a Darcy law with force $(\mu - N_{,\varphi}) \nabla \varphi$.

Thus, the model equations for tumour growth are
\begin{subequations}
\label{Model:CHDarcy}
\begin{align}
\div \bm{v} & = \Gamma_{\bm{v}}, \label{CHDarcy:div} \\
\bm{v} & = -K (\nabla p - \mu \nabla \varphi + N_{,\varphi} \nabla \varphi), \label{CHDarcy:velo} \\
\pd_{t}\varphi + \div (\varphi \bm{v})  &= \div (m(\varphi) \nabla \mu) + \Gamma_{\varphi}, \label{CHDarcy:order} \\
\mu & = A\Psi'(\varphi) - B \Laplace \varphi + N_{,\varphi}, \label{CHDarcy:mu} \\
\pd_{t}\sigma + \div (\sigma \bm{v}) & = \div (n(\varphi) \nabla N_{,\sigma})  -\mathcal{S}, \label{CHDarcy:sigma}
\end{align}
\end{subequations}
where 
\begin{align*}
\Gamma_{\bm{v}} = \bar{\rho}_{1}^{-1}\Gamma_{1} + \bar{\rho}_{2}^{-1} \Gamma_{2}, \quad \Gamma_{\varphi} = \bar{\rho}_{2}^{-1}\Gamma_{2} - \bar{\rho}_{1}^{-1}\Gamma_{1}.
\end{align*}
Supplemented with the boundary conditions
\begin{align}\label{Model:CHDarcy:BC}
\nabla \varphi \cdot \bm{\nu} = \nabla \mu \cdot \bm{\nu} = 0  \text{ on } \pd \Omega,
\end{align}
then the above model satisfies the following energy equality:
\begin{align}
\notag & \; \frac{\dd}{\dt} \int_{\Omega} \left [ A \Psi(\varphi) + \frac{B}{2} \abs{\nabla \varphi}^{2} + N(\varphi, \sigma) \right ] \dx   \\
\notag + & \; \int_{\Omega} m(\varphi) \abs{\nabla \mu}^{2} + n(\varphi) \abs{\nabla N_{,\sigma}}^{2} + \frac{1}{K} \abs{\bm{v}}^{2} \dx + \int_{\Omega} \mathcal{S} N_{,\sigma} - \lambda_{\bm{v}} \Gamma_{\bm{v}} - \mu \Gamma_{\varphi}  \dx \\
+ & \; \int_{\pd \Omega} \bm{v} \cdot \bm{\nu} (N(\varphi, \sigma) + p) - n(\varphi) N_{,\sigma} \nabla N_{,\sigma} \cdot \bm{\nu} \dHaus = 0.\label{Model:CHDarcy:energyeq}
\end{align}
This follows from integrating (\ref{dissipation4}) over $\Omega$ and using the definition of $-\mathcal{D}$ from (\ref{dissipation}), the constitutive assumptions (\ref{assump:constitutive}), (\ref{pressurechoice}), and (\ref{constutitive:velo}), and applying the divergence theorem.  Here, we have not prescribed boundary conditions for $N_{,\sigma}$ and $\bm{v}$.  We will look at suitable boundary conditions for them later.

We point out that using (\ref{proto:mixturemass}), (\ref{mixturedensity}), (\ref{CHDarcy:div}), (\ref{CHDarcy:order}), and the definition of $\Gamma_{\varphi}$ and $\Gamma_{\bm{v}}$, we obtain
\begin{align*}
\Gamma_{1} + \Gamma_{2} & = \pd_{t}\rho + \div (\rho \bm{v}) + \div\bm{\mathcal{J}} \\
& = \frac{\overline{\rho}_{2} - \overline{\rho}_{1}}{2} (\md \varphi + \varphi \div \bm{v}) + \frac{\overline{\rho}_{2} + \overline{\rho}_{1}}{2} \div \bm{v} + \div \bm{\mathcal{J}} \\
& = \frac{\overline{\rho}_{2} - \overline{\rho}_{1}}{2} (\div (m(\varphi) \nabla \mu) + \overline{\rho}_{2}^{-1} \Gamma_{2} - \overline{\rho}_{1}^{-1} \Gamma_{1}) + \frac{\overline{\rho}_{2} + \overline{\rho}_{1}}{2} (\overline{\rho}_{1}^{-1} \Gamma_{1} + \overline{\rho}_{2}^{-1} \Gamma_{2}) + \div \bm{\mathcal{J}} \\
& = \div \left ( \bm{\mathcal{J}} + \tfrac{\overline{\rho}_{2} - \overline{\rho}_{1}}{2} m(\varphi) \nabla \mu \right ) + \Gamma_{1} + \Gamma_{2}.
\end{align*}
Thus, we identify $\bm{\mathcal{J}} = -\tfrac{\overline{\rho}_{2} - \overline{\rho}_{1}}{2} m(\varphi) \nabla \mu$, and the equation for $\rho$ becomes
\begin{equation}\label{PF:densityrho}
\pd_{t}\rho + \div (\rho \bm{v}) = \div \left (\tfrac{\overline{\rho}_{2} - \overline{\rho}_{1}}{2} m(\varphi) \nabla \mu \right ) + \Gamma_{1} + \Gamma_{2}.
\end{equation}

\begin{rmk}[Reformulations of the pressure and Darcy's law]
In the above derivation, we may consider the following pressure-type functions:
\begin{itemize}
\item Let $q := p - A \Psi(\varphi) - \frac{B}{2} \abs{\nabla \varphi}^{2}$ so that $\lambda_{\bm{v}} = q + e - \mu \varphi - N_{,\sigma} \sigma$ and
\begin{align}\label{Darcy:type0}
\bm{v} = K (\nabla (-q - \tfrac{B}{2} \abs{\nabla \varphi}^{2}) - B \Laplace \varphi \nabla \varphi) = - K (\nabla q + B \div (\nabla \varphi \otimes \nabla \varphi)).
\end{align}
\item Let $\hat{p} := p + N(\varphi, \sigma)$ so that $\lambda_{\bm{v}} = \hat{p} - \mu \varphi - N_{,\sigma} \sigma$ and
\begin{align}
\notag \bm{v} & = K \left (\nabla (N(\varphi, \sigma) + A \Psi(\varphi) -\hat{p}) - B \Laplace \varphi \nabla \varphi \right ) \\
& = -K (\nabla \hat{p} - \mu \nabla \varphi - N_{,\sigma} \nabla \sigma). \label{Darcy:type1}
\end{align}
\item Let $\tilde{p} := p + N(\varphi, \sigma) - \mu \varphi - N_{,\sigma} \sigma$ so that $\lambda_{\bm{v}} = \tilde{p}$ and
\begin{align}
\notag \bm{v} & = K (\nabla (N(\varphi, \sigma) + A \Psi(\varphi) - \mu \varphi - N_{,\sigma} \sigma - \tilde{p}) - B \Laplace \varphi \nabla \varphi) \\
& = -K (\nabla \tilde{p} + \varphi \nabla \mu + \sigma \nabla N_{,\sigma}). \label{Darcy:type2}
\end{align}
\end{itemize}
We point out that \eqref{Darcy:type0} can also be obtained from the momentum balance of the Navier--Stokes--Cahn--Hilliard equations 
\begin{align*}
\pd_{t} (\rho \bm{v}) + \div (\rho \bm{v} \otimes \bm{v}) - \div (\eta \, ( \nabla \bm{v} + (\nabla \bm{v})^{\top})) + \nabla q = - \div (B \nabla \varphi \otimes \nabla \varphi)
\end{align*}
by neglecting the inertia terms and replacing the viscous term with a multiple of the velocity.  This is consistent with the classical derivation of Darcy's law.  

Meanwhile, in \eqref{Darcy:type1} we have the gradient of the primary variables $(\varphi, \sigma)$ multiplied by their corresponding chemical potentials $(\mu, N_{,\sigma})$, and vice versa in \eqref{Darcy:type2} (compare with the interfacial term $\bm{K}$ in Section 3 of \cite{article:AbelsGarckeGrun12} and Eq. $\mathrm{(2.34)}$ of \cite{article:LeeLowengrubGoodman01}).  It is common to reformulate the pressure as above to obtain equations of momentum balance in the Navier--Stokes--Cahn--Hilliard equations or the Cahn--Hilliard--Darcy equations that are more amenable to further analysis.  See for instance \cite{article:AbelsDepnerGarcke13,article:FengWise12}.
\end{rmk}

\subsection{Specific models}
\subsubsection{Zero excess of total mass }
Assuming $\Gamma_{2} = - \Gamma_{1} =: \Gamma $, so that there is no source term in (\ref{PF:densityrho}), and let 
\begin{align}\label{defn:alpharhoS}
\alpha := \frac{1}{\bar{\rho}_{2}} - \frac{1}{\bar{\rho}_{1}}, \quad \bar{\rho}_{S} = \frac{1}{\bar{\rho}_{2}} + \frac{1}{\bar{\rho}_{1}},
\end{align}
so that 
\begin{align*}
\Gamma_{\bm{v}} = \alpha \Gamma, \quad \Gamma_{\varphi} = \overline{\rho}_{S} \Gamma.
\end{align*}
Then (\ref{Model:CHDarcy}) becomes
\begin{subequations}
\label{Model:CHDarcy:Gamma}
\begin{align}
\div \bm{v} & = \alpha \Gamma, \\
\bm{v} & = -K (\nabla p - \mu \nabla \varphi + N_{,\varphi} \nabla \varphi), \\
\pd_{t}\varphi + \div(\bm{v} \varphi)  &= \div (m(\varphi) \nabla \mu) + \bar{\rho}_{S} \Gamma, \\
\mu & = A\Psi'(\varphi) - B \Laplace \varphi + N_{,\varphi}, \\
\pd_{t}\sigma + \div (\sigma \bm{v}) & = \div (n(\varphi) \nabla N_{,\sigma})  -\mathcal{S} .
\end{align}
\end{subequations}
In the case that the densities are equal, i.e., $\bar{\rho}_{1} = \bar{\rho}_{2} = \bar{\rho}$, then, $\alpha = 0$ and $\bar{\rho}_{S} = \frac{2}{\bar{\rho}}$, and (\ref{Model:CHDarcy:Gamma}) becomes
\begin{subequations}
\label{Model:CHDarcy:EqualDensities}
\begin{align}
\div \bm{v} & = 0, \\
\bm{v} & = -K (\nabla p - \mu \nabla \varphi + N_{,\varphi} \nabla \varphi), \\
\pd_{t}\varphi + \bm{v} \cdot \nabla \varphi  &= \div (m(\varphi) \nabla \mu) 
+ \tfrac{2}{\bar{\rho}} \Gamma, \\
\mu & = A\Psi'(\varphi) - B \Laplace \varphi + N_{,\varphi}, \\
\pd_{t}\sigma + \bm{v} \cdot \nabla \sigma & = \div (n(\varphi) \nabla N_{,\sigma})  -\mathcal{S} .
\end{align}
\end{subequations}

\subsubsection{Absence of nutrients}
Setting $\sigma = N(\sigma, \varphi) = 0$, then (\ref{Model:CHDarcy}) simplifies to
\begin{subequations}
\label{Model:CHDarcy:NoNurtient}
\begin{align}
\div \bm{v} & = \bar{\rho}_{1}^{-1}\Gamma_{1} + \bar{\rho}_{2}^{-1} \Gamma_{2}, \\
\bm{v} & = -K (\nabla p - \mu \nabla \varphi), \\
\pd_{t}\varphi + \div (\bm{v} \varphi) &= \div (m(\varphi) \nabla \mu) + \bar{\rho}_{2}^{-1}\Gamma_{2} - \bar{\rho}_{1}^{-1}\Gamma_{1}, \\
\mu & = A\Psi'(\varphi) - B \Laplace \varphi.
\end{align}
\end{subequations}

\subsubsection{Zero velocity, zero excess of total mass and equal densities}
Suppose the volume-averaged mixture velocity $\bm{v}$ is zero, the excess of total mass $\Gamma_{1} + \Gamma_{2}$ is zero and the densities are equal.  Then, substituting $\bm{v} = \bm{0}$ in (\ref{Model:CHDarcy:EqualDensities}) and neglecting the Darcy system (\ref{Model:CHDarcy:EqualDensities}a,b), we obtain
\begin{subequations}
\label{Model:CHDarcy:NoVelo:EqualDensities}
\begin{align}
\pd_{t}\varphi  &= \div (m(\varphi) \nabla \mu) + \tfrac{2}{\bar{\rho}} \Gamma, \\
\mu & = A\Psi'(\varphi) - B \Laplace \varphi + N_{,\varphi}, \\
\pd_{t}\sigma  & = \div (n(\varphi) \nabla N_{,\sigma}) -\mathcal{S} .
\end{align}
\end{subequations}

\subsubsection{Boundary conditions for velocity and nutrient}
For the nutrient, we may prescribe a Robin type boundary condition:
\begin{align}\label{bc:nutrient:Robin}
 (n(\varphi) \nabla N_{,\sigma}) \cdot \bm{\nu} = c (\sigma_{\infty} - \sigma) \text{ on } \pd \Omega,
\end{align}
where $c \geq 0$ is a constant, and $\sigma_{\infty}$ denotes a given supply at the boundary.  When $c = 0$, we obtain the zero flux boundary condition:
\begin{align}\label{bc:nutrient:Neumann0}
( n(\varphi) \nabla N_{,\sigma}) \cdot \bm{\nu} = 0 \text{ on } \pd \Omega.
\end{align}
If we formally send $c \to \infty$, then we obtain the Dirichlet boundary condition:
\begin{align}\label{bc:nutrient:Dirichlet}
\sigma = \sigma_{\infty} \text{ on } \pd \Omega.
\end{align}
We may consider a boundary condition for the normal component of the velocity (which corresponds to a Neumann boundary condition for the pressure):
\begin{align}\label{bc:normalvelo}
-\bm{v} \cdot \bm{\nu} = K \nabla p \cdot \bm{\nu} = g_{2} \text{ on } \pd \Omega,
\end{align}
for some given function $g_{2}$.  We point out that a compatibility condition is required to hold if we consider the boundary condition (\ref{bc:normalvelo}) for the Models (\ref{Model:CHDarcy}), (\ref{Model:CHDarcy:Gamma}), (\ref{Model:CHDarcy:EqualDensities}), and (\ref{Model:CHDarcy:NoNurtient}).  Namely, if the mass exchange terms $\Gamma_{1}$ and $\Gamma_{2}$ are given, then we require that $g_{2}$ satisfies 
\begin{align*}
- \int_{\pd \Omega} g_{2} \dHaus & = \int_{\pd \Omega} \bm{v} \cdot \bm{\nu} \dHaus = \int_{\Omega} \div \bm{v} \dx \\
& = \begin{cases} \displaystyle \int_{\Omega} \overline{\rho}_{1}^{-1} \Gamma_{1} + \overline{\rho}_{2}^{-1} \Gamma_{2} \dx & \text{ for Models } (\ref{Model:CHDarcy}), (\ref{Model:CHDarcy:NoNurtient}), \\[2mm]
\displaystyle \int_{\Omega} \alpha \Gamma \dx & \text{ for Model } (\ref{Model:CHDarcy:Gamma}), \\[2mm]
0 & \text{ for Model } (\ref{Model:CHDarcy:EqualDensities}).
\end{cases} 
\end{align*}
However, if the mass source terms $\Gamma_{i}$ depend on $\varphi$, $\sigma$ or $\mu$, then considering (\ref{bc:normalvelo}) as a boundary condition would imply that $\varphi$, $\mu$ and $\sigma$ have to satisfy
\begin{align*}
\int_{\Omega} \overline{\rho}_{1}^{-1} \Gamma_{1}(\varphi, \sigma, \mu) + \overline{\rho}_{2}^{-1} \Gamma_{2} (\varphi, \sigma, \mu) \dx = \int_{\pd \Omega} - g_{2} \dHaus.
\end{align*} 
Alternatively, we can prescribe a boundary condition for the pressure.  Recall the reformulated pressure $\hat{p}$ and the Darcy's law (\ref{Darcy:type1}).  We can prescribe a Dirichlet boundary condition:
\begin{align}\label{bc:pressure:Dirichlet}
\hat{p} = g_{1} \text{ on } \pd \Omega,
\end{align}
for some given function $g_{1}$.  We may also consider the mixed boundary condition as in Section 2.3.3 of \cite{book:ChenHuanMa} (which corresponds to a Robin boundary condition for the pressure):
\begin{align}\label{bc:mixed}
a \hat{p} - b \bm{v} \cdot \bm{\nu} = a\hat{p} + b K \nabla \hat{p} \cdot \bm{\nu} - b K N_{,\sigma} \nabla \sigma \cdot \bm{\nu} =  g_{3} \text{ on } \pd \Omega,
\end{align}
for constants $a,b \geq 0$ and a given function $g_{3}$.

\subsection{Comparison to other models in the literature}
\subsubsection{Absence of nutrients}
Scaling mass, permeability, and mobility appropriately, by setting
\begin{align*}
\Gamma_{1} = 0, \quad \Gamma := \Gamma_{2}, \quad \bar{\rho}_{2} = \bar{\rho}_{1} = 1, \quad K = 1, \quad m(\varphi) = 1
\end{align*}
in (\ref{Model:CHDarcy:NoNurtient}), we obtain the following system
\begin{subequations}
\label{Model:LowengrubTitiZhao}
\begin{align}
\div \bm{v} & = \Gamma, \\
\bm{v} &= -\nabla p + \mu \nabla \varphi, \\
\pd_{t} \varphi + \div (\bm{v} \varphi) & = \Laplace \mu + \Gamma, \\
\mu & = A\Psi'(\varphi) - B \Laplace \varphi.
\end{align}
\end{subequations}
The existence of strong solutions in 2D and 3D have been studied in \cite{article:LowengrubTitiZhao13} for the case $\Gamma = 0$.  For the case where $\Gamma \neq 0$ is prescribed, existence of global weak solutions and unique local strong solutions in both 2D and 3D can be found in \cite{preprint:JiangWuZheng14}.  We also refer the reader to \cite{preprint:BosiaContiGrasselli14} for the study of weak solutions to a related system, denoted as the Cahn--Hilliard--Brinkman system, where an additional viscosity term $\eta \div D(\bm{v})$ is added to the left hand side of the velocity equation (\ref{Model:LowengrubTitiZhao}b) and the mass exchange $\Gamma$ is set to zero.  Here, $D(\bm{v}) = \tfrac{1}{2} (\nabla \bm{v} + (\nabla \bm{v})^{\top})$ is the rate of deformation tensor and $\eta$ is the viscosity.

\subsubsection{Zero velocity, zero excess of total mass and equal densities}\label{sec:0velo0massexcesEqualDens}
We consider the model (\ref{Model:CHDarcy:NoVelo:EqualDensities}) with the rescaled density $\bar{\rho} = 1$.  Let $\mathcal{P}$, $\mathcal{A}$, $\mathcal{C}$, $\chi_{\sigma}$, $\chi_{\varphi}$ be non-negative constants.  For physically relevant values of the model variables, i.e., $\varphi \in [-1,1]$ and $\sigma \geq 0$, we choose
\begin{subequations}\label{specificchoice}
\begin{align}
\Gamma & = (\mathcal{P} \sigma - \mathcal{A}) h(\varphi), \label{specificchoice:Gamma} \\
 N(\varphi, \sigma) & = \frac{\chi_{\sigma}}{2} \abs{\sigma}^{2} + \chi_{\varphi} \sigma(1-\varphi), \label{specificchoice:N} \\
 \mathcal{S} & = \mathcal{C} \sigma h(\varphi), \label{specificchoice:S}
 \end{align}
\end{subequations}
where $h(\varphi)$ is an interpolation function with $h(-1) = 0$ and $h(1) = 1$.

We have elaborated on the physical motivations for the particular forms of $\Gamma$ and $\mathcal{S}$ in Section \ref{sec:Intro}.  For the choice of $N(\varphi, \sigma)$, if both $\chi_{\varphi}$ and $\chi_{\sigma}$ are positive constants, then for physically relevant parameter values, i.e., $\sigma \geq 0$, and $\varphi \in [-1,1]$,
\begin{align}\label{N,sigma}
N_{,\sigma} = \chi_{\sigma} \sigma + \chi_{\varphi} (1-\varphi) \geq 0.
\end{align} 
Thus, this choice of the flux $\nabla N_{,\sigma}$  provides two transport mechanisms for the nutrient $\sigma$.  The first term $\chi_{\sigma} \nabla\sigma$ results in a diffusion process along negative gradients of $\sigma$, while the second term $-\chi_{\varphi} \nabla \varphi$ is a chemotactic term that drives the nutrient towards the tumour cell regions.  In particular, in the tumour cell regions $\{\varphi = +1\}$, the nutrient will only experience diffusion, while in the healthy cell regions $\{ \varphi = -1\}$, the nutrient will experience diffusion and active transport to the tumour.

We point out that for this particular form of $N_{,\sigma}$, together with the zero Neumann boundary condition for $\varphi$, we have
\begin{align*}
\nabla N_{,\sigma} \cdot \bm{\nu} = \chi_{\sigma} \nabla \sigma \cdot \bm{\nu} - \chi_{\varphi} \nabla \varphi \cdot \bm{\nu} = \chi_{\sigma} \nabla \sigma \cdot \bm{\nu} \text{ on } \pd \Omega.
\end{align*}
With these choices,  (\ref{Model:CHDarcy:NoVelo:EqualDensities}) becomes
\begin{subequations}\label{Model:Specific:NoVelo:1}
\begin{align}
\pd_{t}\varphi  &= \div (m(\varphi) \nabla \mu) + 2( \mathcal{P} \sigma - \mathcal{A})  h(\varphi), \label{Model:Specific:NoVelo:1:a} \\
\mu & = A\Psi'(\varphi) - B \Laplace \varphi - \chi_{\varphi} \sigma, \label{Model:Specific:NoVelo:1:b} \\
\pd_{t}\sigma  & = \div (n(\varphi) (\chi_{\sigma} \nabla \sigma - \chi_{\varphi} \nabla \varphi)) - \mathcal{C} \sigma h(\varphi). \label{Model:Specific:NoVelo:1:c}
\end{align}
\end{subequations}
We remark that (\ref{Model:Specific:NoVelo:1}) is similar to Eq. (68)-(73) of \cite{article:CristiniLiLowengrubWise09}, the two-phase diffuse interface tumour model (Eq. 5.27) of \cite{article:OdenHawkinsPrudhomme10}, and model $\mathcal{M}_{2}$ of \cite{article:HawkinsPrudhommevanderZeeOden13}.  The only difference between these three models and (\ref{Model:Specific:NoVelo:1}) is that the flux for the nutrient equation (\ref{Model:Specific:NoVelo:1:c}) consists of an advection term and a Fickian diffusion term for \cite{article:CristiniLiLowengrubWise09}, while in \cite{article:OdenHawkinsPrudhomme10,article:HawkinsPrudhommevanderZeeOden13}, the nutrient is in a quasi-steady state and the flux for the nutrient equation is a Fickian diffusive flux.  We point out that in \cite{article:CristiniLiLowengrubWise09,article:HawkinsPrudhommevanderZeeOden13,article:OdenHawkinsPrudhomme10}, $h(\varphi)$ is replaced by $\varphi$ in the definition of $\Gamma$ and $\mathcal{S}$.  Since, in their notation, $\varphi \in [0,1]$ denotes the tumour volume fraction instead of the difference of volume fractions.

Next, choosing $N(\varphi, \sigma)$ as in (\ref{specificchoice:N}) above, and
\begin{align*}
 \Gamma = \frac{1}{2}P(\varphi)(N_{,\sigma} - \mu), \quad \mathcal{S} = P(\varphi)(N_{,\sigma} - \mu), 
\end{align*}
where $P(\varphi)$ is a non-negative function, then (\ref{Model:CHDarcy:NoVelo:EqualDensities}) becomes
\begin{subequations}\label{Model:Specific:NoVelo:2}
\begin{align}
\pd_{t}\varphi  &= \div (m(\varphi) \nabla \mu) + P(\varphi)(\chi_{\sigma} \sigma + \chi_{\varphi}(1- \varphi)  - \mu), \label{Model:Specific:NoVelo:2:a} \\
\mu & = A\Psi'(\varphi) - B \Laplace \varphi - \chi_{\varphi} \sigma, \label{Model:Specific:NoVelo:2:b} \\
\pd_{t}\sigma  & = \div (n(\varphi) (\chi_{\sigma} \nabla \sigma - \chi_{\varphi} \nabla \varphi)) -P(\varphi)(\chi_{\sigma} \sigma + \chi_{\varphi}(1- \varphi) - \mu) . \label{Model:Specific:NoVelo:2:c}
\end{align}
\end{subequations}
This is similar to the model derived in \cite{article:HawkinsZeeOden12}, where the chemical potentials $N_{,\sigma}$ and $\mu$ enter as source terms in (\ref{Model:Specific:NoVelo:2:a}) and (\ref{Model:Specific:NoVelo:2:c}).  The specific form for $\Gamma$ is motivated by linear phenomenological constitutive laws for chemical reactions.  The non-negative function $P(\varphi)$ takes on the form
\begin{align}\label{HawkinsPvarphi}
P(\varphi) = \begin{cases} \delta P_{0} (1+\varphi), & \text{ if } \varphi \geq -1, \\
0, & \text{ otherwise},
\end{cases}
\end{align}
for positive constants $\delta$ and $P_{0}$.  Subsequently, if we choose 
\begin{align*}
\chi_{\sigma} = 1, \quad \chi_{\varphi} = 0, \quad n(\varphi) = m(\varphi) = 1
\end{align*} 
in (\ref{Model:Specific:NoVelo:2}), we obtain
\begin{subequations}\label{Model:Specific:NoVelo:3}
\begin{align}
\pd_{t}\varphi  &= \Laplace \mu + P(\varphi)(\sigma - \mu), \label{Model:Specific:NoVelo:3:a} \\
\mu & = A\Psi'(\varphi) - B \Laplace \varphi, \label{Model:Specific:NoVelo:3:b} \\
\pd_{t}\sigma  & = \Laplace \sigma -P(\varphi)(\sigma - \mu) . \label{Model:Specific:NoVelo:3:c}
\end{align}
\end{subequations}
This is the model studied in \cite{article:FrigeriGrasselliRocca15}, for a more general function $P(\varphi)$ than (\ref{HawkinsPvarphi}), while a viscosity regularised version of (\ref{Model:Specific:NoVelo:3}) (where there is an extra $\alpha \pd_{t}\mu$ term on the left hand side of (\ref{Model:Specific:NoVelo:3:a}) and an extra $\alpha \pd_{t}\varphi$ term on the right hand side of (\ref{Model:Specific:NoVelo:3:b}) for a positive constant $\alpha$) is studied in \cite{article:ColliGilardiHilhorst15}.  A formal asymptotic limit for the viscosity regularised version of (\ref{Model:Specific:NoVelo:3}) is derived in \cite{article:HilhorstKaampmannNguyenZee15}.

\section{Sharp Interface Asymptotics}\label{sec:Asym}
We consider Model (\ref{Model:CHDarcy}) with the following choices and assumptions:
\begin{assump}
\
\begin{itemize}
\item $A = \frac{\beta}{\eps}$ and $B = \beta \eps$ for positive constants $\beta, \eps > 0$.
\item $N(\varphi, \sigma)$ is chosen as in $(\ref{specificchoice:N})$ with constant positive parameters $\chi_{\sigma}, \chi_{\varphi} > 0$.
\item The mass exchange terms $\Gamma_{i}$, $i = 1,2$, and the nutrient consumption term $\mathcal{S}$ depend only on $\sigma$, $\mu$, and $\varphi$, and not on any derivatives.
\item The mobilities $m(\varphi)$ and $n(\varphi)$ are strictly positive and continuously differentiable.
\item The potential $\Psi$ is chosen to be either the smooth double-well potential $\Psi(\varphi) = \tfrac{1}{4} (1-\varphi^{2})^{2}$ or the double-obstacle potential 
\begin{align}\label{defn:DoubleObstacle}
\Psi(\varphi) := \frac{1}{2}(1-\varphi^{2}) + I_{[-1,1]}(\varphi), \quad I_{[-1,1]}(\varphi) = \begin{cases}
0 & \text{ if } \abs{\varphi} \leq 1, \\
+\infty & \text{ otherwise }.
\end{cases}
\end{align}
\end{itemize}
\end{assump}
With these choices, Model (\ref{Model:CHDarcy}) becomes
\begin{subequations}\label{Model:CHDarcy:Asymp}
\begin{align}
\div \bm{v} & = \overline{\rho}_{1}^{-1} \Gamma_{1}(\sigma, \varphi, \mu) + \overline{\rho}_{2}^{-1} \Gamma_{2}(\sigma, \varphi, \mu), \label{Asymp:div} \\
\bm{v} & = -K (\nabla p - \mu \nabla \varphi - \chi_{\varphi} \sigma \nabla \varphi), \label{Asymp:velo} \\
\pd_{t}\varphi + \div ( \bm{v} \varphi)  &= \div (m(\varphi) \nabla \mu) + \overline{\rho}_{2}^{-1}\Gamma_{2}(\sigma, \varphi, \mu) - \overline{\rho}_{1}^{-1} \Gamma_{1}(\sigma, \varphi, \mu), \label{Asymp:order} \\
\mu & = \frac{\beta}{\eps}\Psi'(\varphi) - \beta \eps \Laplace \varphi - \chi_{\varphi} \sigma, \label{Asym:mu} \\
\pd_{t}\sigma + \div (\sigma \bm{v}) & = \div (n(\varphi) ( \chi_{\sigma} \nabla \sigma - \chi_{\varphi} \nabla \varphi)) -\mathcal{S}(\sigma, \varphi, \mu). \label{Asymp:sigma}
\end{align}
\end{subequations}
We point out that in the case of the double-obstacle potential, the ``derivative'' $\Psi'$ is to be understood in the sense of subdifferentials, i.e.,
\begin{align}
\label{defn:DoubleObstacle:subdiff}
\Psi'(\varphi) = -\varphi + \pd I_{[-1,1]}(\varphi), \quad \pd I_{[-1,1]}(\varphi) = \begin{cases}
(-\infty,0] & \text{ if } \varphi = -1, \\
0 & \text{ if } \abs{\varphi} < 1,  \\
[0, + \infty) & \text{ if } \varphi = +1,
\end{cases}
\end{align}
and (\ref{Asym:mu}) will have to be formulated in terms of the following variational inequality:
\begin{align}
\label{Asym:mu:variationalinequ}
\int_{\Omega} - \mu (\psi - \varphi) - \frac{\beta}{\eps} \varphi(\psi - \varphi) + \beta \eps \nabla \varphi \cdot \nabla (\psi - \varphi) - \chi_{\varphi} \sigma (\psi - \varphi) \dx \geq 0,
\end{align}
for all $\psi \in \mathcal{K} := \{ \eta \in H^{1}(\Omega) : \abs{\eta} \leq 1 \}$.

We perform a formal asymptotic analysis on Model (\ref{Model:CHDarcy:Asymp}) in the limit $\eps \to 0$.  Details of the method can also be found in \cite{article:AbelsGarckeGrun12,article:BhateBowerKumar02,incoll:BloweyElliott93,article:GarckeLamStinner14}.  The following assumptions are considered:
\begin{assump}\label{assump:SharpInterfaceAsym}
\
\begin{itemize}
\item We assume that for small $\eps$, the domain $\Omega$ can be divided into two open subdomains $\Omega^{\pm}(\eps)$, separated by an interface $\Sigma(\eps)$ that does not intersect with $\pd \Omega$.  
\item We assume that there is a family $(\varphi_{\eps}, \bm{v}_{\eps}, p_{\eps}, \mu_{\eps}, \sigma_{\eps})_{\eps > 0}$ of solutions to $(\ref{Model:CHDarcy:Asymp})$, which are sufficiently smooth and have an asymptotic expansion in $\eps$ in the bulk regions away from $\Sigma(\eps)$ (the outer expansion), and another expansion in the interfacial region close to $\Sigma(\eps)$ (the inner expansion). 
\item We assume that the zero level sets of $\varphi_{\eps}$ converge to a limiting hypersurface $\Sigma_{0}$ moving with normal velocity $\mathcal{V}$.
\end{itemize}
\end{assump}
The idea of the method is to plug the outer and inner expansions in the model equations and solve them order by order, in addition we have to define a suitable region where these expansions should match up.

We will use the following notation: $(\ref{Asymp:sigma})_{O}^{\alpha}$ and $(\ref{Asymp:sigma})_{I}^{\alpha}$  denote the terms resulting from the order $\alpha$ outer and inner expansions of (\ref{Asymp:sigma}), respectively.

\subsection{Outer expansion}\label{sec:OuterExp}
We assume that for $f_{\eps} \in \{\varphi_{\eps}, \bm{v}_{\eps}, p_{\eps}, \mu_{\eps}, \sigma_{\eps} \}$, the following outer expansions hold:
\begin{align*}
f_{\eps} = f_{0} + \eps f_{1} + \eps^{2} f_{2} + \ldots.
\end{align*}
To leading order $(\ref{Asym:mu})_{O}^{-1}$ gives
\begin{align}
\label{Outer:mu:-1}
-\beta \Psi'(\varphi_{0}) = 0.
\end{align}
The solutions to (\ref{Outer:mu:-1}) corresponding to minima of $\Psi$ are $\varphi_{0} = \pm 1$, and thus, we can define the tumour tissues and the healthy tissues region by
\begin{align}\label{defn:OmegaTOmegaH}
\Omega_{T} := \{ x \in \Omega : \varphi_{0}(x) = 1 \}, \quad \Omega_{H} := \{ x \in \Omega : \varphi_{0}(x) = - 1\}.
\end{align}
Then, thanks to $\nabla \varphi_{0} = \bm{0}$, we obtain from the equations to zeroth order:
\begin{align}
\div \bm{v}_{0} & = \overline{\rho}_{1}^{-1}\Gamma_{1}(\sigma_{0}, \varphi_{0}, \mu_{0}) + \overline{\rho}_{2}^{-1} \Gamma_{2}(\sigma_{0}, \varphi_{0}, \mu_{0}), \label{Outer:div:0} \\
\bm{v}_{0} & = -K \nabla p_{0}, \label{Outer:velo:0} \\
-\div (m(\varphi_{0}) \nabla \mu_{0}) & = \overline{\rho}_{2}^{-1}(1-\varphi_{0})\Gamma_{2}(\sigma_{0}, \varphi_{0}, \mu_{0}) - \overline{\rho}_{1}^{-1}(1+\varphi_{0}) \Gamma_{1}(\sigma_{0}, \varphi_{0}, \mu_{0}), \label{Outer:order:0} \\
\pd_{t} \sigma_{0} + \div (\sigma_{0} \bm{v}_{0}) & = \div( n(\varphi_{0}) \chi_{\sigma} \nabla \sigma_{0}) - \mathcal{S}(\sigma_{0}, \varphi_{0}, \mu_{0}). \label{Outer:sigma:0}
\end{align}
For the double-obstacle potential, we obtain from $(\ref{Asym:mu:variationalinequ})_{O}^{-1}$,
\begin{align*}
\int_{\Omega} -\beta \varphi_{0}(\psi_{0} - \varphi_{0}) \dx \geq 0 \text{ for all } \psi_{0} \in \mathcal{K} := \{ \eta \in H^{1}(\Omega) : \abs{\eta} \leq 1 \}.
\end{align*}
For this to hold for all $\abs{\psi_{0}} \leq 1$, we require that $\varphi_{0} = \pm 1$, and thus we can define $\Omega_{T}$ and $\Omega_{H}$ as before, and also recover (\ref{Outer:div:0}), (\ref{Outer:velo:0}), (\ref{Outer:order:0}), and (\ref{Outer:sigma:0}).

\subsection{Inner expansions and matching conditions}\label{sec:InnerExp}
By assumption, $\Sigma_{0}$ is the limiting hypersurface of the zero level sets of $\varphi_{\eps}$.  In order to study the limiting behaviour in these parts of $\Omega$ we introduce a new coordinate system.  

We introduce the signed distance function $d(x)$ to $\Sigma_{0}$, and set $z = \frac{d}{\eps}$ as the rescaled distance variable, and use the convention that $d(x) < 0$ in $\Omega_{H}$, and $d(x) > 0$ in $\Omega_{T}$.  Thus, the gradient $\nabla d$ points from $\Omega_{H}$ to $\Omega_{T}$, and we may use $\nabla d$ on $\Sigma_{0}$ to denote the unit normal of $\Sigma_{0}$, pointing from $\Omega_{H}$ to $\Omega_{T}$.

Let $\para(t,s)$ denote a parametrization of $\Sigma_{0}$ by arc-length $s$, and let $\bm{\nu}$ denote the unit normal of $\Sigma_{0}$, pointing into the tumour region.  Then, in a tubular neighbourhood of $\Sigma_{0}$, for sufficiently smooth function $f(x)$, we have
\begin{align*}
f(x) = f(\para(t,s) + \eps z \bm{\nu}(\para(t,s))) =: F(t,s,z).
\end{align*}
In this new $(t,s,z)$-coordinate system, the following change of variables apply, compare \cite{article:GarckeStinner06}:
\begin{align*}
\pd_{t} f & = -\frac{1}{\eps} \mathcal{V} \pd_{z} F + \text{ h.o.t.}, \\
\nabla_{x} f & = \frac{1}{\eps}\pd_{z} F \bm{\nu} + \nabla_{\Sigma_{0}} F + \text{ h.o.t.},
\end{align*}
where $\mathcal{V}$ is the normal velocity of $\Sigma_{0}$, $\nabla_{\Sigma_{0}}g$ denotes the surface gradient of $g$ on $\Sigma_{0}$ and h.o.t. denotes higher order terms with respect to $\eps$.  In particular, we have
\begin{align*}
\Laplace f  = \div_{x} (\nabla_{x} f) & = \frac{1}{\eps^{2}} \pd_{zz}F + \frac{1}{\eps}\underbrace{\div_{\Sigma_{0}} (\pd_{z}F \bm{\nu})}_{= - \kappa \pd_{z}F} + \text{ h.o.t.},
\end{align*}
where $\kappa = - \div_{\Sigma_{0}} \bm{\nu}$ is the mean curvature of $\Sigma_{0}$.  Moreover, if $\bm{v}$ is a vector-valued function with $\bm{V}(t,s,z) = \bm{v}(x)$ for $x$ in a tubular neighbourhood of $\Sigma_{0}$, then we obtain
\begin{align*}
\div_{x} \bm{v} = \frac{1}{\eps}\pd_{z} \bm{V} \cdot \bm{\nu} + \div_{\Sigma_{0}} \bm{V} + \text{ h.o.t.}.
\end{align*}
We denote the variables $\varphi_{\eps}$, $\bm{v}_{\eps}$, $p_{\eps}$, $\mu_{\eps}$, $\sigma_{\eps}$ in the new coordinate system by $\Phi_{\eps}$, $\bm{V}_{\eps}$, $P_{\eps}$, $\Xi_{\eps}$, $C_{\eps}$, respectively.  We further assume that they have the following inner expansions:
\begin{align*}
F_{\eps}(s,z) = F_{0}(s,z) + \eps F_{1}(s,z) + \dots,
\end{align*}
for $F_{\eps} \in \{ \Phi_{\eps}, \bm{V}_{\eps}, P_{\eps}, \Xi_{\eps}, C_{\eps} \}$.  The assumption that the zero level sets of $\varphi_{\eps}$ converge to $\Sigma_{0}$ implies that
\begin{align}\label{CondPhi}
\Phi_{0}(t,s, z=0) = 0.
\end{align}
Furthermore, we make the following assumption:
\begin{assump}\label{assump:DoubleObstacleAsym} 
For the double-obstacle potential, we assume that the inner variable $\Phi_{\eps}$ is monotone increasing with $z$ and the interfacial layer has finite thickness of $2l$, where the value of $l$ will be specified later.  For the double-well potential, we take $l = \infty$.  Furthermore, we assume that
\begin{align}\label{DO:Phi:assumption}
\Phi_{\eps}(t,s,z = +l) = +1, \quad \Phi_{\eps}(t,s,z = -l) = -1.
\end{align}
\end{assump}
In order to match the inner expansions valid in the interfacial region to the outer expansions of Section \ref{sec:OuterExp} we employ the matching conditions, see \cite{article:GarckeStinner06}:
\begin{align}
\label{MatchingCond1}
\lim_{z \to \pm l} F_{0}(t,s,z) &= f_{0}^{\pm}(t,x), \\
\label{MatchingCond2}
\lim_{z \to \pm l} \pd_{z}F_{0}(t,s,z) &= 0 ,\\
\label{MatchingCond3}
\lim_{z \to \pm l} \pd_{z} F_{1}(t,s,z) &= \nabla f_{0}^{\pm}(t,x) \cdot \bm{\nu},
\end{align}
where $f_{0}^{\pm}(t,x):= \lim_{\delta \searrow 0} f_{0}(t,x \pm \delta \bm{\nu})$ for $x \in \Sigma_{0}$.  Moreover, we use the following notation:  Let $\delta > 0$ and for $x \in \Sigma_{0}$ with $x - \delta \bm{\nu} \in \Omega_{H}$ and $x + \delta \bm{\nu} \in \Omega_{T}$, we denote the jump of a quantity $f$ across the interface by
\begin{align}\label{defn:jump}
\jump{f} := \lim_{\delta \searrow 0} f(t,x + \delta \bm{\nu}) - \lim_{\delta \searrow 0} f(t,x - \delta \bm{\nu}).
\end{align}
For convenience, we define the constant $\gamma > 0$ to be
\begin{align}
\label{defn:gammaconst}
\gamma := \begin{cases} \displaystyle
\int_{-\infty}^{\infty} \frac{1}{2} \sech^{4}(z/\sqrt{2}) \dz = \frac{2\sqrt{2}}{3} & \text{ for the double-well potential}, \\[2ex]
\displaystyle \int_{-\frac{\pi}{2}}^{\frac{\pi}{2}} \cos^{2}(z) \dz = \frac{\pi}{2} & \text{ for the double-obstacle potential}.
\end{cases}
\end{align}
\subsubsection{Expansions to leading order}
To leading order $(\ref{Asym:mu})_{I}^{-1}$ gives
\begin{align}
\label{Inner:mu:-1}
\pd_{zz} \Phi_{0} - \Psi'(\Phi_{0}) = 0.
\end{align}
Using (\ref{CondPhi}) we obtain that $\Phi_{0}$ can be chosen to be independent of $s$ and $t$, i.e., $\Phi_{0}$ is only a function of $z$, and solves
\begin{align}
\label{ODE}
\Phi_{0}''(z) - \Psi'(\Phi_{0}(z)) = 0, \quad \Phi_{0}(0) = 0, \quad \Phi_{0} (\pm l) = \pm 1.
\end{align}
For the double-well potential, we have the unique solution
\begin{align}
\label{Phi0profile}
\Phi_{0}(z) = \tanh \left ( \frac{z}{\sqrt{2}} \right ).
\end{align}
Furthermore, multiplying (\ref{ODE}) by $\Phi_{0}'(z)$, integrating and applying the matching conditions (\ref{MatchingCond1}) and (\ref{MatchingCond2}) to $\Phi_{0}$ gives the so-called equipartition of energy:
\begin{align}
\label{equipartition}
\frac{1}{2} \abs{\Phi_{0}'(z)}^{2} = \Psi(\Phi_{0}(z)) \text{ for all } \abs{z} < \infty.
\end{align}
Similarly, for the double-obstacle potential, we obtain from $(\ref{Asym:mu:variationalinequ})_{I}^{-1}$,
\begin{align}\label{Inner:DO:mu:-1}
\int_{\Omega} -\beta (\Phi_{0} + \pd_{zz} \Phi_{0})(\psi - \Phi_{0}) \dx \geq 0 \text{ for all } \abs{\psi} \leq 1.
\end{align}
For (\ref{Inner:DO:mu:-1}) to be satisfied, it suffices to consider $\Phi_{0}$ as a function only in $z$ which solves
\begin{align}\label{ODE:DO}
\Phi_{0}(z) + \Phi_{0}''(z) = 0, \quad \Phi_{0}(0) = 0, \quad \Phi_{0}(\pm l) = \pm 1.
\end{align}
A solution to (\ref{ODE:DO}) is
\begin{align}\label{DO:profile}
\Phi_{0}(z) = \begin{cases}
+1 & \text { if } z \geq \frac{\pi}{2}, \\
\sin(z) & \text{ if } \abs{z} \leq \frac{\pi}{2}, \\
-1 & \text{ if } z \leq -\frac{\pi}{2},
\end{cases}
\end{align}
so that $l = \frac{\pi}{2}$ for the double-obstacle potential, and we deduce from (\ref{DO:Phi:assumption}) that for the double-obstacle potential,
\begin{align}\label{DO:Phi1}
\Phi_{1}(t,s, \pm \tfrac{\pi}{2} ) = 0.
\end{align}
Moreover, we obtain the equipartition of energy (\ref{equipartition}) via a similar argument to the double-well potential.  Thanks to the equipartition of energy (\ref{equipartition}), and the definition of $\gamma$ (\ref{defn:gammaconst}), we point out that
\begin{align}\label{equipartition:gamma}
\int_{-l}^{l} \abs{\Phi_{0}'(z)}^{2} \dz = \int_{-l}^{l} 2 \Psi(\Phi_{0}(z)) \dz = \gamma.
\end{align}
For the rest of this section, we do not differentiate between the two cases of potentials, and use the notation that $l = \frac{\pi}{2}$ represents the case of the double-obstacle potential and $l = \infty$ represents the case of the double-well potential.  

Next, $(\ref{Asymp:div})_{I}^{-1}$ gives
\begin{align}
\label{Inner:div:-1}
\pd_{z} \bm{V}_{0} \cdot \bm{\nu} = 0.
\end{align}
Integrating from $-l$ to $l$ with respect to $z$, and applying the matching condition (\ref{MatchingCond1}) to $\bm{V}_{0}$ yields
\begin{align}
\label{velonormalcomp}
\jump{\bm{v}_{0}} \cdot \bm{\nu} := \bm{v}_{0}^{+} \cdot \bm{\nu} - \bm{v}_{0}^{-} \cdot \bm{\nu} = 0.
\end{align}
We have from $(\ref{Asymp:order})_{I}^{-2}$,
\begin{align}
\label{Inner:mu:-2}
\pd_{z}(m(\Phi_{0}) \pd_{z}\Xi_{0}) = 0.
\end{align}
Upon integrating and using the matching condition (\ref{MatchingCond2}) applied to $\Xi_{0}$, we obtain
\begin{align*}
m(\Phi_{0}) \pd_{z} \Xi_{0}(t,s,z) = 0 \text{ for all } \abs{z} < l.
\end{align*}
Since $\abs{\Phi_{0}(z)} < 1$ for $\abs{z} < l$ and $m(\Phi_{0}) > 0$, we have
\begin{align}
\label{pdzXi0zero}
\pd_{z} \Xi_{0}(t,s,z) = 0 \text{ for all } \abs{z} < l.
\end{align}
Thus, integrating once more with respect to $z$ from $-l$ to $l$, and applying the matching condition (\ref{MatchingCond1}) to $\Xi_{0}$, we obtain 
\begin{align}
\label{mu0cont}
\jump{\mu_{0}} = 0.
\end{align}
To leading order, the nutrient equation $(\ref{Asymp:sigma})_{I}^{-2}$ yields
\begin{align}
\label{Inner:sigma:-2}
\pd_{z}(n(\Phi_{0}) \chi_{\sigma} \pd_{z} C_{0}) - (n(\Phi_{0}) \chi_{\varphi} \Phi_{0}'(z))' = 0.
\end{align}
Integrating and using the matching condition (\ref{MatchingCond2}) applied to both $C_{0}$ and $\Phi_{0}$ leads to
\begin{align*}
n(\Phi_{0}) (\chi_{\sigma} \pd_{z}C_{0} - \chi_{\varphi} \Phi_{0}'(z)) = 0 \text{ for all } \abs{z} < l.
\end{align*}
As $n(\Phi_{0}) > 0$, we see that
\begin{align}
\label{pdzC0pdzPhi0}
\chi_{\sigma} \pd_{z}C_{0}(t,s,z) = \chi_{\varphi} \Phi_{0}'(z) \text{ for all } \abs{z} < l.
\end{align}
Integrating once more with respect to $z$ from $-l$ to $l$, and applying the matching condition (\ref{MatchingCond1}) to $C_{0}$ and $\Phi_{0}$ then gives
\begin{align}
\label{jumpsigma0}
\jump{\sigma_{0}} = \frac{\chi_{\varphi}}{\chi_{\sigma}} \jump{\varphi_{0}} = 2\frac{\chi_{\varphi}}{\chi_{\sigma}}.
\end{align}
Lastly, $(\ref{Asymp:velo})_{I}^{-1}$ yields
\begin{align}
\label{Inner:velo:-1}
\pd_{z} P_{0} = (\Xi_{0}  + \chi_{\varphi} C_{0}) \Phi_{0}'.
\end{align}
Integrating and applying the matching condition (\ref{MatchingCond1}) to $P_{0}$ and $\Xi_{0}$ leads to
\begin{align}
\label{Pressurerelation1}
\jump{p_{0}} = 2 \mu_{0} + \chi_{\varphi} \int_{-l}^{l} C_{0}(t,s,z) \Phi_{0}'(z) \dz.
\end{align}
Thanks to (\ref{pdzC0pdzPhi0}), we see that
\begin{align}
\notag \int_{-l}^{l} C_{0} \Phi_{0}' \dz & = \frac{\chi_{\sigma}}{\chi_{\varphi}} \int_{-l}^{l} C_{0} \pd_{z} C_{0} \dz = \frac{\chi_{\sigma}}{\chi_{\varphi}} \int_{-l}^{l} \pd_{z} \left ( \frac{\abs{C_{0}}^{2}}{2} \right ) \dz \\
& = \frac{\chi_{\sigma}}{2\chi_{\varphi}} \left [ \abs{C_{0}}^{2} \right ]_{-l}^{l} = \frac{\chi_{\sigma}}{2\chi_{\varphi}}  \jump{\abs{\sigma_{0}}^{2}}. \label{int:C0pdzPhi0}
\end{align}
Then, (\ref{Pressurerelation1}) becomes
\begin{align}\label{pressure:kappa:2}
\jump{p_{0}} = 2 \mu_{0} + \frac{\chi_{\sigma}}{2} \jump{\abs{\sigma_{0}}^{2}}.
\end{align}
\subsubsection{Expansions to first order}
For the double-well potential, to first order, we obtain from $(\ref{Asym:mu})_{I}^{0}$,
\begin{align}
\label{Inner:mu:0}
\beta \Psi''(\Phi_{0}) \Phi_{1} - \beta \pd_{zz} \Phi_{1} + \beta \kappa \Phi_{0}' - \chi_{\varphi} C_{0} = \Xi_{0}.
\end{align}
We multiply (\ref{Inner:mu:0}) with $\Phi_{0}'$ and integrate with respect to $z$ from $-\infty$ to $\infty$, which gives
\begin{align}
\notag & \; \int_{-\infty}^{\infty} \Xi_{0}(t,s) \Phi_{0}'(z) \dz \\
= & \; \int_{-\infty}^{\infty} \beta (\Psi'(\Phi_{0}))' \Phi_{1} - \beta \pd_{zz} \Phi_{1} \Phi_{0}' + \beta \kappa \abs{\Phi_{0}'}^{2} - \chi_{\varphi} C_{0} \Phi_{0}' \dz. \label{Inner:mu:0:int}
\end{align}
Applying integration by parts and the matching conditions (\ref{MatchingCond1}) and (\ref{MatchingCond2}) applied to $\Phi_{0}$, and using that $\Psi'(\pm 1) = 0$, we see that
\begin{align*}
& \; \int_{-\infty}^{\infty} (\Psi'(\Phi_{0}))' \Phi_{1} - \pd_{zz} \Phi_{1} \Phi_{0}' \dz \\
= & \; \underbrace{[\Psi'(\Phi_{0}) \Phi_{1} - \pd_{z} \Phi_{1} \Phi_{0}']_{-\infty}^{\infty}}_{=0 \text{ by } (\ref{MatchingCond1}), (\ref{MatchingCond2})} - \int_{-\infty}^{\infty} \pd_{z} \Phi_{1} \underbrace{(\Psi'(\Phi_{0}) - \Phi_{0}'')}_{=0 \text{ by } (\ref{ODE})} \dz,
\end{align*}
and so the first two terms on the right hand side of (\ref{Inner:mu:0:int}) are zero.  Then, using (\ref{equipartition:gamma}),  (\ref{pdzXi0zero}), and (\ref{int:C0pdzPhi0}), we obtain from (\ref{Inner:mu:0:int}),
\begin{align}
\label{Inner:mu:0:int:interfacelaw}
2 \mu_{0} = \beta \kappa \gamma - \chi_{\varphi} \int_{-\infty}^{\infty} C_{0} \Phi_{0}' \dz = \beta \gamma \kappa - \frac{\chi_{\sigma}}{2} \jump{\abs{\sigma_{0}}^{2}}.
\end{align}
Moreover, together with (\ref{pressure:kappa:2}), we obtain
\begin{align}
\label{pressure:kappa}
\jump{p_{0}} =  \beta \gamma \kappa.
\end{align}
Meanwhile, for the double-obstacle potential, to first order, we obtain from $(\ref{Asym:mu:variationalinequ})_{I}^{0}$,
\begin{align}
\label{Inner:DO:mu:0}
\int_{\Omega} (- \Xi_{0} - \chi_{\varphi} C_{0} - \beta \pd_{zz} \Phi_{1} - \beta \Phi_{1} + \kappa \beta \Phi_{0}') (\psi - \Phi_{0}) \dx \geq 0 \text{ for all } \abs{\psi} \leq 1.
\end{align}
Since $\abs{\Phi_{0}(z)} < 1$ for $\abs{z} < \frac{\pi}{2}$, we can test with $\psi = \Phi_{0} + \lambda$ with either non-positive or non-negative $\lambda \in \mathcal{K}$, leading to the equality
\begin{align*}
- \Xi_{0} - \chi_{\varphi} C_{0} - \beta \pd_{zz} \Phi_{1} - \beta \Phi_{1} + \kappa \beta \Phi_{0}' = 0.
\end{align*}
Multiplying with $\Phi_{0}'$ and integrating with respect to $z$ from $-\frac{\pi}{2}$ to $\frac{\pi}{2}$, and applying matching conditions leads to
\begin{align}\label{DO:Phi1Solvability}
-2 \mu_{0} + \beta \kappa \int_{-\frac{\pi}{2}}^{\frac{\pi}{2}} \abs{\Phi_{0}'}^{2} \dx - \int_{-\frac{\pi}{2}}^{\frac{\pi}{2}} \chi_{\varphi} C_{0} \Phi_{0}' \dz = \beta \int_{-\frac{\pi}{2}}^{\frac{\pi}{2}}  \pd_{zz} \Phi_{1} \Phi_{0}' + \Phi_{1} \Phi_{0}' \dz.
\end{align}
Upon integrating by parts and using (\ref{ODE:DO}), the matching conditions (\ref{MatchingCond2}) for $\Phi_{0}$, (\ref{MatchingCond3}) for $\Phi_{1}$, and (\ref{DO:Phi1}), we see that
\begin{align*}
\int_{-\frac{\pi}{2}}^{\frac{\pi}{2}} \pd_{zz} \Phi_{1} \Phi_{0}' + \Phi_{1} \Phi_{0}' \dz = [\pd_{z}\Phi_{1} \Phi_{0}' + \Phi_{1}\Phi_{0}]_{z=-\frac{\pi}{2}}^{z=\frac{\pi}{2}} - \int_{-\frac{\pi}{2}}^{\frac{\pi}{2}} (\Phi_{0}'' + \Phi_{0})\pd_{z}\Phi_{1} \dz = 0.
\end{align*}
Then, using (\ref{equipartition:gamma}), and (\ref{int:C0pdzPhi0}), we obtain from (\ref{DO:Phi1Solvability}) the following solvability condition for $\Phi_{1}$:
\begin{align*}
2 \mu_{0} = \beta \gamma \kappa - \frac{\chi_{\sigma}}{2} \jump{\abs{\sigma_{0}}^{2}}.
\end{align*}
Lastly, thanks to (\ref{pdzXi0zero}), we obtain from $(\ref{Asymp:order})_{I}^{-1}$ and $(\ref{Asymp:sigma})_{I}^{-1}$, respectively,
\begin{align}
(-\mathcal{V} + \bm{V}_{0} \cdot \bm{\nu}) \Phi_{0}' & = \pd_{z}(m(\Phi_{0}) \pd_{z} \Xi_{1}), \label{Inner:order:-1} 
\end{align}
and
\begin{align}
\notag & \; (-\mathcal{V} + \bm{V}_{0} \cdot \bm{\nu}) \pd_{z} C_{0} \\
\notag = & \; \pd_{z}( n(\Phi_{0}) (\chi_{\sigma} \pd_{z} C_{1} - \chi_{\varphi} \pd_{z} \Phi_{1})) \\
\notag + & \; \pd_{z} (n'(\Phi_{0}) \Phi_{1} \underbrace{(\chi_{\sigma} \pd_{z} C_{0} - \chi_{\varphi} \Phi_{0}')}_{=0 \text{ by } (\ref{pdzC0pdzPhi0})}) + \div_{\Sigma_{0}} (n(\Phi_{0}) \underbrace{(\chi_{\sigma} \pd_{z} C_{0} - \chi_{\varphi} \Phi_{0}')}_{=0 \text{ by } (\ref{pdzC0pdzPhi0})} \bm{\nu}) \\
= & \; \pd_{z}( n(\Phi_{0}) (\chi_{\sigma} \pd_{z} C_{1} - \chi_{\varphi} \pd_{z} \Phi_{1})). \label{Inner:sigma:-1}
\end{align}
Thanks to (\ref{Inner:div:-1}), upon integrating from $-l$ to $l$ with respect to $z$, and applying the matching condition (\ref{MatchingCond3}) to $\Xi_{1}$, we obtain from (\ref{Inner:order:-1})
\begin{align}
\label{jumpnablamu}
2(-\mathcal{V} + \bm{v}_{0} \cdot \bm{\nu}) =  \jump{m(\varphi_{0}) \nabla \mu_{0}} \cdot \bm{\nu}.
\end{align}
Similarly, thanks to $\nabla \varphi_{0} = \bm{0}$, upon integrating from $-l$ to $l$ with respect to $z$, and applying the matching condition (\ref{MatchingCond3}) to $C_{1}$ and $\Phi_{1}$, we obtain from (\ref{Inner:sigma:-1})
\begin{align}
\label{jumpnablasigma}
(-\mathcal{V} + \bm{v}_{0} \cdot \bm{\nu}) \jump{\sigma_{0}} = \chi_{\sigma} \jump{n(\varphi_{0}) \nabla \sigma_{0}} \cdot \bm{\nu}.
\end{align}

\subsubsection{Sharp interface limit for Model (\ref{Model:CHDarcy:Asymp})}
In summary, we obtain the following sharp interface limit from Model (\ref{Model:CHDarcy:Asymp}):  
\begin{subequations}\label{SI:CHDarcy}
\begin{alignat}{3}
\bm{v}_{0} & = - K \nabla p_{0} && \text{ in } \Omega_{T} \cup \Omega_{H}, \label{SI:CHDarcy:velo} \\
\div \bm{v}_{0} &= \overline{\rho}_{1}^{-1} \Gamma_{1}(\sigma_{0}^{T}, 1, \mu_{0}^{T}) + \overline{\rho}_{2}^{-1} \Gamma_{2} (\sigma_{0}^{T}, 1, \mu_{0}^{T}) && \text{ in } \Omega_{T}, \\
\div \bm{v}_{0} &= \overline{\rho}_{1}^{-1} \Gamma_{1}(\sigma_{0}^{H}, -1, \mu_{0}^{H}) + \overline{\rho}_{2}^{-1} \Gamma_{2} (\sigma_{0}^{H}, -1, \mu_{0}^{H}) && \text{ in } \Omega_{H},
\label{SI:CHDarcy:pressure} \\
-m(1) \Laplace \mu_{0}^{T} &= - 2 \overline{\rho}_{1}^{-1}\Gamma_{1}(\sigma_{0}^{T}, 1, \mu_{0}^{T}) && \text{ in } \Omega_{T} , \label{SI:CHDarcy:mu} \\
-m(-1) \Laplace \mu_{0}^{H} &= 2 \overline{\rho}_{2}^{-1}\Gamma_{2}(\sigma_{0}^{H}, -1, \mu_{0}^{H}) && \text{ in } \Omega_{H} , \\
\pd_{t} \sigma_{0}^{T} + \div (\sigma_{0}^{T} \bm{v}_{0}) &= n(1) \chi_{\sigma} \Laplace \sigma_{0}^{T} - \mathcal{S}(\sigma_{0}^{T}, 1, \mu_{0}^{T}) && \text{ in } \Omega_{T} , \label{SI:CHDarcy:sigma} \\
\pd_{t} \sigma_{0}^{H} + \div (\sigma_{0}^{H} \bm{v}_{0}) &= n(-1) \chi_{\sigma} \Laplace \sigma_{0}^{H} - \mathcal{S}(\sigma_{0}^{H}, -1, \mu_{0}^{H}) && \text{ in } \Omega_{H} ,
\end{alignat}
\end{subequations}
together with the free boundary conditions
\begin{subequations}
\begin{align}
\jump{\bm{v}_{0}} \cdot \bm{\nu} = 0, \quad \jump{\sigma_{0}} = 2 \frac{\chi_{\varphi}}{\chi_{\sigma}}, \quad \jump{p_{0}} = \beta \gamma \kappa   & \text{ on } \Sigma_{0}, \label{SI:CHDarcy:jump} \\
\jump{\mu_{0}} = 0, \quad 2 \mu_{0} = \beta \gamma \kappa  - \frac{\chi_{\sigma}}{2} \jump{\abs{\sigma_{0}}^{2}} & \text{ on } \Sigma_{0}, \label{SI:CHDarcy:jumpmu} \\
2(-\mathcal{V} + \bm{v}_{0} \cdot \bm{\nu})  = (m(1)\nabla \mu_{0}^{T} - m(-1) \nabla \mu_{0}^{H}) \cdot \bm{\nu} & \text{ on } \Sigma_{0}, \label{SI:CHDarcy:nablamu} \\
2 \frac{\chi_{\varphi}}{\chi_{\sigma}} (-\mathcal{V} + \bm{v}_{0} \cdot \bm{\nu}) =  \chi_{\sigma} (n(1)\nabla \sigma_{0}^{T} - n(-1) \nabla \sigma_{0}^{H}) \cdot \bm{\nu} & \text{ on } \Sigma_{0}, \label{SI:CHDarcy:nablasigma}
\end{align}
\end{subequations}
where $\gamma$ is defined in (\ref{defn:gammaconst}).  Note that we can write
\begin{align}\label{jumpsigmasquare}
\jump{\abs{\sigma_{0}}^{2}} = \jump{\sigma_{0}} (\sigma_{0}^{T} + \sigma_{0}^{H}) =: 2 \overline{\sigma}_{0} \jump{\sigma_{0}},
\end{align}
where $\overline{\sigma}_{0}:= \tfrac{1}{2} (\sigma_{0}^{T} + \sigma_{0}^{H})$ denotes the average of the nutrient concentrations from both sides of $\Sigma_{0}$.  Thus, using (\ref{jumpsigma0}), we can rewrite (\ref{Inner:mu:0:int:interfacelaw}) and $(\ref{SI:CHDarcy:jumpmu})_{2}$ as
\begin{align}\label{Sitka:mu:interfacelaw}
 \mu_{0} = \frac{1}{2} \beta \gamma \kappa - \overline{\sigma}_{0} \chi_{\varphi} \text{ on } \Sigma_{0}.
\end{align}

\subsection{Specific sharp interface models}
In this section, we take $\Psi$ as the double-well potential.
\subsubsection{Sharp interface limit of the new active transport model}
Choosing as before $N(\varphi, \sigma) = \frac{\chi_{\sigma}}{2} \abs{\sigma}^{2} + \chi_{\varphi} \sigma(1-\varphi)$ and
\begin{align}\label{Sitka:choice}
\Gamma(\sigma, \varphi) = (\mathcal{P} \sigma - \mathcal{A}) h(\varphi), \quad \mathcal{S}(\sigma, \varphi) = \mathcal{C} \sigma h(\varphi), \\
m(\varphi) = m_{0} > 0, \quad n(\varphi) = n_{0} > 0, 
\end{align}
for some positive constants $\mathcal{P}$, $\mathcal{A}$ and $\mathcal{C}$ in Model (\ref{Model:CHDarcy:Gamma}), we obtain the Cahn--Hilliard--Darcy model (\ref{Intro:CHDarcy:zeroexcesstotalmass}) in Section \ref{sec:Intro}.  Then, the sharp interface limit of Model (\ref{Model:CHDarcy:Gamma}) with (\ref{Sitka:choice}) is given by
\begin{subequations}\label{Sitka:SI:CHDarcy:Gamma}
\begin{align}
-\Laplace p_{0} & = \begin{cases}
\frac{\alpha}{K}(\mathcal{P} \sigma_{0} - \mathcal{A}) & \text{ in } \Omega_{T}, \\
0 & \text{ in } \Omega_{H},
\end{cases} \\
-m_{0} \Laplace \mu_{0} & = \begin{cases}
(\overline{\rho}_{S} - \alpha)(\mathcal{P} \sigma_{0} - \mathcal{A}) & \text{ in } \Omega_{T}, \\
0 & \text{ in } \Omega_{H}, 
\end{cases} \\
\pd_{t} \sigma_{0} - \div (K\sigma_{0} \nabla p_{0}) - n_{0} \chi_{\sigma} \Laplace \sigma_{0} & = \begin{cases}
- \mathcal{C} \sigma_{0} & \text{ in } \Omega_{T}, \\
0 & \text{ in } \Omega_{H}, 
\end{cases} \label{Sitka:SI:CHDarcy:Gamma:sigma} \\
K\jump{\nabla p_{0}} \cdot \bm{\nu} = 0, & \quad \jump{\sigma_{0}} = 2 \frac{\chi_{\varphi}}{\chi_{\sigma}}, \quad \jump{p_{0}} = \beta \gamma \kappa   \text{ on } \Sigma_{0},  \label{Sitka:SI:CHDarcy:Gamma:jump} \\
\jump{\mu_{0}} = 0, & \quad \mu_{0} + \overline{\sigma}_{0} \chi_{\varphi} = \frac{1}{2} \beta \gamma \kappa  \text{ on } \Sigma_{0}, \label{Sitka:SI:CHDarcy:Gamma:interface:mu} \\
-2(\mathcal{V} + K \nabla p_{0} \cdot \bm{\nu}) & = m_{0} \jump{\nabla \mu_{0}} \cdot \bm{\nu} \text{ on } \Sigma_{0}, \label{Sitka:SI:CHDarcy:Gamma:nablamu} \\
-2 \frac{\chi_{\varphi}}{\chi_{\sigma}} (\mathcal{V} + K \nabla p_{0} \cdot \bm{\nu})& = n_{0} \chi_{\sigma} \jump{\nabla \sigma_{0}} \cdot \bm{\nu}  \text{ on } \Sigma_{0}. \label{Sitka:SI:CHDarcy:Gamma:nablasigma}
\end{align}
\end{subequations}
The active transport term $n(\varphi) \nabla (\chi_{\varphi} \varphi)$ in the flux for the nutrient results in the jump term $2 \frac{\chi_{\varphi}}{\chi_{\sigma}}$ in $(\ref{Sitka:SI:CHDarcy:Gamma:jump})_{2}$ which is a new feature of the proposed model.

\subsubsection{Linear constitutive laws for chemical reactions}
Let us consider Model (\ref{Model:Specific:NoVelo:2}) with $m(\varphi) = n(\varphi) = 1$, and $P(\varphi)$ defined as in (\ref{HawkinsPvarphi}).  Then, we obtain that
\begin{equation}
\begin{aligned}
P(\varphi_{0}) (\chi_{\sigma} \sigma_{0} + \chi_{\varphi}(1-\varphi_{0}) - \mu_{0}) = \begin{cases}
2 \delta P_{0}(\chi_{\sigma} \sigma_{0} - \mu_{0}) & \text{ in } \Omega_{T}, \\
0 & \text{ in } \Omega_{H}.
\end{cases}
\end{aligned}
\end{equation}
This was the setting introduced in \cite{article:HawkinsZeeOden12}.  Hence, we obtain from (\ref{Model:Specific:NoVelo:2}) the following sharp interface limit:
\begin{subequations}
\begin{align}
- \Laplace \mu_{0} & = \begin{cases}
2\delta P_{0}(\chi_{\sigma} \sigma_{0} - \mu_{0}), & \text{ in } \Omega_{T}, \\
0 & \text{ in } \Omega_{H}, \\
\end{cases} \\
\pd_{t} \sigma_{0} - \chi_{\sigma} \Laplace \sigma_{0} & = \begin{cases} -2 \delta P_{0} (\chi_{\sigma} \sigma_{0} - \mu_{0}), & \text{ in } \Omega_{T}, \\
0 & \text{ in } \Omega_{H}, \\
\end{cases} \\
\mu_{0} = \frac{1}{2} \beta \gamma \kappa - \overline{\sigma}_{0} \chi_{\varphi},  \quad \jump{\mu_{0}} = 0, & \quad \jump{\sigma_{0}} = 2\frac{\chi_{\varphi}}{\chi_{\sigma}} \text{ on } \Sigma_{0}, \\
-2 \mathcal{V} = \jump{\nabla \mu_{0}}\cdot \bm{\nu} , \quad -2\frac{\chi_{\varphi}}{\chi_{\sigma}} \mathcal{V} & = \chi_{\sigma} \jump{\nabla \sigma_{0}} \cdot \bm{\nu} \text{ on } \Sigma_{0}.
\end{align}
\end{subequations}
We point out that the diffuse interface model studied in \cite{article:HilhorstKaampmannNguyenZee15} takes a different mass transition term $\Gamma$ and a different consumption term $\mathcal{S}$.  More precisely, the following choices are considered:
\begin{align}\label{HilhorstGammaS}
\Gamma = \frac{1}{2 \eps} P(\varphi)(\sigma - \delta \mu), \quad \mathcal{S} = \frac{1}{\eps} P(\varphi) (\sigma - \delta \mu),
\end{align}
where \footnote{In \cite{article:HilhorstKaampmannNguyenZee15}, the choice of $P(\varphi)$ is actually
\begin{align}\label{HilhorstPvarphi}
P(\varphi) := \begin{cases}
2 P_{0} (\Psi(\varphi))^{\frac{1}{2}} & \text{ if } \varphi \in [-1,1], \\
0 & \text{ otherwise }.
\end{cases}
\end{align}
This presents some difficulties in the analysis of the outer expansions, as $P'(\pm 1) \neq 0$ for the choice $\Psi(\varphi) = (1-\varphi^{2})^{2}$.  Thus, we do not recover (\ref{SI:Hilhorst:mu}) and (\ref{SI:Hilhorst:sigma}).  However, the formal analysis in \cite{article:HilhorstKaampmannNguyenZee15} is different compared to what we present here, and it turns out that the analysis in \cite{article:HilhorstKaampmannNguyenZee15} has to be modified and will only work for (\ref{alternate:Hilhorst:choiceP}).
}
\begin{align}\label{alternate:Hilhorst:choiceP}
P(\varphi) := \begin{cases}
2 P_{0} \Psi(\varphi) & \text{ if } \varphi \in [-1,1], \\
0 & \text{ otherwise }.
\end{cases}
\end{align}
The term $\frac{1}{\eps} \Psi(\varphi)$ acts as a regularisation on the Hausdorff measure restricted to the limiting hypersurface $\Sigma_{0}$, and hence the reaction term $\sigma - \delta \mu$ will appear in the interfacial relations of $\mu$ and $\sigma$ rather than in the bulk equations.  More precisely, we obtain
\begin{subequations}\label{SI:Hilhorst}
\begin{align}
\Laplace \mu_{0} = 0 & \text{ in } \Omega_{T} \cup \Omega_{H}, \label{SI:Hilhorst:mu} \\
\pd_{t} \sigma_{0} = \Laplace \sigma_{0} & \text{ in } \Omega_{T} \cup \Omega_{H}, \label{SI:Hilhorst:sigma} \\
\jump{\mu_{0}} = 0, \quad \jump{\sigma_{0}} = 0, \quad 2 \mu_{0} = \beta \gamma \kappa & \text{ on } \Sigma_{0}, \\
-2 \mathcal{V} = \jump{\nabla \mu_{0}} \cdot \bm{\nu} + P_{0} \gamma (\sigma_{0} - \delta \mu_{0}) & \text{ on } \Sigma_{0}, \label{SI:Hilhorst:nablamu} \\
0 = \jump{\nabla \sigma_{0}} \cdot \bm{\nu} - P_{0} \gamma (\sigma_{0} - \delta \mu_{0}) & \text{ on } \Sigma_{0}, \label{SI:Hilhorst:nablasigma}
\end{align}
\end{subequations}
as a sharp interface limit of Model (\ref{Model:CHDarcy:NoVelo:EqualDensities}) with $\overline{\rho} = 1$, $m(\varphi) = n(\varphi) = 1$, $\Gamma$ and $\mathcal{S}$ as in (\ref{HilhorstGammaS}) with $P(\varphi)$ chosen as in (\ref{alternate:Hilhorst:choiceP}) and $N(\varphi, \sigma) = \frac{1}{2} \abs{\sigma}^{2}$.  This is similar to the sharp interface limit (Eq. (1.9)) of \cite{article:HilhorstKaampmannNguyenZee15} with $\alpha = 0$.

\subsubsection{The limit of vanishing active transport}\label{sec:compare:Cristini}
We consider Model (\ref{Model:Specific:NoVelo:1}) with a quasi-steady nutrient (i.e., neglecting the left hand side of (\ref{Model:Specific:NoVelo:1:c})), with positive constants $D$ and $\lambda$, and the interpolation function $h(\varphi) = \frac{1}{2}(1+\varphi)$, we set
\begin{subequations}\label{Cristini:Specific}
\begin{align}
\mathcal{D}(\varphi) := \frac{1 + \varphi}{2} + D \frac{1-\varphi}{2} = \frac{1}{2} (1 + D) + \frac{\varphi}{2} (1-D), \\
m(\varphi) = \frac{1}{2}(1+\varphi)^{2}, \quad n(\varphi) = \lambda \mathcal{D}(\varphi) \chi_{\varphi}^{-1}, \quad \chi_{\sigma} = \lambda^{-1} \chi_{\varphi},
\end{align}
\end{subequations}
so that, we obtain
\begin{subequations}\label{Model:Cristini}
\begin{align}
\pd_{t} \varphi & = \div (\tfrac{1}{2}(1+\varphi)^{2} \nabla \mu) + \mathcal{P} \sigma (\varphi + 1) - \mathcal{A} (\varphi + 1) \label{Model:Cristini:order}, \\
\mu & = \frac{\beta}{\eps} \Psi'(\varphi) - \beta \eps \Laplace \varphi - \chi_{\varphi} \sigma \label{Model:Cristini:mu}, \\
0 & = \div(\mathcal{D}(\varphi) \nabla \sigma) - \lambda \div(\mathcal{D}(\varphi) \nabla \varphi) - \frac{1}{2} \mathcal{C} \sigma (\varphi + 1) \label{Model:Cristini:sigma}.
\end{align}
\end{subequations}
The specific choice (\ref{Cristini:Specific}) allows us to control the influence of the active transport term $n(\varphi) \chi_{\varphi} \nabla \varphi$ via the parameter $\lambda$, while preserving the chemotaxis term $-\chi_{\varphi} \sigma$ in (\ref{Model:Specific:NoVelo:1:b}).  Hence, we have ``decoupled'' chemotaxis and active transport.  

Moreover, if we consider Eq. (68)-(70) of \cite{article:CristiniLiLowengrubWise09} with the choice $\phi = \frac{1}{2} (1+\varphi)$, $\mathcal{G} = 1$, and the rescaling $\mu \mapsto \eps \mu$, the resulting phase field model almost coincides with Model (\ref{Model:Cristini}) with the exception of the additional term $\lambda  \div (\mathcal{D}(\varphi) \nabla \varphi)$ in (\ref{Model:Cristini:sigma}).

We briefly state the derivation of the sharp interface limit for Model (\ref{Model:Cristini}).  From $(\ref{Model:Cristini:mu})_{O}^{-1}$ we have $\varphi_{0} = \pm 1$ and the domains $\Omega_{T}$ and $\Omega_{H}$.  From $(\ref{Model:Cristini:order})_{O}^{0}$ and $(\ref{Model:Cristini:sigma})_{O}^{0}$ we obtain
\begin{align*}
0 = \div (\tfrac{1}{2}(1 + \varphi_{0})^{2} \nabla \mu_{0}) + \mathcal{P} \sigma_{0} (\varphi_{0} + 1) - \mathcal{A} (\varphi_{0} + 1) & \text{ in } \Omega_{T} \cup \Omega_{H}, \\
0 = \div (\mathcal{D}(\varphi_{0}) \nabla \sigma_{0}) - \tfrac{1}{2} \mathcal{C} \sigma_{0}(\varphi_{0} + 1) & \text{ in } \Omega_{T} \cup \Omega_{H}.
\end{align*}
From the leading order inner expansion $(\ref{Model:Cristini:mu})_{I}^{-1}$, we obtain (\ref{Inner:mu:-1}), and subsequently the profile (\ref{Phi0profile}) and the equipartition of energy (\ref{equipartition}).  From $(\ref{Model:Cristini:sigma})_{I}^{-2}$ we have 
\begin{align*}
\pd_{z}(\mathcal{D}(\Phi_{0}) \pd_{z} C_{0} - \lambda \mathcal{D}(\Phi_{0}) \Phi_{0}') = 0.
\end{align*}
Integrating and using the matching conditions (\ref{MatchingCond2}), we obtain
\begin{align*}
\mathcal{D}(\Phi_{0}) (\pd_{z} C_{0} - \lambda\Phi_{0}') = 0.
\end{align*}
Since $\mathcal{D}(\Phi_{0}) > 0$ for $\abs{\Phi_{0}} < 1$, we obtain
\begin{align}\label{pdzC0lambdapdzPhi0}
\pd_{z}C_{0}(t,s,z) = \lambda \Phi_{0}'(z) \text{ for all } \abs{z} < \infty,
\end{align}
and upon matching, we obtain
\begin{align}\label{Cristini:jumpsigmalambda}
\jump{\sigma_{0}} = 2 \lambda.
\end{align}
While from $(\ref{Model:Cristini:order})_{I}^{-2}$ we obtain
\begin{align*}
 \pd_{z} ((1+ \Phi_{0})^{2} \pd_{z} \Xi_{0}) = 0.
\end{align*}
Integrating and using the matching condition (\ref{MatchingCond2}) applied to $\Xi_{0}$ we deduce that
\begin{align*}
(1 + \Phi_{0}(z))^{2} \pd_{z} \Xi_{0}(t,s,z) = 0 \text{ for all } \abs{z} < \infty.
\end{align*}
Since $\abs{\Phi_{0}(z)} < 1$ for $\abs{z} < \infty$, we obtain that $\pd_{z} \Xi_{0}(t,s,z) = 0$ for $\abs{z} < \infty$.

To first order, we obtain from $(\ref{Model:Cristini:mu})_{I}^{0}$,
\begin{align*}
\Xi_{0} = \beta \Psi'(\Phi_{0}) - \beta \pd_{zz} \Phi_{1} + \beta \kappa \Phi_{0}' - \chi_{\varphi} C_{0}.
\end{align*}
Multiplying by $\Phi_{0}'$, using that $\Xi_{0}$ is independent of $z$, and applying integration by parts and matching conditions (\ref{MatchingCond1}) and (\ref{MatchingCond2}) to $\Phi_{0}$, we obtain, in the same spirit as (\ref{Inner:mu:0:int}), 
\begin{align*}
2\mu_{0} = \beta \gamma \kappa  - \int_{\R} \chi_{\varphi} C_{0} \Phi_{0}' \dz = \beta \gamma \kappa  - \frac{\chi_{\varphi}}{\lambda}\frac{1}{2} \jump{\abs{\sigma_{0}}^{2}},
\end{align*}
where we have used (\ref{pdzC0lambdapdzPhi0}).  Applying (\ref{jumpsigmasquare}), and (\ref{Cristini:jumpsigmalambda}), we see that
\begin{align*}
2 \mu_{0} =  \beta \gamma \kappa - \frac{\chi_{\varphi}}{\lambda} \overline{\sigma_{0}} \jump{\sigma_{0}} =  \beta \gamma \kappa - 2 \chi_{\varphi} \overline{\sigma}_{0},
\end{align*}
where we recall that $\overline{\sigma}_{0} := \frac{1}{2} (\sigma_{0}^{T} + \sigma_{0}^{H})$ is the average of the nutrient concentration at the interface.  Meanwhile, thanks to (\ref{pdzC0lambdapdzPhi0}) we obtain from $(\ref{Model:Cristini:sigma})_{I}^{-1}$,
\begin{align*}
0 & = \pd_{z}(\mathcal{D}(\Phi_{0}) (\pd_{z} C_{1} - \lambda \pd_{z}\Phi_{1}) + \mathcal{D}'(\Phi_{0})\Phi_{1}(\pd_{z}C_{0} - \lambda \Phi_{0}'))\\
&  =  \pd_{z}(\mathcal{D}(\Phi_{0}) (\pd_{z} C_{1} - \lambda \pd_{z}\Phi_{1})).
\end{align*}
Integrating with respect to $z$ from $-\infty$ to $\infty$ and applying the matching condition (\ref{MatchingCond3}) to $C_{1}$ and $\Phi_{1}$ leads to
\begin{align*}
0 =  \jump{\mathcal{D}(\varphi_{0}) \nabla \sigma_{0}} \cdot \bm{\nu}.
\end{align*}
Lastly, thanks to $\pd_{z} \Xi_{0} = 0$, we obtain from $(\ref{Model:Cristini:order})_{I}^{-1}$,
\begin{align*}
-\mathcal{V} \Phi_{0}' = \frac{1}{2} \pd_{z} ((1 + \Phi_{0})^{2} \pd_{z} \Xi_{1}).
\end{align*}
Integrating from $-\infty$ to $\infty$ with respect to $z$, and applying the matching condition (\ref{MatchingCond3}) to $\Xi_{1}$ gives
\begin{align*}
-2 \mathcal{V} = 2 \nabla \mu_{0}^{T} \cdot \bm{\nu}.
\end{align*}
Thus, the sharp interface limit of Model (\ref{Model:Cristini}) is
\begin{subequations}\label{SI:Cristini}
\begin{align}
- \Laplace \mu_{0}^{T} & = \mathcal{P} \sigma_{0}^{T} - \mathcal{A} \text{ in } \Omega_{T}, \\
\Laplace \sigma_{0} & = \begin{cases}
\mathcal{C} \sigma_{0} & \text{ in } \Omega_{T}, \\
0 & \text{ in } \Omega_{H},
\end{cases} \\
\jump{\sigma_{0}} = 2\lambda, \quad 2 \mu_{0} & = \beta \gamma \kappa - \chi_{\varphi} (\sigma_{0}^{T} + \sigma_{0}^{H}) \text{ on } \Sigma_{0}, \label{SI:Cristini:jump} \\
0 = (\nabla \sigma_{0}^{T} - D \nabla \sigma_{0}^{H}) \cdot \bm{\nu}, \quad -\mathcal{V} & =  \nabla \mu_{0}^{T} \cdot \bm{\nu} \text{ on } \Sigma_{0}. \label{SI:Cristini:nablamu}
\end{align}
\end{subequations}
In addition, we can use $(\ref{SI:Cristini:jump})_{1}$ to rewrite
$(\ref{SI:Cristini:jump})_{2}$ as
\begin{align}
\label{SI:Cristini:kappa}
2 \mu_{0} = \beta \gamma \kappa- \chi_{\varphi} (2 \sigma_{0}^{T} - 2 \lambda).
\end{align}
Next, sending $\lambda \to 0$ in (\ref{SI:Cristini}) leads to
\begin{align}\label{SI:Cristini:lambdazerolimit}
\jump{\sigma_{0}} = 0, \quad 2 \mu_{0} = \beta \gamma \kappa  - 2 \chi_{\varphi} \sigma_{0},
\end{align}
and we define the bulk velocity and pressure via the relations: 
\begin{align}
\label{defn:Cristini:bulkveloPressure}
\bm{v} := -\nabla (p - \chi_{\varphi} \sigma_{0}), \quad p := \mu_{0}^{T} + \chi_{\varphi} \sigma_{0}.
\end{align}
Then, we deduce that
\begin{align}\label{SI:Cristini:velo}
\bm{v} =  -\nabla \mu_{0}^{T}, \quad \div \bm{v} =  -\Laplace \mu_{0}^{T} = (\mathcal{P} \sigma_{0}^{T} - \mathcal{A}) \text{ in } \Omega_{T},
\end{align}
and from (\ref{SI:Cristini:nablamu}) and \eqref{SI:Cristini:lambdazerolimit},
\begin{align}
- \mathcal{V} =  \nabla \mu_{0}^{T} \cdot \bm{\nu} = -\bm{v} \cdot \bm{\nu} = \nabla (p - \chi_{\varphi} \sigma_{0}) \cdot \bm{\nu} & \text{ on } \Sigma_{0}, \label{SI:Cristini:NormalVelo} \\
p = \mu_{0} + \chi_{\varphi} \sigma_{0} = \frac{1}{2} \beta \gamma \kappa & \text{ on } \Sigma_{0}. \label{SI:Cristini:pressurekappa}
\end{align}
Thus, we obtain
\begin{subequations}\label{SI:Cristini:newlimit}
\begin{align}
\div \bm{v} & = \mathcal{P} \sigma_{0} - \mathcal{A} \text{ in } \Omega_{T}, \\
\bm{v} & = - \nabla ( p - \chi_{\varphi} \sigma_{0}) \text{ in } \Omega_{T}, \\
\Laplace \sigma_{0} & = \begin{cases} \mathcal{C} \sigma_{0} & \text{ in } \Omega_{T}, \\
0 & \text{ in } \Omega_{H},
\end{cases} \\
\jump{\sigma_{0}} = 0, \quad (\nabla \sigma_{0}^{T} - D \nabla \sigma_{0}^{H}) \cdot \bm{\nu} & = 0 \text{ on } \Sigma_{0}, \\
p & = \frac{1}{2} \beta \gamma \kappa \text{ on } \Sigma_{0}, \label{SI:Cristini:newlimit:pressure} \\
- \nabla p \cdot \bm{\nu} + \chi_{\varphi}  \nabla \sigma_{0} \cdot \bm{\nu} & = \mathcal{V} \text{ on } \Sigma_{0},
\end{align}
\end{subequations}
which coincides with the sharp interface model (Eq. (79)-(81), (83)-(86)) of \cite{article:CristiniLiLowengrubWise09}.  We point out that the same sharp interface limit (\ref{SI:Cristini:newlimit}) can be recovered if we set $\lambda = \eps$ in (\ref{Model:Cristini}).  We introduce the parameter $\lambda$ in (\ref{Model:Cristini}) in order to study the effect of active transport on the linear stability of radial solutions to (\ref{SI:Cristini}), see Section \ref{sec:stabAna} below.

Let us also remark that the mobility $m(\varphi) = \frac{1}{2}(1+\varphi)^{2}$ is degenerate in the region $\{ \varphi = -1\}$, and thus the bulk equation for $\mu_{0}^{H}$ in $\Omega_{H}$ and the interfacial condition for $\nabla \mu_{0}^{H} \cdot \bm{\nu}$ on $\Sigma_{0}$ remain undetermined in (\ref{SI:Cristini}).  Furthermore, if $m(\varphi)$ is chosen to be degenerate in the bulk regions $\{ \varphi =  \pm 1 \}$, then we obtain from the outer expansion $(\ref{Model:Cristini:order})_{O}^{0}$ the following equations
\begin{align*}
0 = \mathcal{P} \sigma_{0} (\varphi_{0} + 1) - \mathcal{A} (\varphi_{0} + 1) \text{ in } \Omega_{T} \cup \Omega_{H}.
\end{align*}
In particular, we see that $\sigma_{0}^{T} = \frac{\mathcal{A}}{\mathcal{P}}$ is a constant in $\Omega_{T}$, which is inconsistent with $(\ref{Model:Cristini:sigma})_{O}^{0}$.  Hence, it is necessary that the mobility $m(\varphi)$ is not degenerate in the tumour region $\{ \varphi = 1 \}$.

\section{Linear Stability Analysis}\label{sec:stabAna}
Let us consider the sharp interface model (\ref{SI:Cristini}).  By sending the active transport parameter $\lambda$ to zero, we recover the sharp interface model (Eq. (79)-(81), (83)-(86)) of \cite{article:CristiniLiLowengrubWise09}.  In this section, we extend the linear stability analysis of \cite{article:CristiniLiLowengrubWise09,thesis:Li} to include the effects of active transport.  For the linear stability analysis of a one-phase model, we refer to \cite{article:CristiniLowengrubNie02,article:LiCristinitNieLowengrub07}.

\subsection{Radial solutions}
We now drop the index $0$ in (\ref{SI:Cristini}), and let $\Omega = B_{R}(0)$ denote the $d$-dimensional ball, $d = 2,3$, of radius $R$ centered at the origin.  We assume that the interface $\Sigma$ is a $(d-1)$-sphere of radius $q(t)$, partitioning the domain $\Omega$ into $\Omega_{T}$ and $\Omega_{H}$ as follows:
\begin{align*}
\Sigma = \pd B_{q(t)}, \quad \Omega_{T} = B_{q(t)}(0), \quad \Omega_{H} = B_{R}(0) \setminus \overline{B_{q(t)}(0)}.
\end{align*}
The outer unit normal $\bm{\nu}(\bm{p})$ at a point $\bm{p} \in \Sigma$ is given as
\begin{align}\label{defn:curvatureandnormal}
\bm{\nu}(\bm{p}) = \frac{\bm{p}}{\abs{\bm{p}}} = \frac{\bm{p}}{q(t)},
\end{align}
while the normal velocity $\mathcal{V}$ is given as
\begin{align}
\label{defn:normalvelo}
\mathcal{V} = \frac{\dd q}{\dt}.
\end{align}
The mean curvature $\kappa$ for a $(d-1)$-sphere radius $r_{0}$ is given by 
\begin{align}\label{meancurvature:sphere}
\kappa = \frac{d-1}{r_{0}},
\end{align}
where $d$ denotes the dimension.  Then, for radially symmetric solutions $\varphi(\abs{\bm{x}}) = \varphi(r)$, $\mu(\abs{\bm{x}}) = \mu(r)$, $\sigma(\abs{\bm{x}}) = \sigma(r)$, (\ref{SI:Cristini}) becomes
\begin{subequations}\label{SI:Cristini:radial}
\begin{align}
\mu_{T}'' + \frac{d-1}{r} \mu_{T}' & = \mathcal{A} - \mathcal{P} \sigma_{T} \text{ in } r < q(t), \\
\sigma'' + \frac{d-1}{r} \sigma' & = \begin{cases}
\mathcal{C} \sigma & \text{ in } r < q(t), \\
0 & \text{ in } r > q(t),
\end{cases} \\
\jump{\sigma} = 2\lambda, \quad 2 \mu_{T} & = \beta \gamma \frac{d-1}{q(t)} - \chi_{\varphi} (\sigma_{T} + \sigma_{H}) \text{ on } r = q(t),  \\
\sigma_{T}' = D \sigma_{H}', \quad -\frac{\dd q}{\dt} & = \mu_{T}' \text{ on } r = q(t) \label{SI:Cristini:radial:sigmajumpvelo}. 
\end{align}
\end{subequations}
We complete (\ref{SI:Cristini:radial}) with the following boundary conditions:
\begin{align}
\label{Radial:bdycond}
\sigma_{H}(r = R,t) = \sigma_{\infty}, \quad \sigma_{T}(r = 0, t) < \infty, \quad \mu_{T}(r = 0, t) < \infty,
\end{align}
where $\sigma_{\infty}$ denotes the concentration of a nutrient supply from the boundary $\pd \Omega$.  

Upon solving the differential equations and applying the interface and boundary conditions, we arrive at the following radial solutions:
\begin{subequations}\label{radialsoln}
\begin{align}
\sigma_{H}(r,t) & = \begin{cases}
\sigma_{\infty} + a_{2}(t) (\log(r) - \log(R)) & \text{ for } d = 2, \\
\sigma_{\infty} + a_{3}(t) \left ( \frac{1}{r} - \frac{1}{R} \right ) & \text{ for } d = 3, \\
\end{cases} \\
\sigma_{T}(r,t) & =  \begin{cases}
b_{2}(t) I_{0}(\Lambda r) & \text{ for } d = 2, \\
\displaystyle b_{3}(t) \frac{\sinh(\Lambda r)}{r} & \text{ for } d = 3, \\
\end{cases} \\
\mu_{T}(r, t) & = \begin{cases}
\displaystyle \frac{\mathcal{A}}{4} r^{2} - \frac{\mathcal{P}}{\mathcal{C}} b_{2}(t) I_{0}(\Lambda r) + c_{2}(t) & \text{ for } d = 2, \\[2ex]
\displaystyle \frac{\mathcal{A}}{6} r^{2} -\frac{\mathcal{P}}{\mathcal{C}} b_{3}(t) \frac{\sinh(\Lambda r)}{r} + c_{3}(t) & \text{ for } d = 3, \\
\end{cases}
\end{align}
\end{subequations}
where, for $\alpha \in \R$, $I_{\alpha}(x)$ denote the modified Bessel functions of the first kind :
\begin{align}
\label{defn:ModifiedBessel}
I_{\alpha}(x) = \sum_{k=0}^{\infty} \frac{1}{k! \Gamma(k + \alpha + 1)} \left ( \frac{x}{2} \right)^{2k + \alpha}.
\end{align}
Here, $\Gamma(\cdot)$ denotes the Gamma function.  Together with the modified Bessel functions of the second kind, $ K_{\alpha}(x) := \frac{\pi}{2} \frac{I_{-\alpha}(x) - I_{\alpha}(x)}{\sin(\alpha \pi)}$, the pairs $\{I_{\alpha}, K_{\alpha} \}$ are the two linearly independent solutions to the modified Bessel's equation:
\begin{align}\label{ModifiedBesselequ}
x^{2}\frac{\dd^{2}y}{\dx} + x \frac{\dd y}{\dx} - x^{2}y = \alpha^{2} y.
\end{align}
Moreover, for the case $\alpha = 0$, the modified Bessel functions $I_{0}(x)$, $K_{0}(x)$ satisfy the following properties
\begin{align}\label{Bessel0properties}
I_{0}(0) = 1, \quad \lim_{x \to 0} K_{0}(x) = +\infty, \quad \frac{\dd}{\dx} I_{0}(x) = I_{1}(x), \quad \int x I_{0}(x) \dx = x I_{1}(x).
\end{align}
Furthermore, the coefficients in (\ref{radialsoln}) are given as
\begin{subequations}\label{radialsolution:constantsIntegration}
\begin{align}
\Lambda^{2} & = \mathcal{C}, \\
a_{2}(t) & = \frac{q(t) \Lambda I_{1}(\Lambda q(t)) (\sigma_{\infty} + 2 \lambda)}{D I_{0}(\Lambda q(t)) - \Lambda q(t) \log(q(t)/R) I_{1}(\Lambda q(t))}  \; , \\
a_{3}(t) & = (\sigma_{\infty} + 2 \lambda) \frac{Rq(t) (1-q(t) \Lambda \coth(\Lambda q(t)))}{(R-q(t))(q(t) \Lambda \coth(\Lambda q(t)) - 1) + D R} \; , \\
b_{2}(t) & = \frac{D(\sigma_{\infty} + 2 \lambda)}{DI_{0}(\Lambda q(t)) - q(t) \Lambda \log (q(t)/R)  I_{1}(\Lambda q(t))} \; , \\
b_{3}(t) & = \frac{(\sigma_{\infty} + 2 \lambda)}{\sinh(\Lambda q(t))} \frac{DR q(t)}{(R-q(t))(q(t) \Lambda \coth(\Lambda q(t)) - 1) + D R} \; , \\
c_{2}(t) & = - \frac{\mathcal{A}}{4} q(t)^{2} + \frac{\beta \gamma}{2 q(t)} + \chi_{\varphi} \lambda + \left ( \frac{\mathcal{P}}{\Lambda} - \chi_{\varphi} \right ) b_{2}(t) I_{0}(\Lambda q(t)) \; , \\
c_{3}(t) & = - \frac{\mathcal{A}}{6} q(t)^{2} + \frac{\beta \gamma}{q(t)} + \chi_{\varphi} \lambda + \left ( \frac{\mathcal{P}}{\mathcal{C}} - \chi_{\varphi} \right ) b_{3}(t) \frac{\sinh(\Lambda q(t))}{q(t)} \; ,
\end{align}
\end{subequations}
and the differential equation satisfied by $q(t)$ is
\begin{equation}\label{radial:ODE}
\frac{\dd q}{\dt}  = \begin{cases}
\displaystyle  -\frac{\mathcal{A}}{2} q + \frac{\mathcal{P}}{\Lambda} b_{2}(t) I_{1}(\Lambda q)  & \text{ for } d = 2, \\[2ex]
\displaystyle -\frac{\mathcal{A}}{3} q + b_{3}(t)\frac{ \mathcal{P}}{\mathcal{C}} \left ( \frac{\Lambda \cosh(\Lambda q)}{q} - \frac{\sinh(\Lambda q)}{q^{2}} \right ) & \text{ for } d = 3.
\end{cases}
\end{equation}
We point out that, thanks to the boundary condition $\sigma_{T}(r = 0, t) < \infty$, the solution $\sigma_{T}$ does not contain any terms involving $K_{0}(\Lambda r)$ (in $d = 2$) and $\cosh(\Lambda r)/r$ (in $d = 3$).

\subsection{Perturbation of radial solutions}
We now consider a perturbation of a radially symmetric tumour, whose radius $w$ is given by
\begin{align}\label{defn:radius:perturb}
w(r, \theta, \phi, t) = q(t) + \delta(t) Z(\theta, \phi), \quad Z(\theta, \phi) = \begin{cases}
\cos(l \theta) & \text{ for } d = 2, \\
Y_{l,m}(\theta, \phi) & \text{ for } d = 3,
\end{cases}
\end{align}
where $q(t)$ is the radius of the unperturbed interface, $\delta(t)$ is a dimensionless perturbation size, $Y_{l,m}$ is a spherical harmonic with $l$ and $\theta$ denoting the polar wavenumber and angle, and $m$ and $\phi$ denoting the azimuthal wavenumber and angle, respectively.  We will denote the radial solutions in (\ref{radialsoln}) by $\sigma_{H}^{*}$, $\sigma_{T}^{*}$, and $\mu_{T}^{*}$, and consider 
\begin{subequations}
\begin{align}
 \sigma_{H}(r, \theta, \phi, t) & = \sigma_{H}^{*}(r,t) + U(r,t) \delta(t) Z(\theta, \phi), \\
 \sigma_{T}(r, \theta, \phi, t) & = \sigma_{T}^{*}(r,t) + V(r,t) \delta(t) Z(\theta, \phi), \\
 \mu_{T}(r, \theta, \phi, t)  & = \mu_{T}^{*}(r,t) + W(r,t) \delta(t) Z(\theta, \phi),
\end{align}
\end{subequations}
where we assume that $(\sigma_{H}, \sigma_{T}, \mu_{T})$ solve (\ref{SI:Cristini}).  Therefore, we get
\begin{subequations}\label{SI:perturbed}
\begin{alignat}{2}
\Laplace (\mu_{T}^{*} + W \delta Z)  & = \mathcal{A} - \mathcal{P} (\sigma_{T}^{*} + V \delta Z) && \text{ in } r < w, \label{SI:perturbed:muT} \\
\Laplace (\sigma_{T}^{*} + V \delta Z) & = \mathcal{C} (\sigma_{T}^{*} + V \delta Z) && \text{ in } r < w, \label{SI:perturbed:sigmaT} \\
\Laplace (\sigma_{H}^{*} + U \delta Z)  &= 0 && \text{ in } r > w, \label{SI:perturbed:sigmaH} \\
\sigma_{T}^{*} - \sigma_{H}^{*} + (V - U)\delta Z  &= 2\lambda && \text{ on } r = w, \label{SI:perturbed:sigmajump} \\
2 (\mu_{T}^{*} + W \delta Z) & = \beta \kappa \gamma - 2\chi_{\varphi} (\sigma_{T}^{*} + V \delta Z) + 2 \chi_{\varphi} \lambda && \text{ on } r = w, \label{SI:perturbed:kappa} \\
(\sigma_{T}^{*})_{r} + \delta \nabla (VZ) \cdot \bm{\nu} &=  D ((\sigma_{H}^{*})_{r} + \delta \nabla (UZ) \cdot \bm{\nu}) && \text{ on } r = w, \\
-\frac{\dd q}{\dt} - Z \frac{\dd \delta}{\dt} &  =  (\mu_{T}^{*})_{r} + \delta \nabla (WZ) \cdot \bm{\nu} && \text{ on } r = w.
\end{alignat}
\end{subequations}
Here, we used the more convenient form (\ref{SI:Cristini:kappa}) of $(\ref{SI:Cristini:jump})_{2}$.  Next, we linearise (\ref{SI:perturbed}) about the original interface $r = q$ to derive the equations satisfied by $U$, $V$, $W$ and $\delta$.  We introduce the Laplace--Beltrami operator on the $(d-1)$-sphere, for $d = 2, 3$:
\begin{align}
\mathcal{L}_{d} := 
\begin{cases} 
\displaystyle \frac{\pd^{2}}{\pd \theta^{2}} & \text{ for } d = 2, \\[2ex]
\displaystyle  \frac{\pd^{2}}{\pd \theta^{2}} + \cot(\theta) \frac{\pd}{\pd \theta} + \frac{1}{\sin(\theta)^{2}} \frac{\pd^{2}}{\pd \phi^{2}} & \text{ for } d = 3, 
\end{cases}
\end{align}
so that the Laplace operator can be decomposed into
\begin{align}
\Laplace f = f_{rr} + \frac{d-1}{r} f_{r} + \frac{1}{r^{2}} \mathcal{L}_{d} f.
\end{align}
Moreover, the function $Z(\theta, \phi)$ defined in (\ref{defn:radius:perturb}) satisfies
\begin{align}\label{defn:zetald}
\mathcal{L}_{d} Z(\theta, \phi) = \zeta_{l,d} Z(\theta, \phi), \quad \zeta_{l,d} =  \begin{cases}
-l^{2} & \text{ for } d = 2, \\
-l(l+1) & \text{ for } d = 3.
\end{cases}
\end{align}
From the bulk equation (\ref{SI:perturbed:muT}) we obtain
\begin{align*}
\Laplace \mu_{T}^{*} + \delta \Laplace (WZ) = \mathcal{A} - \mathcal{P} \sigma_{T}^{*} - \delta \mathcal{P} V Z,
\end{align*}
and so, using that $\Laplace \mu_{T}^{*} = \mathcal{A} - \mathcal{P} \sigma_{T}^{*}$, we deduce that
\begin{align*}
\mathcal{P} V + W_{rr} + \frac{d-1}{r} W_{r} + \frac{1}{r^{2}} \zeta_{l,d} W = 0 \text{ in } r < q.
\end{align*}
For the interface conditions, we employ Taylor's expansion and neglect terms of order $\mathcal{O}(\delta^{2})$.  For instance, from (\ref{SI:perturbed:sigmajump}), we see that 
\begin{align*}
2 \lambda = \sigma_{T}^{*}(q) + (\sigma_{T}^{*})_{r}(q) (w-q) - \sigma_{H}^{*}(q) - (\sigma_{H}^{*})_{r}(q) (w-q) + (V - U) \delta Z + \mathcal{O}(\delta^{2}).
\end{align*}
Then, by $(\ref{SI:Cristini:jump})_{1}$, $(\ref{SI:Cristini:nablamu})_{1}$, and (\ref{defn:radius:perturb}), we obtain
\begin{align*}
U(q,t) - V(q,t) = (\sigma_{T}^{*})_{r}(q) - (\sigma_{H}^{*})_{r}(q) = (D - 1) (\sigma_{H}^{*})_{r}(q) \text{ on } r = q.
\end{align*}
We use the following expansion for the mean curvature (compare with Eq. (4.12) of \cite{article:EscherMatioc13}, page 647 of \cite{article:CuiEscher08} and page 12 of \cite{article:FriedmanReitich01}, where instead of \eqref{meancurvature:sphere}, the mean curvature of a $(d-1)$-sphere radius $r_{0}$ is defined to be $\frac{1}{r_{0}}$):
\begin{align*}
\kappa (r = w) = \frac{d-1}{q} - \frac{d-1}{q^{2}} \delta \left ( 1 + \frac{\zeta_{l,d}}{d-1} \right ) Z(\theta, \phi) + \mathcal{O}(\delta^{2}),
\end{align*}
so that the linearisation of (\ref{SI:perturbed:kappa}) about $r = q$ is
\begin{align*}
(\mu_{T}^{*})_{r}(q) + W(q,t) = - \frac{\beta \gamma}{2} \frac{d-1}{q^{2}} \left ( 1 + \frac{\zeta_{l,d}}{d-1} \right ) - \chi_{\varphi} ((\sigma_{T}^{*})_{r}(q) + V(q,t)) \text{ on } r = q.
\end{align*}
Finally, by the relation $\nabla f (\abs{\bm{x}}) \cdot \bm{\nu} = f'(r)$ for $ \bm{x} \in \Sigma$, we have that
\begin{align*}
\nabla (VZ) \cdot \bm{\nu} \vert_{r=q} = \pd_{r} (V(r,t) Z(\theta, \phi)) \vert_{r=q} = V_{r}(q,t) Z(\theta, \phi),
\end{align*}
and so we obtain the following system for the perturbations $U, V, W$ and $\delta$ from linearising (\ref{SI:perturbed}) about the unperturbed interface $r = q$:
\begin{subequations}\label{Perturbation:system}
\begin{alignat}{2}
W_{rr} + \frac{d-1}{r} W_{r} + \frac{\zeta_{l,d}}{r^{2}} W  &= - \mathcal{P} V && \text{ in } r < q, \label{Perturbed:W} \\
V_{rr} + \frac{d-1}{r} V_{r} + \frac{\zeta_{l,d}}{r^{2}} V &= \mathcal{C} V && \text{ in } r < q, \label{Perturbed:V} \\
U_{rr} + \frac{d-1}{r} U_{r} + \frac{\zeta_{l,d}}{r^{2}} U &= 0 && \text{ in } r > q, \label{Perturbed:U} \\
U - V &= (D-1) (\sigma_{H}^{*})_{r}(q) && \text{ on } r = q, \label{Perturbed:sigmajump} \\
(\mu_{T}^{*} + \chi_{\varphi} \sigma_{T}^{*})_{r}(q) + W + \chi_{\varphi} V &= - \frac{\beta \gamma}{2} \frac{d-1}{q^{2}} \left ( 1 + \frac{\zeta_{l,d}}{d-1} \right )  && \text{ on } r = q, \label{Perturbed:kappa} \\
(\sigma_{T}^{*} - D \sigma_{H}^{*})_{rr}(q) &= DU_{r} - V_{r} && \text{ on } r = q, \label{Perturbed:nablasigma} \\
\frac{\dd \delta}{\dt} & =  -(\mu_{T}^{*})_{rr}(q) \delta - \delta W_{r} && \text{ on } r = q. \label{Perturbed:velo}
\end{alignat}
\end{subequations}
We complete (\ref{Perturbation:system}) with the following boundary conditions:
\begin{align}
\label{Perturbation:system:bdy}
W(r = 0,t) < \infty, \quad V(r = 0, t) < \infty, \quad U(r = R, t) = 0.
\end{align}

\subsection{Solutions to the perturbed system}
Recalling the definition of $\zeta_{l,d}$ in (\ref{defn:zetald}), we see that the general solution for (\ref{Perturbed:U}) is
\begin{align}\label{generalsoln:U}
U(r,t) = \begin{cases}
F_{0}(t) r^{l} + F_{1}(t) r^{-l} & \text{ for } d = 2, \\
F_{0}(t) r^{l} + F_{1}(t) r^{-l-1} & \text{ for } d = 3.
\end{cases}
\end{align}
We observe that the ODE (\ref{Perturbed:V}) in $d = 2$ is a scaled modified Bessel's equation (see (\ref{ModifiedBesselequ})), while (\ref{Perturbed:V}) in $d = 3$ is a scaled modified spherical Bessel's equation.  Due to the boundary condition (\ref{Perturbation:system:bdy}), we see that the general solution to (\ref{Perturbed:V}) is given by
\begin{align}\label{generalsoln:V}
V(r,t) = \begin{cases}
F_{2}(t)I_{l}(\Lambda r) & \text{ for } d = 2, \\
F_{2}(t) i_{l}(\Lambda r) & \text{ for } d = 3,
\end{cases}
\end{align}
where $I_{l}(x)$ is the modified Bessel function of the first kind, defined in (\ref{defn:ModifiedBessel}), while $i_{l}(x)$ is the modified spherical Bessel function of first kind that satisfies
\begin{align*}
\frac{\dd^{2}}{\dd r^{2}} i_{l}(r) + \frac{2}{r} \frac{\dd}{\dd r} i_{l}(r) - \frac{l(l+1)}{r^{2}} i_{l}(r) = i_{l}(r), \quad i_{l}(0) < \infty \quad \forall l \geq 0.
\end{align*}
Again, due to the boundary condition $(\ref{Perturbation:system:bdy})_{2}$, $V(r,t)$ does not contain any terms involving the modified spherical Bessel function of the second kind.

For (\ref{Perturbed:W}), we see that $W$ is a sum of the solution to the homogeneous equation (\ref{Perturbed:U}) and the particular solution $\frac{\mathcal{P}}{\mathcal{C}} V$.  Due to the boundary condition (\ref{Perturbation:system:bdy}) for $W$, we find that the general solution to (\ref{Perturbed:W}) is
\begin{align}\label{generalsoln:W}
W(r,t) = F_{3}(t) r^{l} - \frac{\mathcal{P}}{\mathcal{C}} V(r,t).
\end{align}
With these solutions (\ref{generalsoln:U}), (\ref{generalsoln:V}), and (\ref{generalsoln:W}), we use the relations (\ref{Perturbed:sigmajump}), (\ref{Perturbed:kappa}), (\ref{Perturbed:nablasigma}) in order to simplify the resulting differential equation (\ref{Perturbed:velo}) for $\delta$.  Let
\begin{align*}
Q(\Lambda, q) := \frac{\Lambda q \cosh(\Lambda q) - \sinh(\Lambda q)}{q^{2}} \; ,
\end{align*} 
then from (\ref{Perturbed:sigmajump}), (\ref{Perturbed:kappa}), and (\ref{Perturbed:nablasigma}) we obtain the following relations:  
\begin{subequations}\label{Perturbation:ConstRelations}
\allowdisplaybreaks[3]
\begin{align}
\begin{cases}
F_{0} q^{l} + F_{1} q^{-l} - F_{2} I_{l}(\Lambda q) = (D-1) \frac{a_{2}}{q} & \text{ for } d = 2, \\
& \\
F_{0} q^{l} + F_{1} q^{-l-1} - F_{2} i_{l}(\Lambda q) = (1-D) \frac{a_{3}}{q^{2}} & \text{ for } d = 3, \\
\end{cases} \label{relation:perturbation:sigmajump} \\
\notag \\
\begin{cases}\left ( \chi_{\varphi} - \frac{\mathcal{P}}{\mathcal{C}} \right ) (b_{2} \Lambda I_{1}(\Lambda q) + F_{2} I_{l}(\Lambda q)) + \frac{\mathcal{A}}{2}q + F_{3} q^{l} = \frac{\beta \gamma}{2} \frac{l^{2} - 1}{q^{2}} & \text{ for } d = 2, \\
& \\
\left ( \chi_{\varphi} - \frac{\mathcal{P}}{\mathcal{C}} \right ) (b_{3} Q(\Lambda, q) + F_{2} i_{l}(\Lambda q)) + \frac{\mathcal{A}}{3}q + F_{3} q^{l} = \frac{\beta \gamma}{2} \frac{(l+2)(l-1)}{q^{2}} & \text{ for } d = 3, \\
\end{cases} \\
\notag \\
\begin{cases}
\mathcal{C} b_{2} I_{0}(\Lambda q) = D (l F_{0} q^{l-1} - F_{1} l q^{-l-1}) - F_{2} \Lambda I_{l}'(\Lambda q) & \text{ for } d = 2, \\
& \\
\mathcal{C} b_{3} \frac{\sinh(\Lambda q)}{q} = D(l F_{0} q^{l-1} - F_{1} (l+1) q^{-l-2}) - F_{2} \Lambda i_{l}'(\Lambda q)  & \text{ for } d = 3, \\
\end{cases}
\label{relation:perturbation:sigmanabla}
\end{align}
\end{subequations}
where we have used that
\begin{align}
(\mu_{T}^{*} + \chi_{\varphi} \sigma_{T}^{*})_{r}(q) = \begin{cases}
\frac{\mathcal{A}}{2} q + \left ( \chi_{\varphi} - \frac{\mathcal{P}}{\mathcal{C}} \right ) \Lambda b_{2} I_{1}(\Lambda q) & \text{ for } d = 2, \\
\frac{\mathcal{A}}{3} q + \left ( \chi_{\varphi} - \frac{\mathcal{P}}{\mathcal{C}} \right ) b_{3}(t) Q(\Lambda, q) & \text{ for } d = 3, \\
\end{cases}
\end{align}
and by $(\ref{SI:Cristini:radial:sigmajumpvelo})_{1}$,
\begin{align}
(\sigma_{T}^{*} - D \sigma_{H}^{*})''(q) = \mathcal{C} \sigma_{T}^{*}(q) - \frac{d-1}{q} (\sigma_{T}^{*})'(q) + \frac{d-1}{q} D (\sigma_{H}^{*})'(q) = \mathcal{C} \sigma_{T}^{*}(q).
\end{align}
Also, from $(\ref{SI:Cristini:radial:sigmajumpvelo})_{1}$, we observe that the following relations hold
\begin{align}\label{sigmaderivative:relation}
(\sigma_{T}^{*})'(q) = D (\sigma_{H}^{*})'(q)  \Rightarrow \begin{cases} b_{2} \Lambda I_{1}(\Lambda q) = D \frac{a_{2}}{q} & \text{ for } d = 2, \\
b_{3} Q(\Lambda, q) = -D \frac{a_{3}}{q^{2}} & \text{ for } d = 3. \\
\end{cases}
\end{align}
Together with the relation
\begin{align*}
(\mu_{T}^{*})''(q) + \mathcal{P} \sigma_{T}^{*}(q) & = \mathcal{A} - \frac{d-1}{q} (\mu_{T}^{*})'(q) \\
& = \mathcal{A} - \frac{d-1}{q} \left ( \frac{\mathcal{A}}{d} q - \frac{\mathcal{P}}{\mathcal{C}} (\sigma_{T}^{*})'(q) \right )  = \begin{cases}
\frac{\mathcal{A}}{2} + \frac{\mathcal{P}}{\mathcal{C}} D \frac{a_{2}}{q^{2}} & \text{ for } d = 2, \\
\frac{\mathcal{A}}{3} - \frac{2 \mathcal{P}}{\mathcal{C}} D \frac{a_{3}}{q^{3}} & \text{ for } d = 3, \\
\end{cases}
\end{align*}
and the relations (\ref{Perturbation:ConstRelations}), we can simplify (\ref{Perturbed:velo}) in order to obtain the following differential equation for the perturbation size $\delta$:
\begin{equation}\label{delta:3d:ODE}
\begin{aligned}
\frac{1}{\delta} \frac{\dd \delta}{\dt} & =\frac{\mathcal{A}}{3} (l-1) - \frac{a_{3}}{q^{3}} \left ( l \chi_{\varphi} - (l + 2 D)\frac{\mathcal{P}}{\mathcal{C}}  \right )  - \beta \gamma \frac{l(l+2)(l-1)}{2 q^{3}} \\
& + F_{0} q^{l-1} \left ( l\chi_{\varphi} + l(D - 1) \frac{\mathcal{P}}{\mathcal{C}} \right ) + \frac{F_{1}}{q^{l+2}} \left (l \chi_{\varphi} - (l + lD + D) \frac{\mathcal{P}}{\mathcal{C}} \right ) \text{ for } d = 3,
\end{aligned}
\end{equation}
and
\begin{equation}\label{delta:2d:ODE}
\begin{aligned}
\frac{1}{\delta} \frac{\dd \delta}{\dt} & = \frac{\mathcal{A}}{2} (l-1) + \frac{a_{2}}{q^{2}} \left ( l \chi_{\varphi} - (l + D) \frac{\mathcal{P}}{\mathcal{C}} \right ) - \beta \gamma \frac{l(l^{2} - 1)}{2 q^{3}} \\
& + F_{0} q^{l-1} \left ( l \chi_{\varphi} + l(D-1) \frac{\mathcal{P}}{\mathcal{C}} \right ) + \frac{F_{1}}{q^{l+1}} \left ( l \chi_{\varphi} -(l +l D)\frac{\mathcal{P}}{\mathcal{C}} \right ) \text{ for } d = 2.
\end{aligned}
\end{equation}
Consequently, using (\ref{radial:ODE}), and (\ref{sigmaderivative:relation}), we obtain the following differential equations for the shape perturbation $\frac{\delta}{q}$:
\begin{equation}\label{shapeperturbation:ODE}
\begin{aligned}
& \; \frac{q}{\delta} \frac{\dd}{\dt} \left ( \frac{\delta}{q} \right ) = \frac{1}{\delta} \frac{\dd \delta}{\dt} - \frac{1}{q} \frac{\dd q}{\dt} \\
= & \;  \begin{cases}
\begin{aligned}
&l \frac{\mathcal{A}}{2}  + \frac{a_{2}}{q^{2}} \left ( l \chi_{\varphi} - (l + 2D) \frac{\mathcal{P}}{\mathcal{C}} \right ) - \beta \gamma \frac{l(l^{2} - 1)}{2 q^{3}} \\
&+ F_{0} q^{l-1} \left ( l \chi_{\varphi} + l(D-1) \frac{\mathcal{P}}{\mathcal{C}} \right ) + \frac{F_{1}}{q^{l+1}} \left ( l \chi_{\varphi} -(l +l D)\frac{\mathcal{P}}{\mathcal{C}} \right )\end{aligned}
 &  \text{ for } d = 2, \\
 & \\
\begin{aligned}
& l \frac{\mathcal{A}}{3}  - \frac{a_{3}}{q^{3}} \left ( l \chi_{\varphi} - (l + 3 D)\frac{\mathcal{P}}{\mathcal{C}}  \right )  - \beta \gamma \frac{l(l+2)(l-1)}{2 q^{3}} \\
& + F_{0} q^{l-1} \left ( l\chi_{\varphi} + l(D - 1) \frac{\mathcal{P}}{\mathcal{C}} \right ) + \frac{F_{1}}{q^{l+2}} \left (l \chi_{\varphi} - (l + lD + D) \frac{\mathcal{P}}{\mathcal{C}} \right )
\end{aligned}
 & \text{ for } d = 3.
\end{cases}
\end{aligned}
\end{equation}
Finally, we mention that the time-dependent constants $F_{0}$ and $F_{1}$ can be computed as follows:  Due to $(\ref{Perturbation:system:bdy})_{3}$, we have
\begin{align}
F_{0} = \begin{cases}
-F_{1} R^{-2l} & \text{ for } d = 2, \\
-F_{1} R^{-2l-1} & \text{ for } d = 3.
\end{cases}
\end{align}
Moreover, by (\ref{relation:perturbation:sigmajump}) and (\ref{relation:perturbation:sigmanabla}), we obtain
\begin{subequations}\label{K1:equation}
\begin{align}
\notag & \mathcal{C} b_{2} I_{0}(\Lambda q) + \Lambda \frac{I_{l}'(\Lambda q)}{I_{l}(\Lambda q)} (1-D) \frac{a_{2}}{q} \\
= & \; - F_{1} \left ( \frac{D l q^{l-1}}{R^{2l}} +\frac{Dl}{q^{l+1}} + \Lambda \frac{I_{l}'(\Lambda q)}{I_{l}(\Lambda q)} \left ( \frac{1}{q^{l}} - \frac{q^{l}}{R^{2l}} \right ) \right ) \text{ for } d = 2, \\
\notag & \;  \frac{\mathcal{C} b_{3}\sinh(\Lambda q)}{q} + \Lambda \frac{i_{l}'(\Lambda q)}{i_{l}(\Lambda q)} (D-1) \frac{a_{3}}{q^{2}} \\
= & \;  - F_{1} \left ( \frac{l D q^{l-1}}{R^{2l+1}} + \frac{(l+1)D}{q^{l+2}} + \Lambda \frac{i_{l}'(\Lambda q)}{i_{l}(\Lambda q)} \left ( \frac{1}{q^{l+1}} - \frac{q^{l}}{R^{2l+1}} \right ) \right ) \text{ for } d = 3,
\end{align}
\end{subequations}
respectively.

We observe that the active transport parameter $\lambda$ enters into the radial solutions (\ref{radialsoln}), the differential equations (\ref{radial:ODE}), (\ref{delta:3d:ODE}), (\ref{delta:2d:ODE}), and (\ref{shapeperturbation:ODE}) only via the time-dependent constants $a_{2}, a_{3}, b_{2}, b_{3}, c_{2}$ and $c_{3}$.

\subsection{Effect of active transport on linear stability}
We now investigate the effect of active transport on the linear stability of the system.  To compare with \cite{article:CristiniLiLowengrubWise09}, we consider the choices
\begin{align*}
\mathcal{C} = 1, \quad \Lambda = 1, \quad \sigma_{\infty} = 1,
\end{align*} 
and neglect $F_{0}$ in (\ref{shapeperturbation:ODE}).  This implies that \eqref{K1:equation} becomes
\begin{equation*}
\begin{alignedat}{3}
\mathcal{C} b_{2} I_{0}(\Lambda q) + \Lambda \frac{I_{l}'(\Lambda q)}{I_{l}(\Lambda q)} (1-D) \frac{a_{2}}{q} & = - F_{1} \left (\frac{Dl}{q^{l+1}} + \Lambda \frac{I_{l}'(\Lambda q)}{I_{l}(\Lambda q)}  \frac{1}{q^{l}}  \right ) && \text{ for } d = 2, \\
\frac{\mathcal{C} b_{3}\sinh(\Lambda q)}{q} + \Lambda \frac{i_{l}'(\Lambda q)}{i_{l}(\Lambda q)} (D-1) \frac{a_{3}}{q^{2}} & = - F_{1} \left ( \frac{(l+1)D}{q^{l+2}} + \Lambda \frac{i_{l}'(\Lambda q)}{i_{l}(\Lambda q)} \frac{1}{q^{l+1}}  \right ) && \text{ for } d = 3.
\end{alignedat}
\end{equation*}
We define
\begin{equation*}
\begin{alignedat}{3}
\overline{a_{2}} & = \frac{q I_{1}(q)}{D I_{0}(q) - q \log(q/R) I_{1}(q)} \; , && \quad \overline{b_{2}} && = \frac{D}{DI_{0}(q) - q \log(q/R) I_{1}(q)} \; ,  \\
\overline{a_{3}} & = \frac{Rq (1-q \coth(q))}{(R-q)(q \coth(q) - 1) + D R} \; , &&\quad \overline{b_{3}} && = \frac{D R q}{(R-q)(q \cosh(q) - \sinh(q)) + DR \sinh(q)} \; ,
\end{alignedat}
\end{equation*}
so that $a_{2} = \overline{a_{2}} (1+2 \lambda)$, $b_{2} = \overline{b_{2}} (1 + 2 \lambda)$, $a_{3} = \overline{a_{3}} (1 + 2 \lambda)$ and $b_{3} = \overline{b_{3}}(1+ 2 \lambda)$, where $a_{2}$, $a_{3}$, $b_{2}$ and $b_{3}$ are as defined in \eqref{radialsolution:constantsIntegration}.  A short computation yields that
\begin{align*}
C_{2} := \frac{D \overline{a_{2}}}{q} = \frac{D I_{1}(q)/I_{0}(q)}{D - q \log(q/R) I_{1}(q)/I_{0}(q)} \; , \; C_{3} := \frac{-D \overline{a_{3}}}{q^{2}} = \frac{D (\coth(q) - \frac{1}{q})}{D + q \frac{R-q}{R} (\coth(q) - \frac{1}{q})} \; .
\end{align*}
Using the following relations for the modified Bessel functions and modified spherical Bessel functions of the first kind:
\begin{align*}
I_{l}'(z) = \frac{l}{z} I_{l}(z) + I_{l+1}(z), \quad i_{l}'(z) = \frac{l}{z} i_{l}(z) + i_{l+1}(z), \quad i_{l}(z) = \sqrt{\frac{\pi}{2 z}} I_{l+\frac{1}{2}}(z),
\end{align*}
and the relations
\begin{align*}
\overline{b_{2}} I_{0}(q) & = \frac{D I_{0}(q)}{D I_{0}(q) - q \log(q/R) I_{1}(q)} = \frac{I_{0}(q)}{I_{1}(q)} \frac{D I_{1}(q)}{D I_{0} - q \log(q/R) I_{1}(q)} = C_{2} \frac{I_{0}(q)}{I_{1}(q)}, \\
\overline{b_{3}} \frac{\sinh(q)}{q} & = \frac{DR}{(R-q) (q\coth(q) - 1) + DR}  = \frac{C_{3}}{\coth(q) - \frac{1}{q}},
\end{align*}
we find that
\begin{equation}\label{K1computed}
\begin{aligned}
F_{1} = \begin{cases}
\displaystyle - (1+2 \lambda) q^{l+1} C_{2} \frac{\left ( \frac{I_{0}(q)}{I_{1}(q)} + \frac{1-D}{D} \left ( \frac{l}{q} + \frac{I_{l+1}(q)}{I_{l}(q)} \right ) \right )}{\left ( Dl + l + q\frac{I_{l+1}(q)}{I_{l}(q)} \right )} & \text{ for } d = 2, \\
& \\
\displaystyle -(1 + 2 \lambda) q^{l+2} C_{3}\frac{\left ( \frac{1}{\coth(q) - 1/q} + \frac{1-D}{D} \left ( \frac{I_{l+3/2}(q)}{I_{l+1/2}(q)} + \frac{l}{q} \right ) \right ) }{ \left ( (l+1)D + l + q \frac{I_{l+3/2}(q)}{I_{l+1/2}(q)} \right )} & \text{ for } d = 3.
\end{cases}
\end{aligned}
\end{equation}
Substituting $F_{0} = 0$ and $\lambda = 0$ in (\ref{K1computed}) and (\ref{shapeperturbation:ODE}), we obtain the differential equation for the shape perturbation as derived in Eq. (89) of \cite{article:CristiniLiLowengrubWise09} with the notation $\tilde{\mathcal{G}}^{-1} :=\frac{1}{2} \beta \gamma$.  

Next, we find, for given $\mathcal{P}$, $D$, $\chi_{\varphi}$, and $\beta$, a critical value $\mathcal{A}_{c}$ such that $\frac{\dd}{\dt} \frac{\delta}{q} = 0$, i.e., the shape perturbation $(\frac{\delta}{q})$ is a constant.  This critical value $\mathcal{A}_{c}$ is given by
\begin{align*}
\mathcal{A}_{c} & = \beta \gamma \frac{(l^{2}-1)}{q^{3}} + (1 + 2 \lambda)2 C_{2} \frac{ (1 + \frac{2D}{l} ) \mathcal{P} - \chi_{\varphi}}{Dq} \\[2ex]
& + (1 + 2 \lambda)2C_{2} ( \chi_{\varphi} - (1+D) \mathcal{P}) \frac{\left ( \frac{I_{0}(q)}{I_{1}(q)} + \frac{1-D}{D} \left ( \frac{l}{q} + \frac{I_{l+1}(q)}{I_{l}(q)} \right ) \right )}{\left ( Dl + l + q\frac{I_{l+1}(q)}{I_{l}(q)} \right )} \text{ for } d = 2,
\end{align*}
and
\begin{align*}
\mathcal{A}_{c} & =  \beta \gamma \frac{3(l+2)(l-1)}{2 q^{3}} + (1 + 2\lambda) 3 C_{3} \frac{ (1 + \frac{3D}{l}) \mathcal{P} - \chi_{\varphi}}{Dq}  \\[2ex]
& + (1+2 \lambda)3 C_{3} (\chi_{\varphi} - (1 + D - \tfrac{D}{l}) \mathcal{P}) \frac{\left ( \frac{1}{\coth(q) - 1/q} + \frac{1-D}{D} \left ( \frac{I_{l+3/2}(q)}{I_{l+1/2}(q)} + \frac{l}{q} \right ) \right ) }{ \left ( (l+1)D + l + q \frac{I_{l+3/2}(q)}{I_{l+1/2}(q)} \right )} \text{ for } d = 3.
\end{align*}
We point out that, when $\lambda = 0$, the expression for $\mathcal{A}_{c}$ coincides with Eq. (90) of \cite{article:CristiniLiLowengrubWise09} with $\tilde{\mathcal{G}}^{-1} := \mathcal{G}^{-1}\tau = \frac{1}{2} \beta \gamma$.  We now look at $\mathcal{A}_{c}$ as a function of $q$ for the following parameter values:
\begin{align*}
\tilde{\mathcal{G}}^{-1} = \frac{1}{2} \beta \gamma = 0.05, \quad \mathcal{P} = 0.1, \quad D = 1, \quad l = 2, \quad R = 13.
\end{align*}
With these choices, we obtain 
\begin{equation}\label{criticalA:specificchoice}
\begin{aligned}
\mathcal{A}_{c} = \begin{cases}
\displaystyle \frac{0.3}{q^{3}} + (1+2\lambda)(\chi_{\varphi} - 0.2) 2 C_{2} \left (X - \frac{1}{q} \right )  & \text{ for } d = 2, \\[2ex]
\displaystyle \frac{0.6}{q^{3}} + 3 C_{3} (1+2 \lambda) \left ( \chi_{\varphi} \left ( Y - \frac{1}{q} \right ) + \left ( \frac{1}{4q}  - 0.15 Y \right ) \right ) & \text{ for } d = 3,
\end{cases}
\end{aligned}
\end{equation}
where 
\begin{alignat*}{3}
C_{2} & = \frac{I_{1}(q)/I_{0}(q)}{1 - q \log(q/R) I_{1}(q)/I_{0}(q)} \; , && \quad C_{3} && = \frac{ (\coth(q) - \frac{1}{q})}{1 + q \frac{R-q}{R} (\coth(q) - \frac{1}{q})} \; , \\
X & = \frac{1}{4 + q\frac{I_{3}(q)}{I_{2}(q)}}\frac{I_{0}(q)}{I_{1}(q)} \; , && \quad Y && = \frac{1}{5 + q \frac{I_{7/5}(q)}{I_{5/2}(q)}}\frac{1}{\coth(q) - 1/q}\; .
\end{alignat*}
Numerically, we find that $C_{2}$, $C_{3}$, $X$ and $Y$ are positive for $q \in (0, 13]$.  Moreover, 
\begin{align}\label{signXY}
X - \frac{1}{q} < 0 , \quad Y - \frac{1}{q} < 0, \quad \frac{1}{4q} - 0.15 Y > 0 \quad \forall q \in (0,13].
\end{align}
\begin{figure}[h]
\centering
\subfigure[\label{fig:Ac2dlambda}]{\includegraphics[width = 0.48\textwidth]{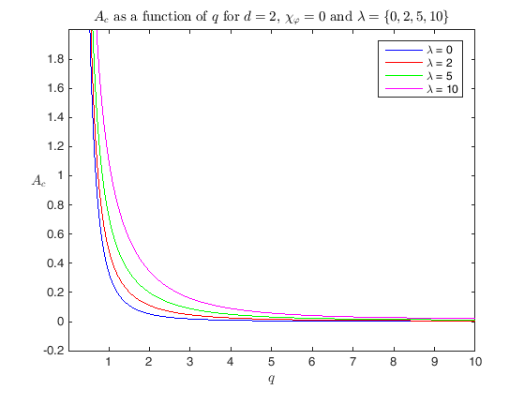}} 
\subfigure[\label{fig:Ac3dlambda}]{\includegraphics[width = 0.48\textwidth]{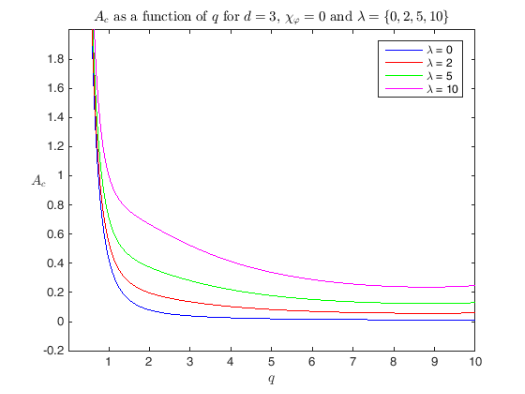}}  \\
\subfigure[\label{fig:Ac2d:varychi}]{\includegraphics[width = 0.48\textwidth]{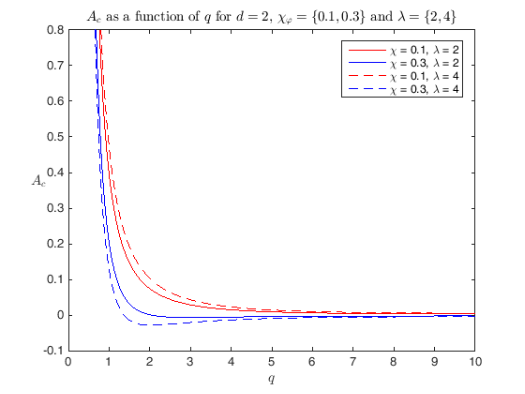}} 
\subfigure[\label{fig:Ac3d:varychi}]{\includegraphics[width = 0.48\textwidth]{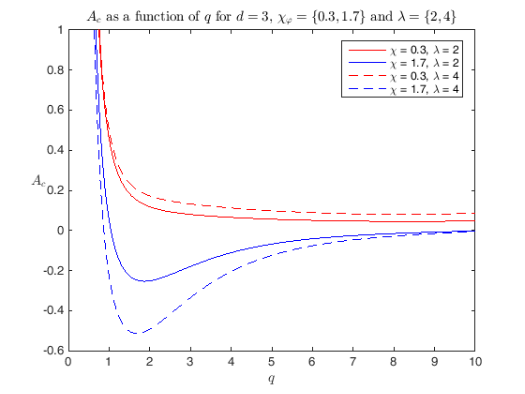}} 
\caption{Effects of $\lambda$ on the critical apoptosis parameter $A_{c}$ as a function of the unperturbed radius $q$ in 2d and 3d with $\beta \gamma = 0.1$, $\mathcal{P} = 0.1$, $D = 1$, $l = 2$, $R = 13$.}
\label{fig:CriticalAc}
\end{figure}
We note that $\mathcal{A}$ is the apoptosis parameter and $\mathcal{A}_{c}$ divides the phase portrait into regions of stable growth for low apoptosis (the region $\mathcal{A} < \mathcal{A}_{c}$) and regions of unstable growth for high apoptosis (the region $\mathcal{A} > \mathcal{A}_{c}$) for a given mode $l$.  Thus, from (\ref{criticalA:specificchoice}) and (\ref{signXY}), we observed the following:
\begin{enumerate}
\item In the absence of chemotaxis, $\chi_{\varphi} = 0$, increasing $\lambda$ will increase the value of $\mathcal{A}_{c}$.  From Figures \ref{fig:Ac2dlambda} and \ref{fig:Ac3dlambda}, the curves are pushed upwards, and so the region of stable growth for low apoptosis is enlarged.  In particular, active transport has a stabilising effect on the perturbations in the absence of chemotaxis.  
\item In dimension $d = 2$, while $\chi_{\varphi} < 0.2$, active transport has a stabilising effect on the perturbations.  When $\chi_{\varphi} > 0.2$, the perturbations are now amplified by the presence of active transport.  In Figure \ref{fig:Ac2d:varychi}, we see that, as $\lambda$ increases, the curves are pushed up for $\chi_{\varphi} = 0.1$, while the curves are pulled down for $\chi_{\varphi} = 0.3$.  Similarly, in dimension $d = 3$, we find that
\begin{align*}
\frac{0.25/q - 0.15Y(q)}{1/q - Y(q)} \in (0.400, 1.459) \text{ for } q \in [0.01, 13],
\end{align*}
and from Figure \ref{fig:Ac3d:varychi}, we see that, as $\lambda$ increases, the curves are pushed up for $\chi_{\varphi} = 0.3$, while the curves are pulled down for $\chi_{\varphi} = 1.7$.
\end{enumerate}

\section{Numerical Computations}\label{sec:Numerics}
In this section we first derive a finite element approximation of (\ref{Model:Cristini}) and then we 
display some numerical results obtained using this approximation. We concentrate on (\ref{Model:Cristini}), however approximations of other variations of the model follow in a natural way.  In the approximation we take $\Psi(\varphi)$ to be the double obstacle potential given in (\ref{defn:DoubleObstacle}).  This choice of potential leads to (\ref{Model:Cristini:mu}) taking the form of a variational inequality (\ref{Asym:mu:variationalinequ}).  

\subsection*{Finite element approximation}
Let $\mathcal{T}$ be a 
regular triangulation of $\Omega$ into disjoint open simplices, associated with $\mathcal{T}$ is the piecewise linear finite element space
\begin{align*}
S_{h} :=  \left \{ \varphi \in C^{0}(\overline\Omega) 
\Big| \,  \varphi_{|_{T}} \in P_{1}(T) \; \forall ~T \in \mathcal{T} \right \} \subset H^1(\Omega),
\end{align*}
where we denote by $P_{1}(T)$ the set of all affine linear functions on $T$.  We now introduce a finite element approximation of (\ref{Model:Cristini}) in which we have taken homogeneous Neumann boundary conditions for $\varphi$ and $\mu$, and the Dirichlet boundary condition $\sigma = \sigma_B\in \mathbb{R}$ on $\partial\Omega$:  Find 
\begin{align*}
\varphi_{h}^{n} \in K_{h} := \{ \chi \in S_{h} |~ \abs{\chi} \leq 1 \}, \quad \mu_{h}^{n} \in S_{h}, \quad
\sigma_{h}^{n} \in S_{h}^{B} := \{ \chi \in S_{h} |~\chi = \sigma_{B} \text{ on } \pd \Omega \}
\end{align*}
such that for all $\eta_{h} \in S_{h}$, $\zeta_{h} \in K_{h}$ and $\chi_{h} \in S_{h}^{0} := \{ \chi \in S_{h} |~\chi = 0 \text{ on } \pd \Omega\}$,
\begin{subequations}\label{FE:Cristini}
\begin{align}
\frac{1}{\tau} ( \varphi_{h}^{n} - \varphi_{h}^{n-1}, \eta_{h})_{h} + (m(\varphi_{h}^{n-1}) \nabla \mu_{h}^{n}, \nabla \eta_{h})_{h} & = ((\mathcal{P} \sigma_{h}^{n-1} - \mathcal{A})(\varphi_{h}^{n} + 1), \eta_{h})_{h}, \label{FE:Cristini:order} \\
\left ( \mu_{h}^{n} + \frac{\beta}{\eps} \varphi_{h}^{n-1} + \chi_{\varphi} \sigma_{h}^{n-1}, \zeta_{h} - \varphi_{h}^{n} \right )_{h} & \leq \beta \eps (\nabla \varphi_{h}^{n}, \nabla (\zeta_{h} - \varphi_{h}^{n})), \label{FE:Cristini:mu} \\
(\mathcal{D}(\varphi_{h}^{n}) \nabla \sigma_{h}^{n}, \nabla \chi_{h})_{h} - \lambda (\mathcal{D}(\varphi_{h}^{n}) \nabla \varphi_{h}^{n}, \nabla \chi_{h})_{h} &= -\frac{1}{2} \mathcal{C} (\sigma_{h}^{n}(\varphi_{h}^{n} + 1), \chi_{h})_{h}, \label{FE:Cristini:sigma}
\end{align}
\end{subequations}
where $m(\varphi) = \frac{1}{2}(1 + \varphi)^{2}$, $\tau$ denotes the time step, $(\eta_{1},\eta_{2})$ denotes the $L^{2}$ inner product and $(\eta_{1},\eta_{2})_{h} := \int_{\Omega}\pi_{h}(\eta_{1}(x) \eta_{2}(x)) \dx$ where on each triangle $\pi_{h}$ is taken to be an affine interpolation of the values of $\eta_{1} \eta_{2}$ at the nodes of the triangle. 

We note that since the interfacial thickness is proportional to $\eps$ in order to resolve the interfacial layer we need to choose $h \ll \eps$, see \cite{DDE} for details. Away from the interface $h$ can be chosen larger and hence adaptivity in space can heavily speed up computations. In fact we use the finite element toolbox Alberta 2.0, see \cite{alberta}, for adaptivity and we implemented the same mesh refinement strategy as in \cite{BNS}, i.e., a fine mesh is constructed where $\abs{\varphi_{h}^{n-1}} < 1$ with a coarser mesh present in the bulk regions $\abs{\varphi_{h}^{n-1}} = 1$.

We begin our numerical results by following the authors in \cite{article:HilhorstKaampmannNguyenZee15} 
in comparing solutions obtained from a simplified form of the diffuse interface model 
with exact solutions to a sharp interface limit model. 

\subsection{Comparison with a sharp interface limit solution}
\begin{figure}[!h]
\centering
\subfigure{\includegraphics[width = 0.28\textwidth]{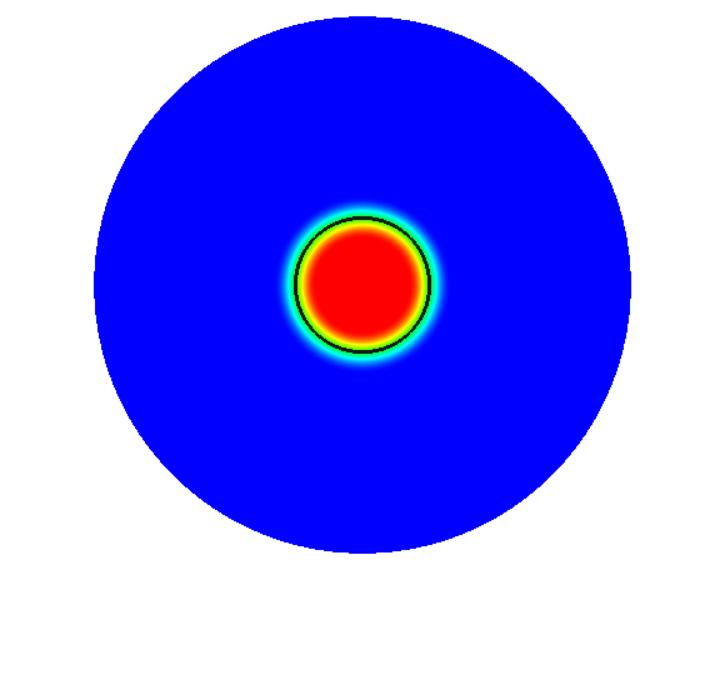}} 
\subfigure{\includegraphics[width = 0.28\textwidth]{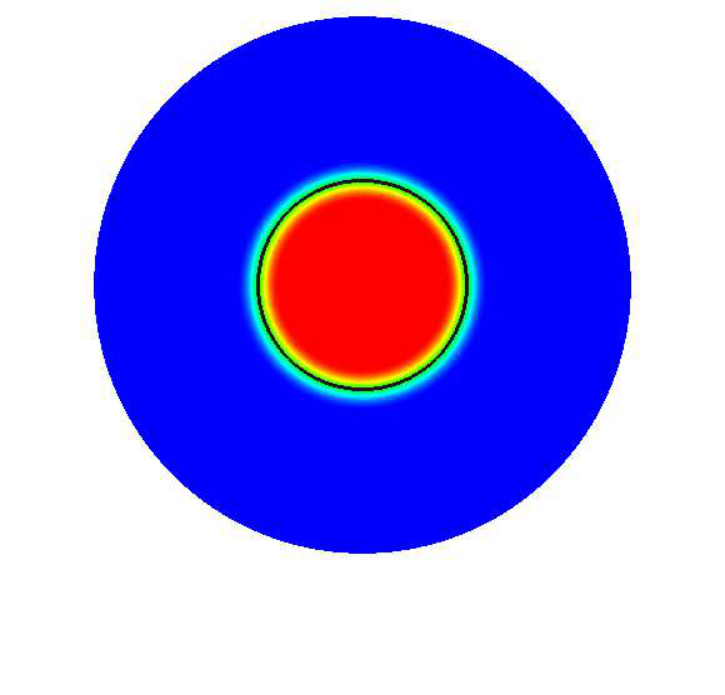}}  
\subfigure{\includegraphics[width = 0.28\textwidth]{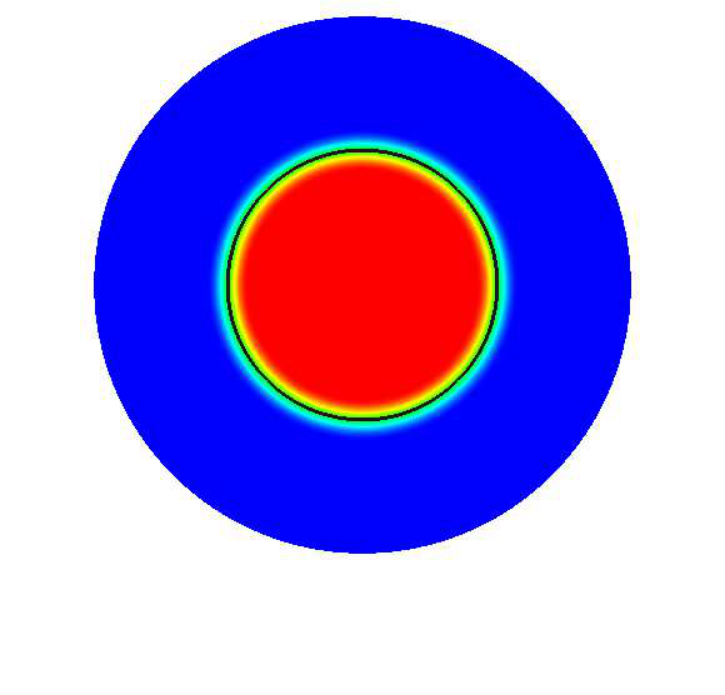}} \\[-6mm]
\subfigure{\includegraphics[width = 0.28\textwidth]{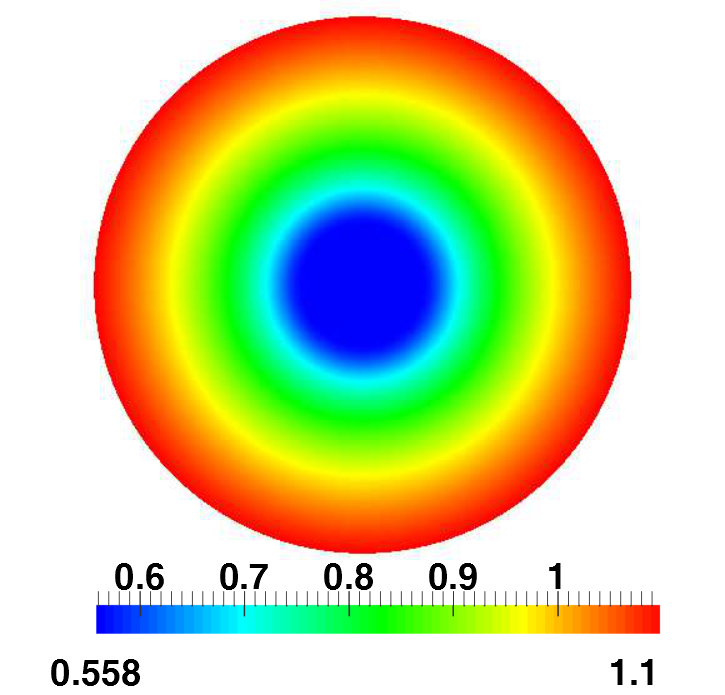}} 
\subfigure{\includegraphics[width = 0.28\textwidth]{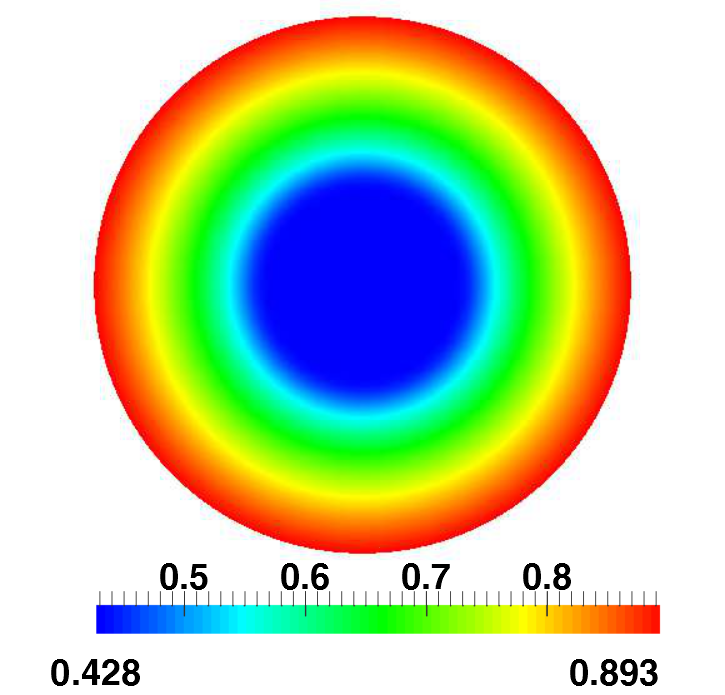}}  
\subfigure{\includegraphics[width = 0.28\textwidth]{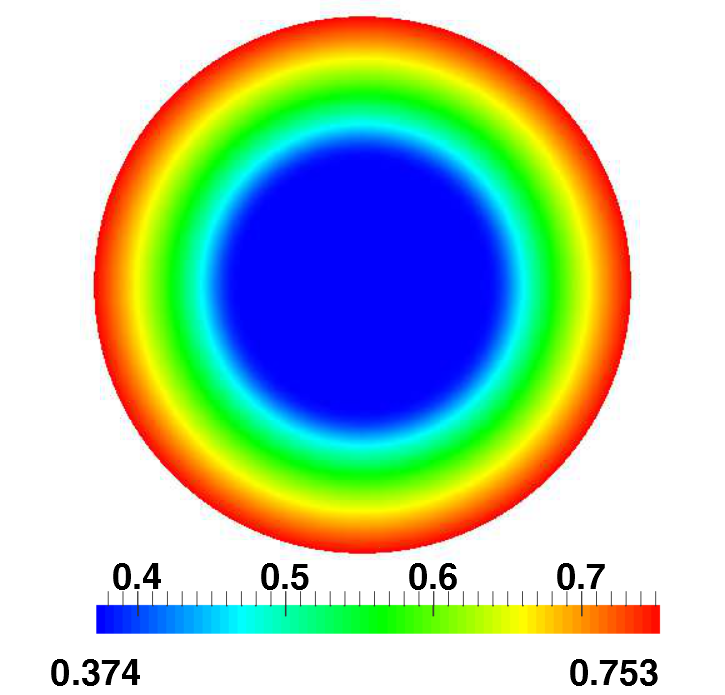}} 
\caption{Approximate solutions of (\ref{Numerical:Hilhorst}) at $t=0$ (left), $t=0.2$ (centre) and $t=0.4$, $\varphi$ (top row), $\sigma$ bottom row. The black 
line in the $\varphi$ solutions denotes the corresponding sharp interface solution.}
\label{f:ct_pf_sol}
\end{figure}

\begin{figure}[!h]
\centering
\subfigure[radius versus time]{\includegraphics[width = 0.45\textwidth]{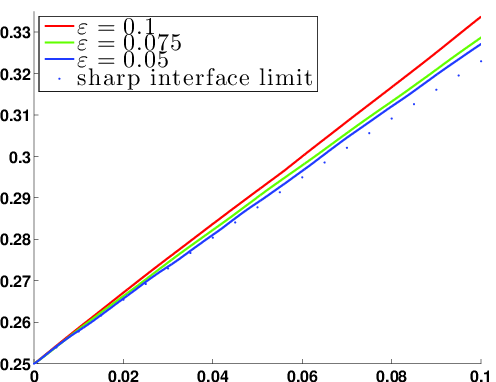}} \hspace{2mm}
\subfigure[$\sigma$ at $t=0.1$]{\includegraphics[width = 0.45\textwidth]{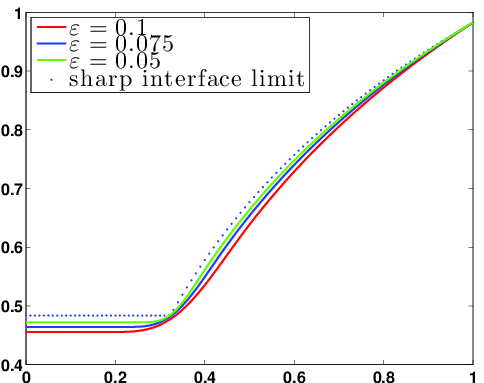}}  
\caption{Comparison of diffuse interface model (\ref{Numerical:Hilhorst}) with the sharp interface solution.}
\label{f:ct_eps_comp}
\end{figure}

In Figures \ref{f:ct_pf_sol} and \ref{f:ct_eps_comp} we display results obtained from 
the growing circle tumour test case introduced in Section 4.2 of \cite{article:HilhorstKaampmannNguyenZee15}.  To this end we consider the simplified model on a circular domain $\Omega$ with radius $R$:
\begin{subequations}\label{Numerical:Hilhorst}
\begin{align}
\pd_{t}\varphi  & = \Laplace \mu +  \frac{1}{\eps} \frac{4\sqrt{2}}{\pi}(1-\varphi^{2}) \sigma, \label{gcta} \\
\mu & = \frac{1}{\eps} \Psi'(\varphi) - \eps \Laplace \varphi, \label{gctb}  \\
0&= \Laplace \sigma - \frac{1}{\eps} \frac{4\sqrt{2}}{\pi}(1-\varphi^{2})\sigma. \label{gctc}
\end{align}
\end{subequations}
Here $\varphi$ and $\mu$ satisfy homogeneous Neumann boundary conditions, and $\sigma$ satisfies the Dirichlet boundary condition $\sigma = \sigma_{R} \in \R$ on $\pd \Omega$. 
We take the radially symmetric case of an initial circular tumour with initial radius $0.25$.  From \cite{article:HilhorstKaampmannNguyenZee15} we have that the solution to the sharp interface limit of (\ref{Numerical:Hilhorst}) is given by
\begin{equation}\label{eq:aa}
\begin{aligned}
\sigma(r,t)= \begin{cases} \sigma_{\rho(t)} & r \leq \rho(t),\\
\sigma_{R} - \frac{\log(r/R)}{\log(\rho(t)/R)}(\sigma_{R} - \sigma_{\rho(t)}) & r>\rho(t),\end{cases}
\end{aligned}
\end{equation}
where $\sigma_{\rho(t)} = \frac{\sigma_{R}}{1-2\sqrt{2}\rho(t)\log(\rho(t)/R)}$, with $\mu$ being constant and $\rho(t)$, which is the radius of the tumour, is determined by numerically solving the ODE $\rho'(t) = \sqrt{2} \sigma(\rho(t),t)$ with initial condition $\rho(0) = 0.25$. 

We set $R = 10$ and $\sigma_{R} = 2$, however for the diffuse interface computations we did not solve the problem in the whole of $\Omega$ instead we solved it on a circular domain with radius $2$ with the time dependent Dirichlet boundary condition $\sigma(x,t) = \sigma_{D}(\abs{x},t)$ computed from (\ref{eq:aa}) with $r = 1$.  We set $\tau = 1.0e^{-4}$, the minimal diameter of an element $h_{min} = 7.8125 \cdot 10^{-3}$ and the maximal diameter $h_{max} = 3.125 \cdot 10^{-2}$. 

In Figure \ref{f:ct_pf_sol} we display the diffuse interface solutions $\varphi$ and $\sigma$ at $t = 0,~0.2,~0.4$ obtained with $\eps = 0.05$. In the plots of $\varphi$ we include the sharp interface limit solution of the tumour position.  In Figure \ref{f:ct_eps_comp} we examine the convergence of the diffuse interface solution to the sharp interface limit solution as $\eps$ tends to zero. In Figure \ref{f:ct_eps_comp}(a) we plot the radius of the growing tumour for the diffuse interface model with $\eps = 0.1,~0.075,~0.05$ together with the sharp interface limit solution $\rho(t)$. In Figure \ref{f:ct_eps_comp}(b) we plot the solution $\sigma$ of the diffuse interface model with $\eps = 0.1,~0.075,~0.05$ together with the sharp interface limit solution $\sigma$ at $t=0.1$.  From this figure we see that as $\eps$ decreases the diffuse interface 
solution converges to the sharp interface limit solution. 

\subsection{Solutions of \eqref{FE:Cristini}}\label{sec:FE:Cristini}
We now investigate the influence of the parameters $\mathcal{P}$, $\chi_{\varphi}$ and $\lambda$ in Model (\ref{Model:Cristini}).  In all computations we set $\Omega = (-12.5,12.5)^{2}$, $\mathcal{A} = 0$, $D = 1$, $\beta = 0.1$, $\mathcal{C} = 2$, $\sigma_{B} = 1$, $\tau = 1.0e^{-3}$, the minimal diameter of an element $h_{min} = 4.888 \cdot 10^{-4}$ and the maximal diameter $h_{max} = 5 \cdot 10^{-1}$. Unless otherwise specified we take $\eps = 0.01$. 

\subsubsection*{Influence of the proliferation rate $\mathcal{P}$}
\begin{figure}[h]
\begin{center}
\subfigure{\includegraphics[width=.24\textwidth,angle=0]{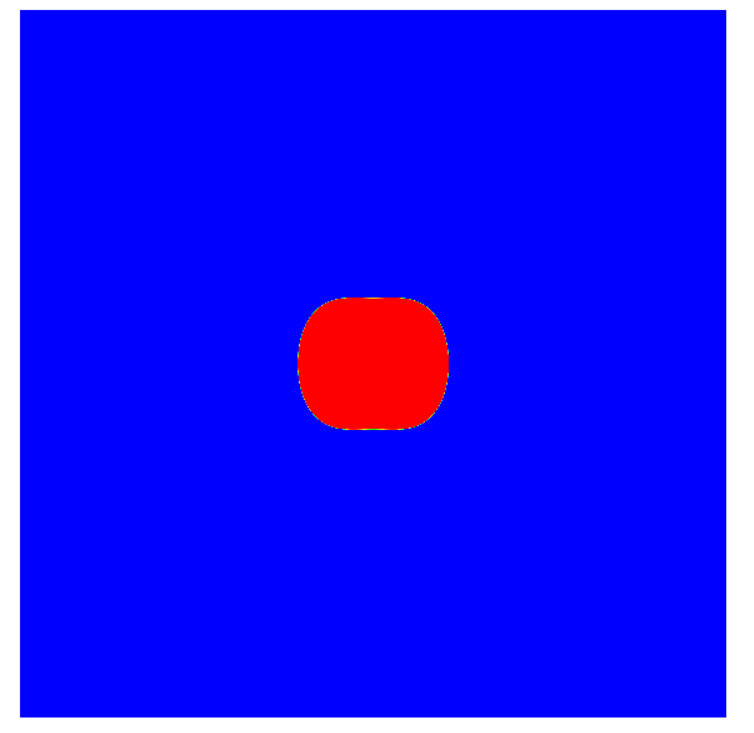}}\hspace{4mm}
\subfigure{\includegraphics[width=.24\textwidth,angle=0]{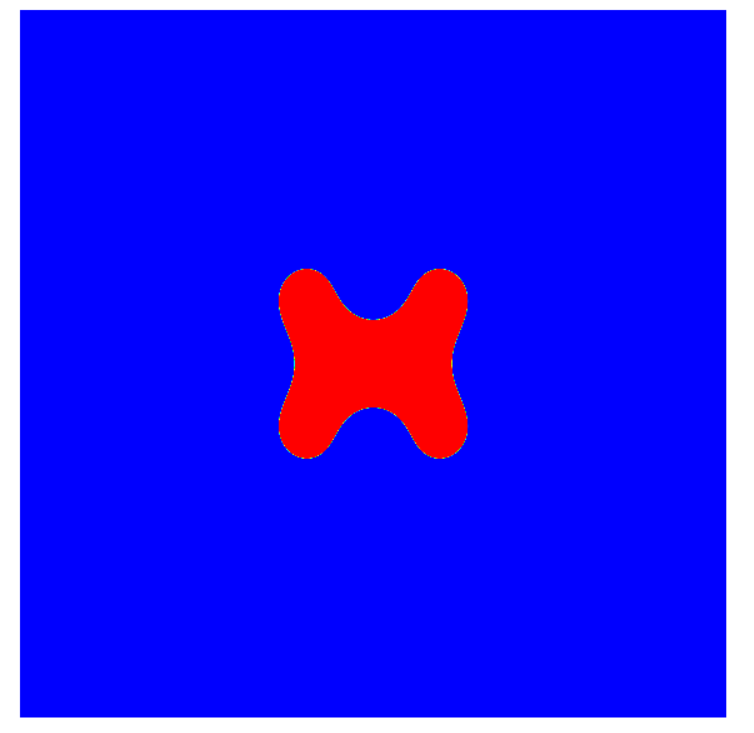}}\hspace{4mm}
\subfigure{\includegraphics[width=.24\textwidth,angle=0]{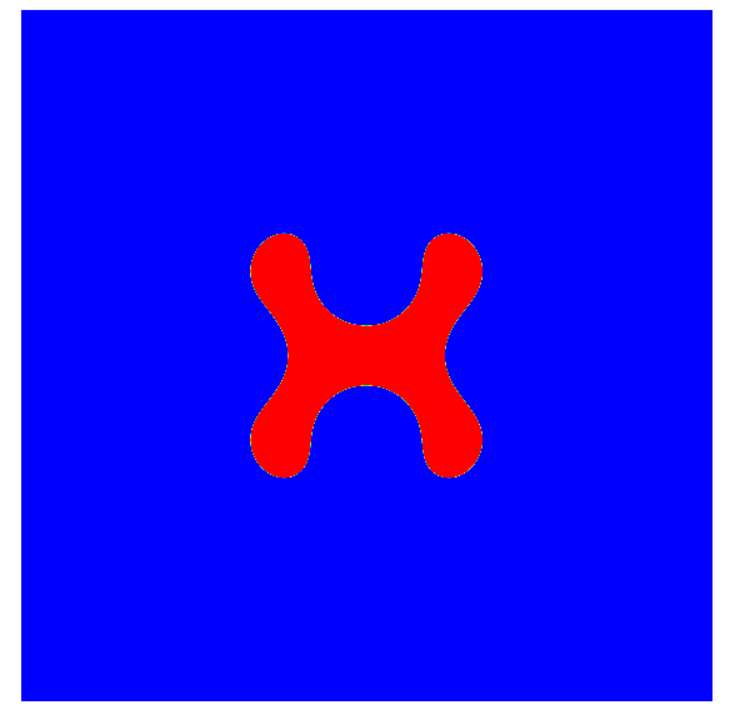}}\\[0mm]
\subfigure{\includegraphics[width=.27\textwidth,angle=0]{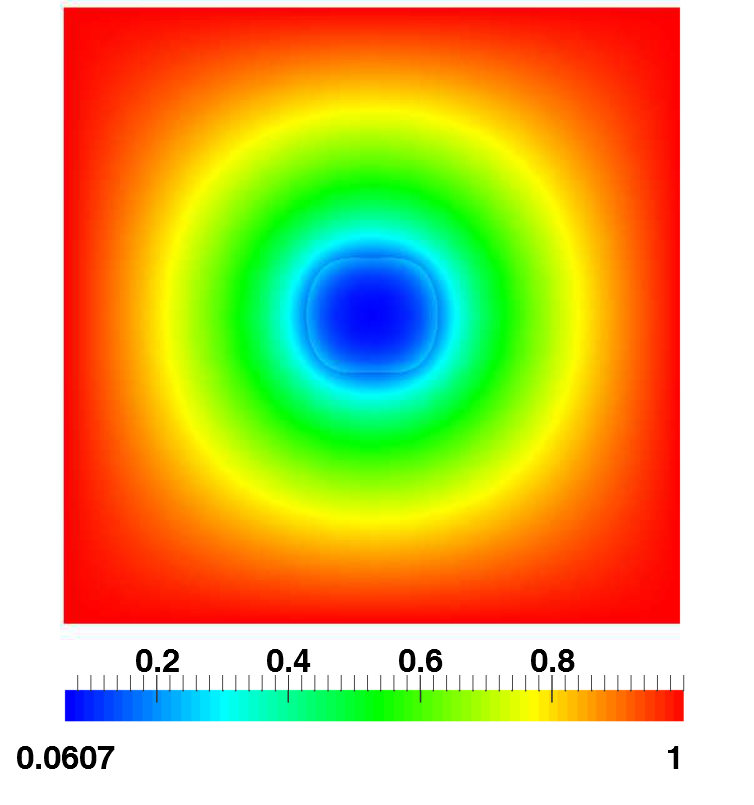}}\hspace{0mm}
\subfigure{\includegraphics[width=.27\textwidth,angle=0]{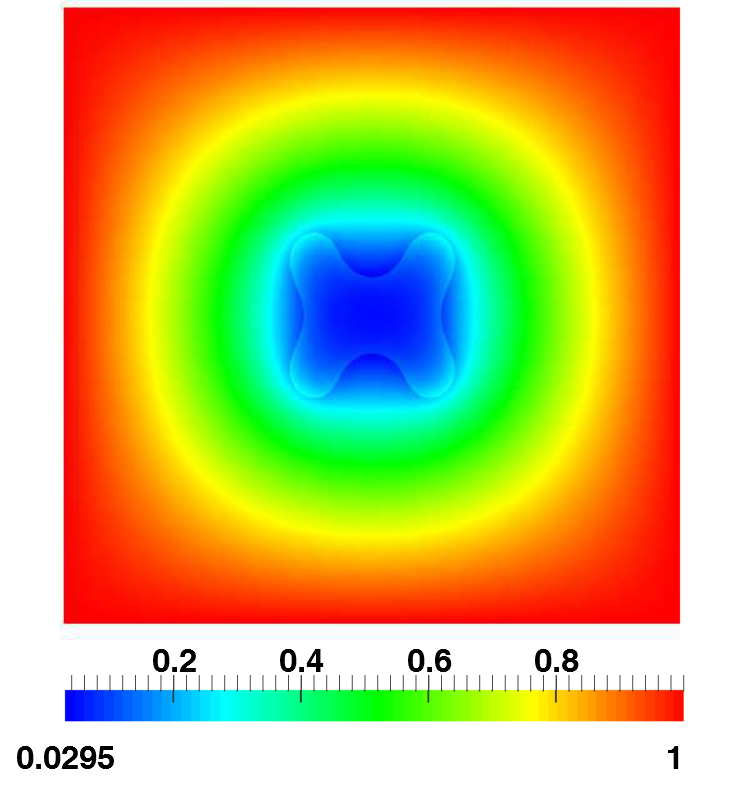}}
\subfigure{\includegraphics[width=.27\textwidth,angle=0]{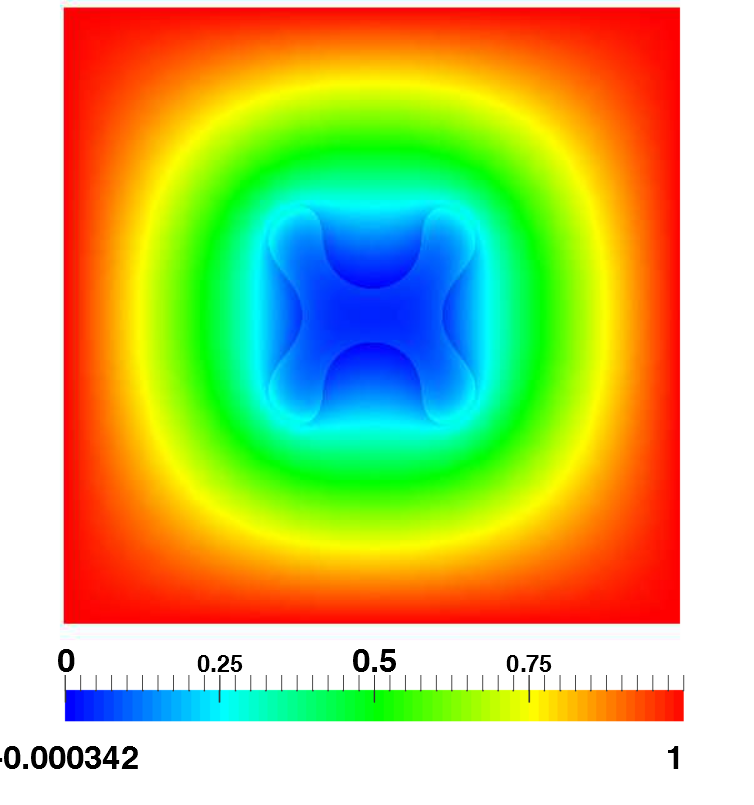}}
\caption{Solutions of (\ref{FE:Cristini}) with $\lambda = 0.03$, $\chi_{\varphi} = 5$, $\mathcal{P} = 0.5$, at $t=5,10,13$.}
\label{f:s2p12}
\end{center}
\end{figure}

\begin{figure}[h]
\begin{center}
\subfigure{\includegraphics[width=.24\textwidth,angle=0]{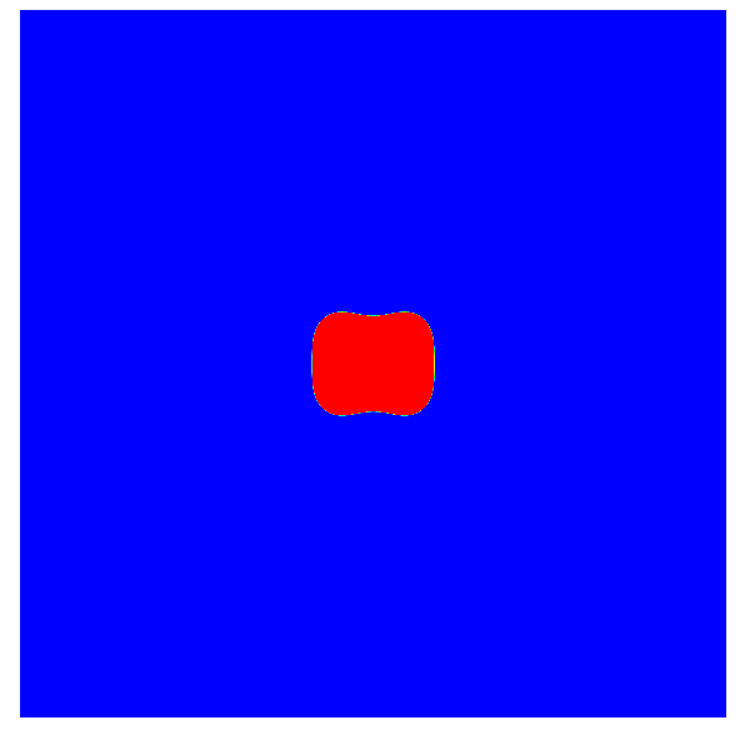}}\hspace{4mm}
\subfigure{\includegraphics[width=.24\textwidth,angle=0]{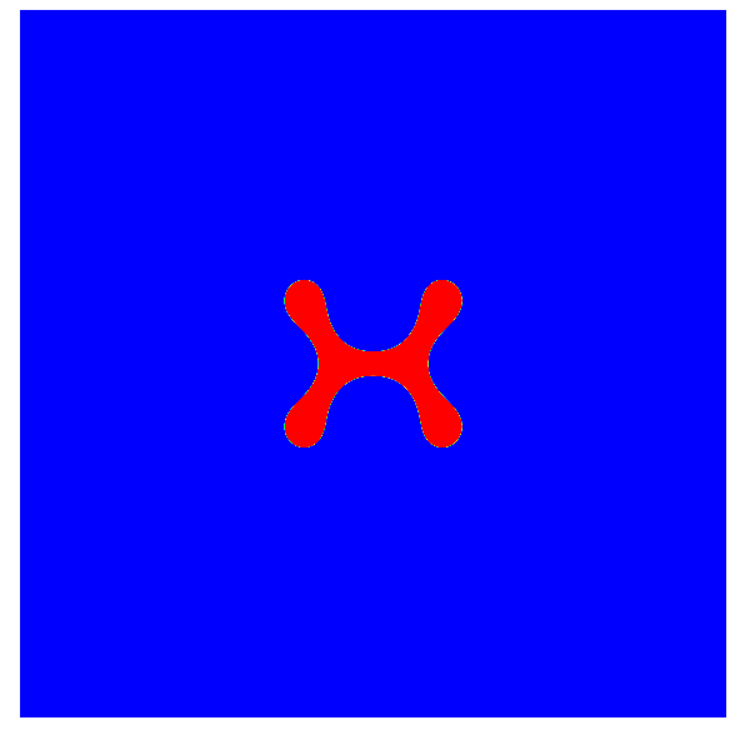}}\hspace{4mm}
\subfigure{\includegraphics[width=.24\textwidth,angle=0]{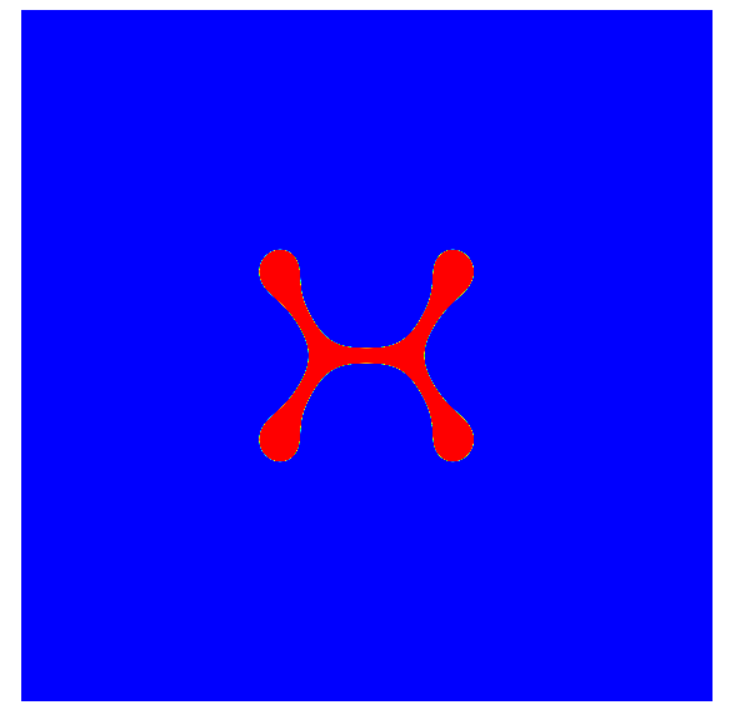}}\\[0mm]
\subfigure{\includegraphics[width=.27\textwidth,angle=0]{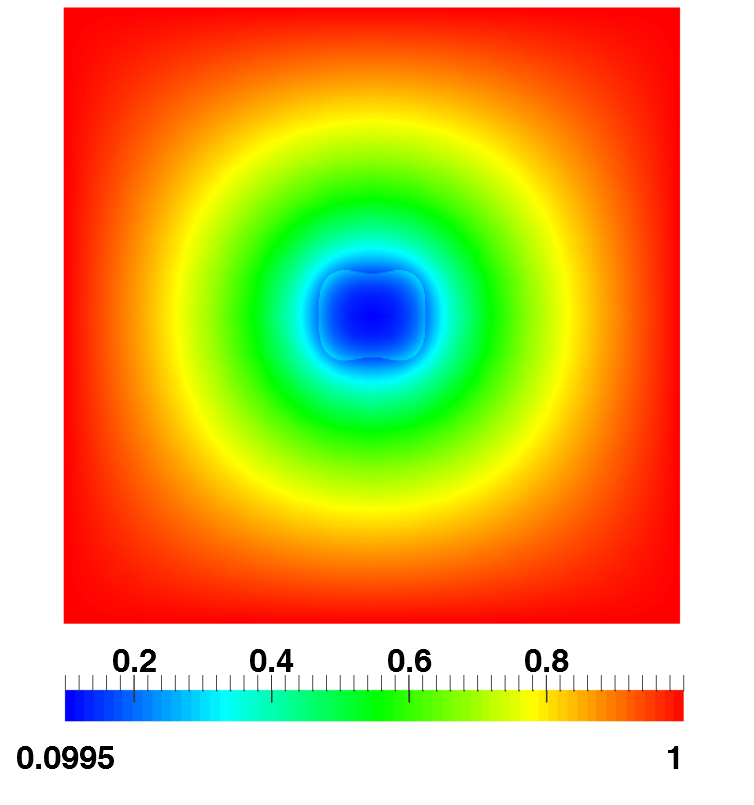}}\hspace{0mm}
\subfigure{\includegraphics[width=.27\textwidth,angle=0]{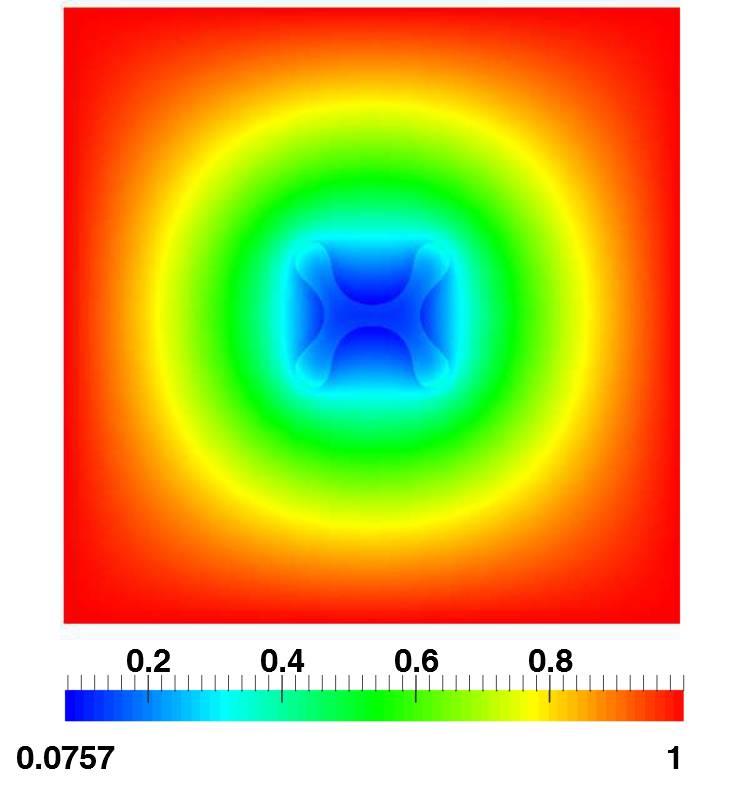}}
\subfigure{\includegraphics[width=.27\textwidth,angle=0]{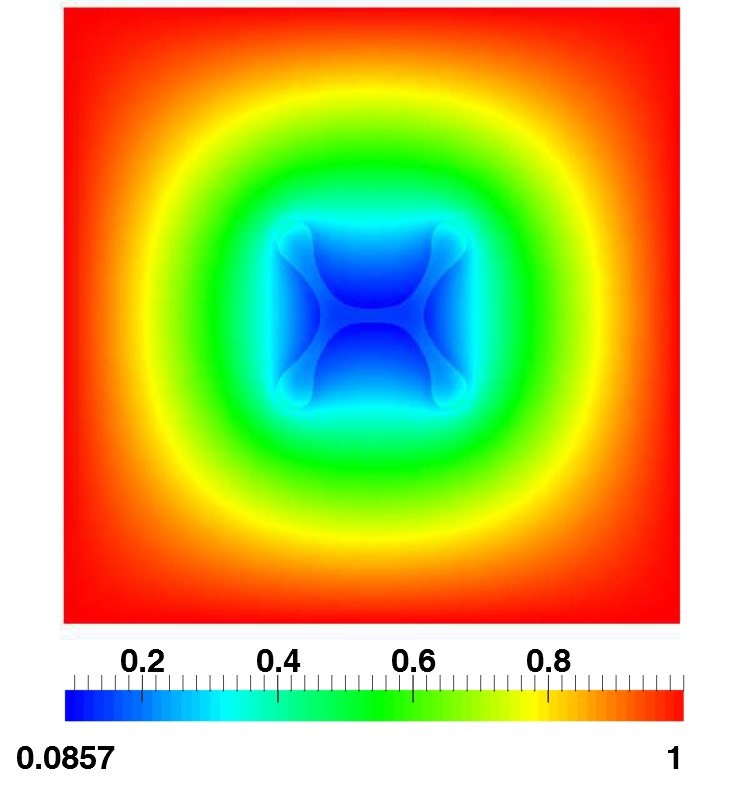}}
\caption{Solutions of (\ref{FE:Cristini}) with $\lambda=0.03$, $\chi_\varphi=5$, $\mathcal{P}=0.1$ at $t=5,10,13$.}
\label{f:s2p2}
\end{center}
\end{figure}

In Figures \ref{f:s2p12} and \ref{f:s2p2} we investigate the influence of $\mathcal{P}$. We set $\chi_{\varphi} = 5$ and $\lambda = 0.03$.  In Figure \ref{f:s2p12} we set $\mathcal{P} = 0.5$ while in Figure \ref{f:s2p2} we set $\mathcal{P} = 0.1$, and in both sets of figures we display $\varphi$ (top row) and $\sigma$ (bottom row) at times $t = 5,10,13$.  From this figure we see taking $\mathcal{P} = 0.5$ gives rise to fingers that are thicker than the ones resulting from $\mathcal{P} = 0.1$.

\subsubsection*{Influence of the chemotaxis parameter $\chi_{\varphi}$}
\begin{figure}[h]
\begin{center}
\subfigure{\includegraphics[width=.24\textwidth,angle=0]{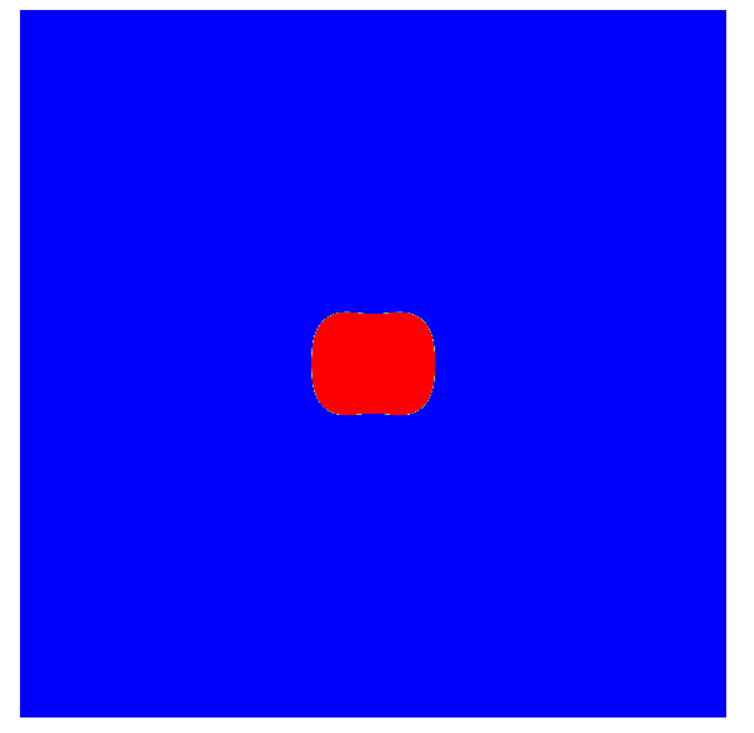}}\hspace{4mm}
\subfigure{\includegraphics[width=.24\textwidth,angle=0]{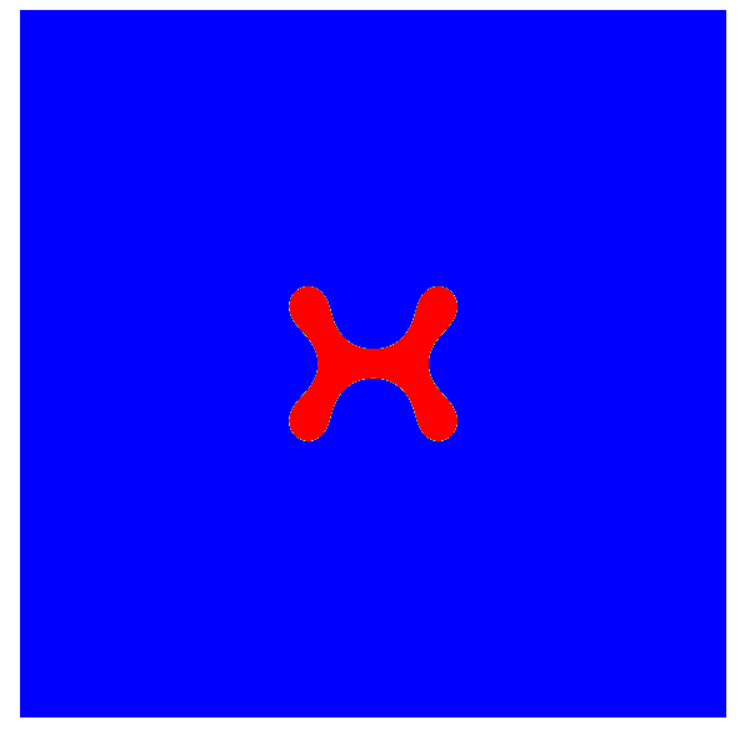}}\hspace{4mm}
\subfigure{\includegraphics[width=.24\textwidth,angle=0]{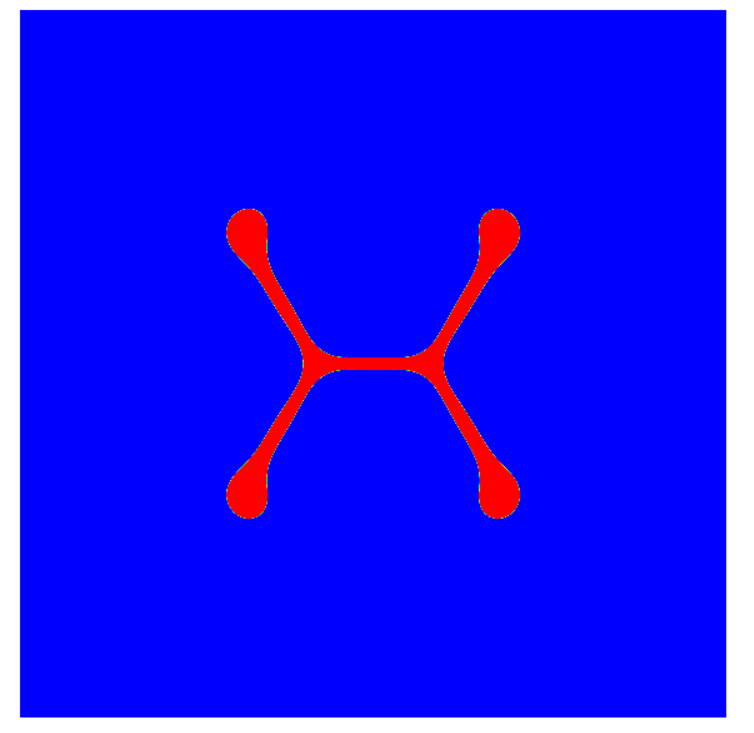}}\\[0mm]
\subfigure{\includegraphics[width=.27\textwidth,angle=0]{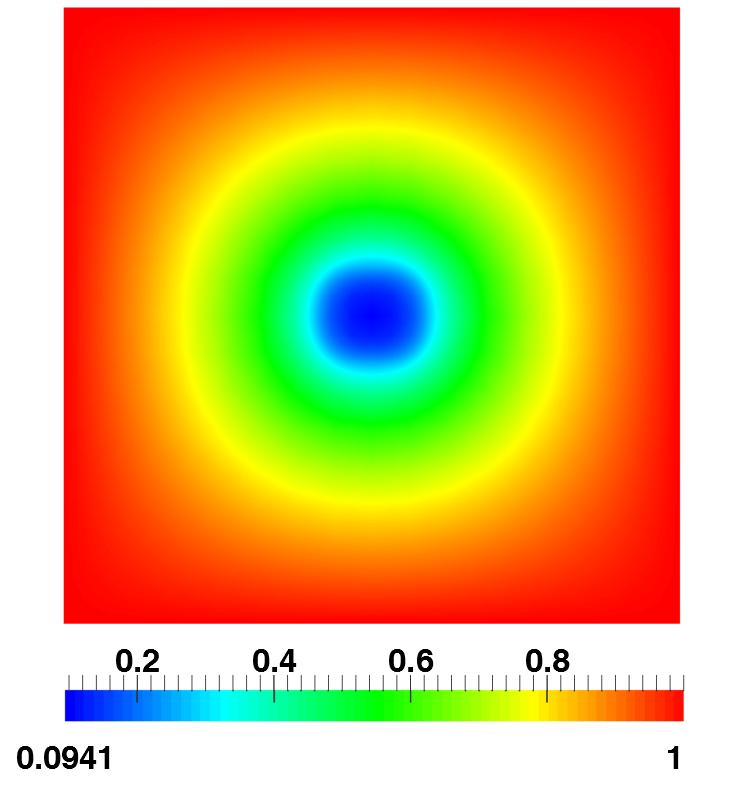}}\hspace{0mm}
\subfigure{\includegraphics[width=.27\textwidth,angle=0]{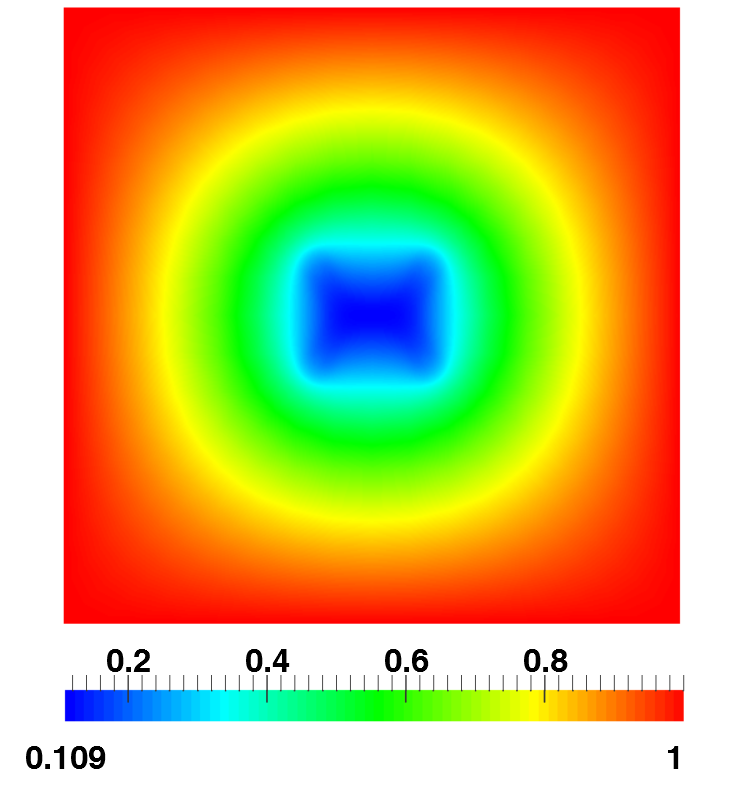}}
\subfigure{\includegraphics[width=.27\textwidth,angle=0]{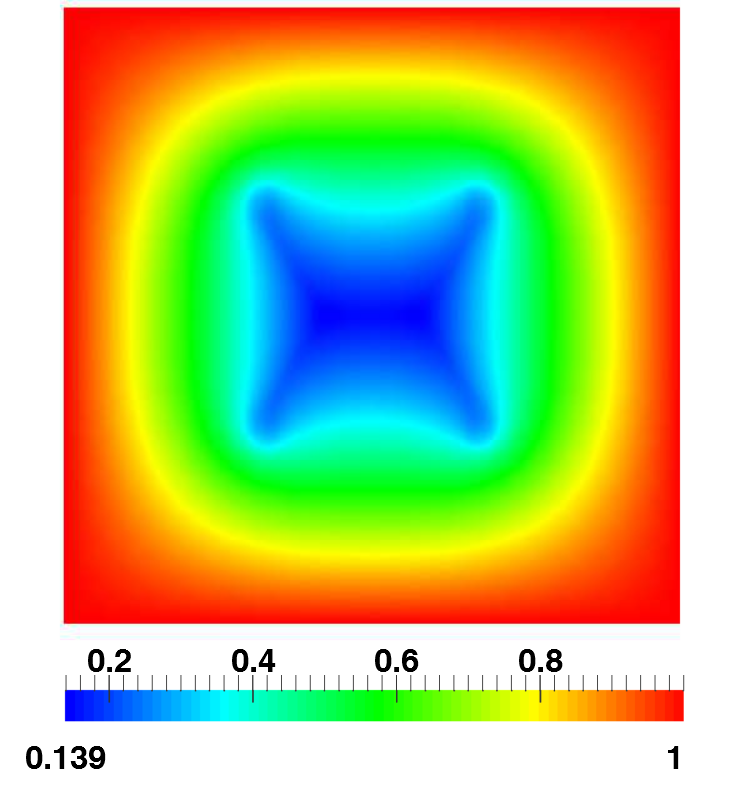}}
\caption{Solutions of (\ref{FE:Cristini}) with $\mathcal{P} = 0.1$, $\lambda = 0$, $\chi_{\varphi} = 5$ at $t=5,10,20$.}
\label{f:s3p14}
\end{center}
\end{figure}

\begin{figure}[h]
\begin{center}
\subfigure{\includegraphics[width=.24\textwidth,angle=0]{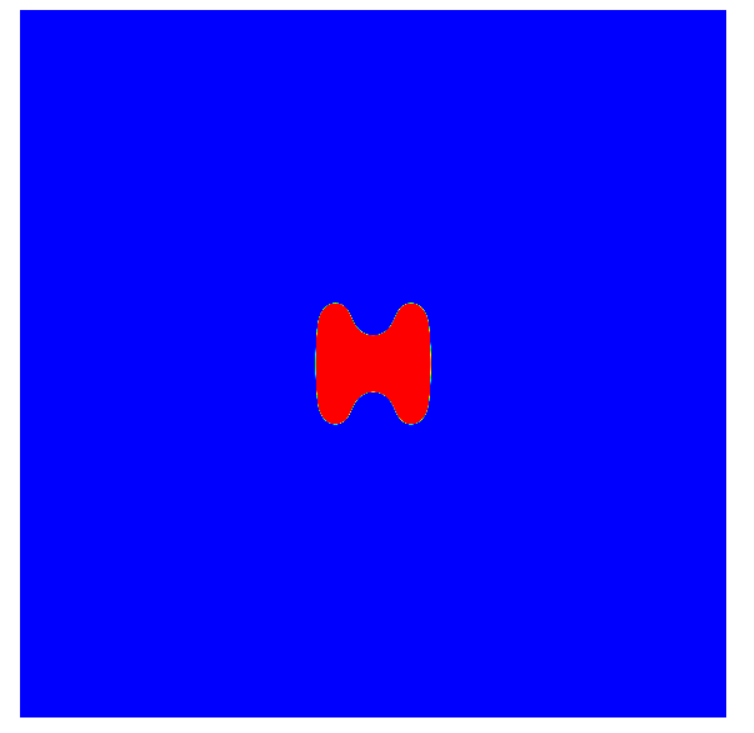}}\hspace{4mm}
\subfigure{\includegraphics[width=.24\textwidth,angle=0]{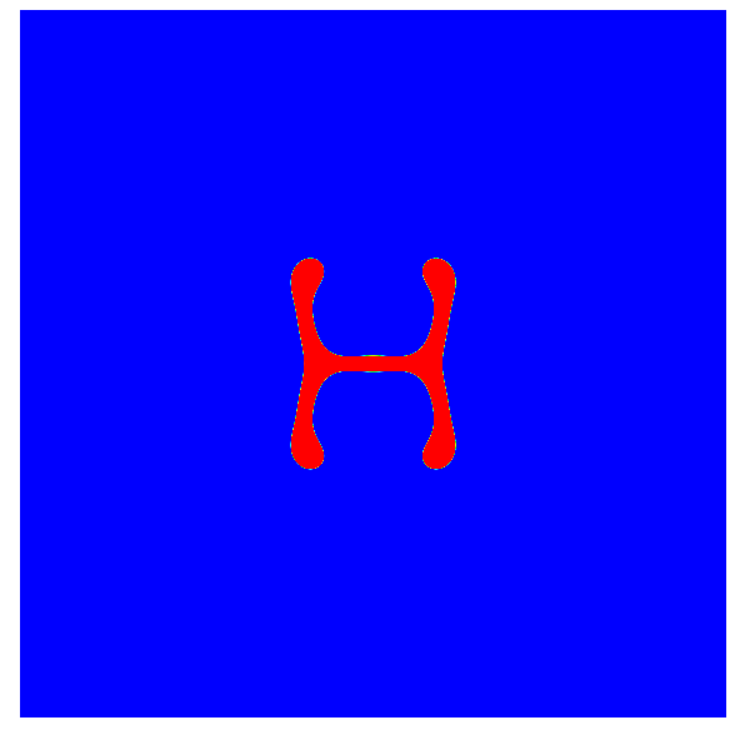}}\hspace{4mm}
\subfigure{\includegraphics[width=.24\textwidth,angle=0]{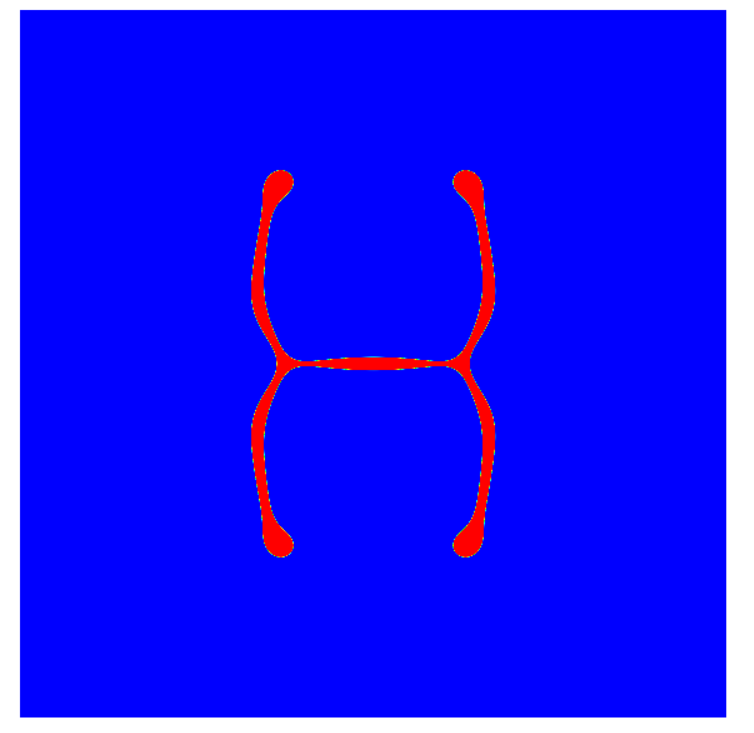}}\\[0mm]
\subfigure{\includegraphics[width=.27\textwidth,angle=0]{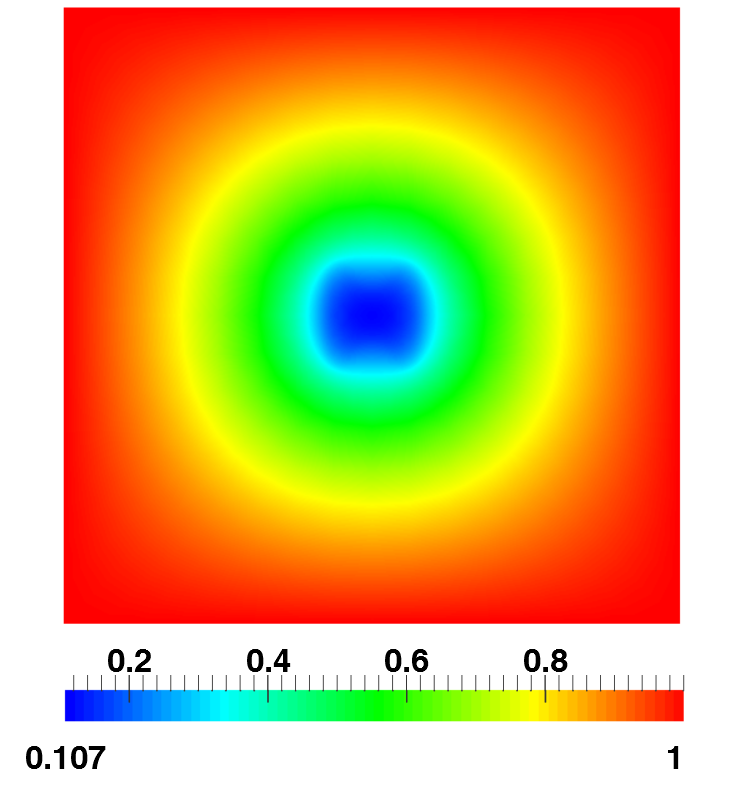}}\hspace{0mm}
\subfigure{\includegraphics[width=.27\textwidth,angle=0]{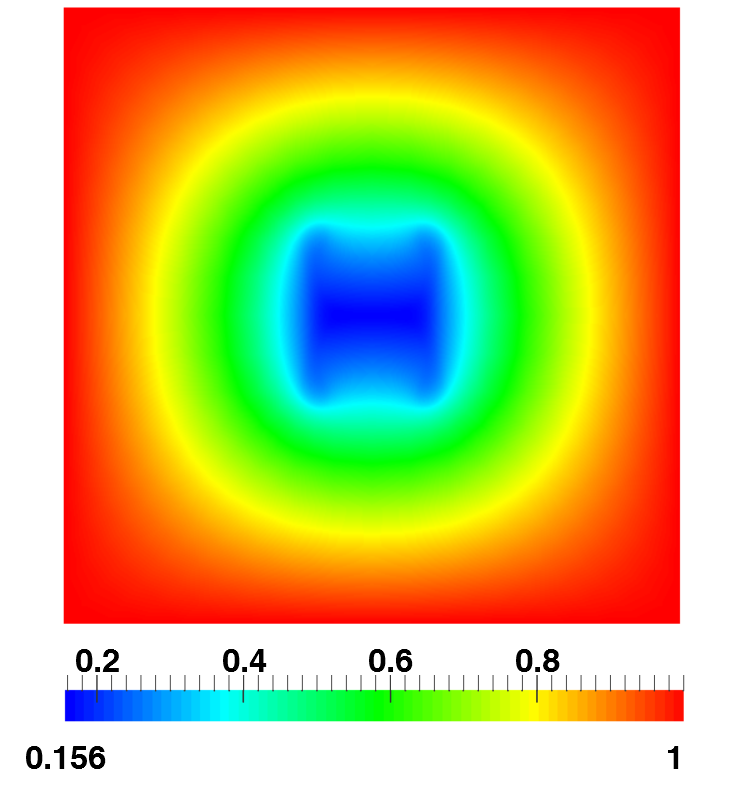}}
\subfigure{\includegraphics[width=.27\textwidth,angle=0]{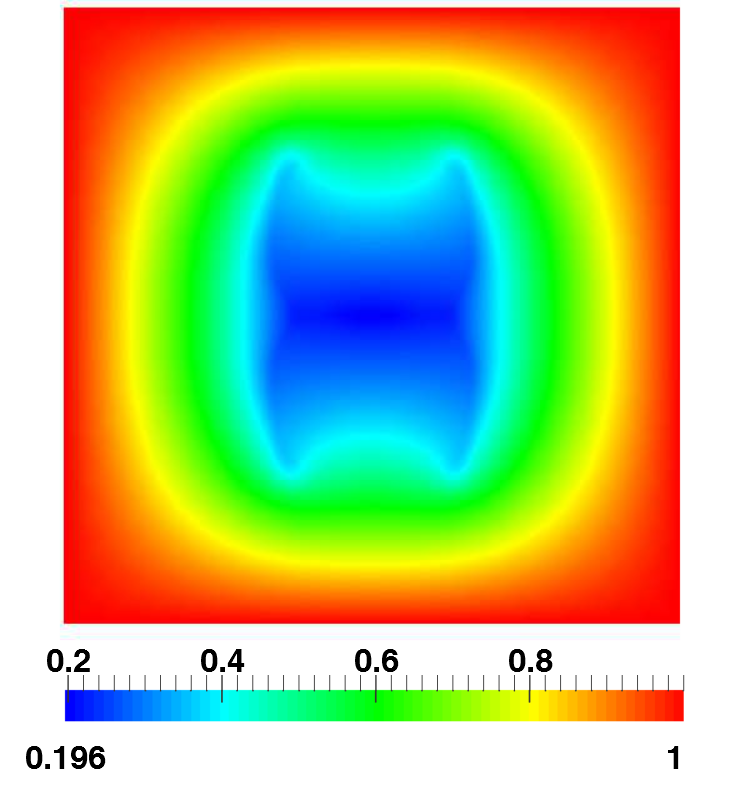}}
\caption{Solutions of (\ref{FE:Cristini}) with $\mathcal{P} = 0.1$, $\lambda = 0$, $\chi_{\varphi} = 10$ at $t=2.5,5,10$.}
\label{f:s3p24}
\end{center}
\end{figure}

In Figures \ref{f:s3p14} and \ref{f:s3p24} we investigate the influence of $\chi_{\varphi}$.  We set $\mathcal{P} = 0.1$ and $\lambda = 0$. 
In Figure \ref{f:s3p14} we set $\chi_{\varphi} = 5$ while in Figure \ref{f:s3p24} we set $\chi_{\varphi} = 10$, and in both sets of figures we display $\varphi$ (top row) and $\sigma$ (bottom row). The results for $\chi_{\varphi} = 5$ are displayed at times $t = 5,10,20$, while the results for $\chi_{\varphi} = 10$ are displayed at times $t = 2.5,5,10$.  From these figures we see that, akin to the results in \cite{article:CristiniLiLowengrubWise09}, 
for both values of $\chi_{\varphi}$ after some time fingers develop, and thereby increasing the surface area of the tumour to allow for better access to the nutrient.  For the larger value of $\chi_{\varphi}$ the formation and evolution of the fingers is quicker and the fingers are slimmer.  

\begin{figure}[!h]
\begin{center}
\subfigure{\includegraphics[width = 0.232\textwidth]{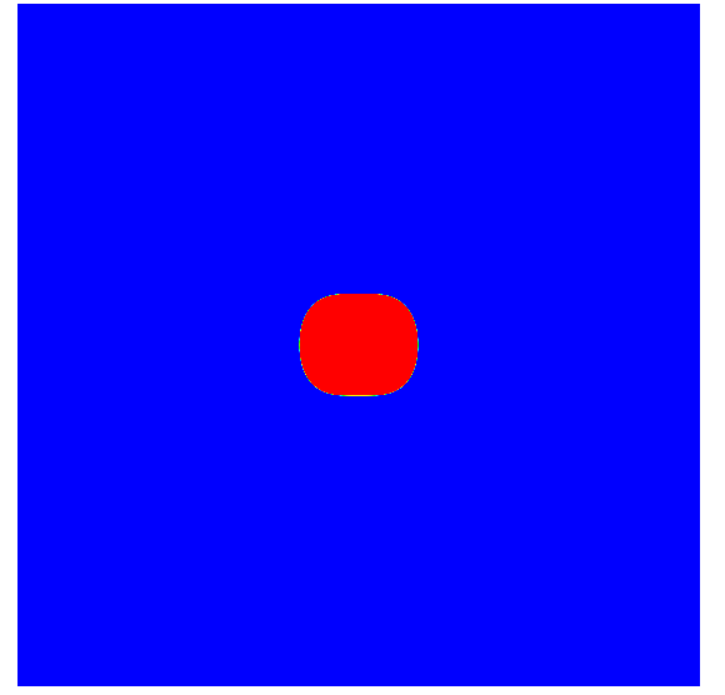}}\hspace{4mm}
\subfigure{\includegraphics[width = 0.232\textwidth]{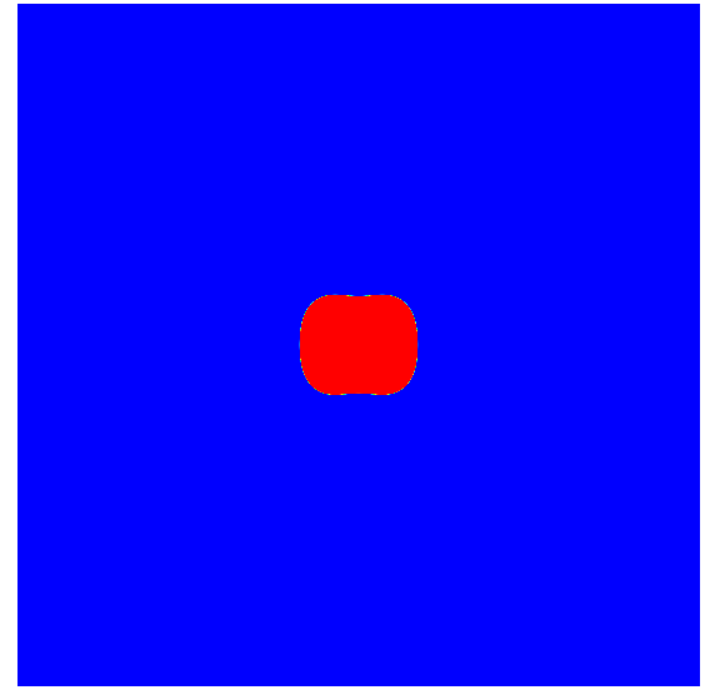}}  \hspace{4mm}
\subfigure{\includegraphics[width = 0.232\textwidth]{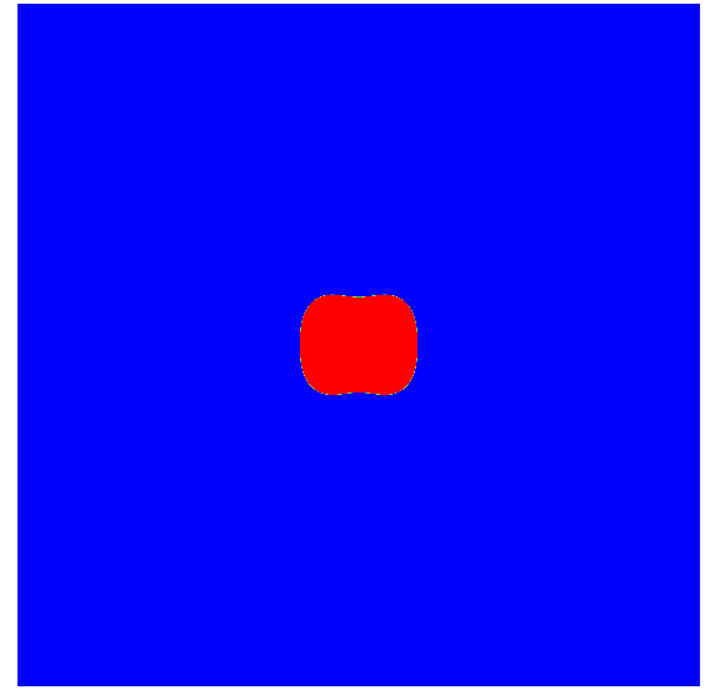}}  \\[0mm]
\subfigure{\includegraphics[width=.27\textwidth,angle=0]{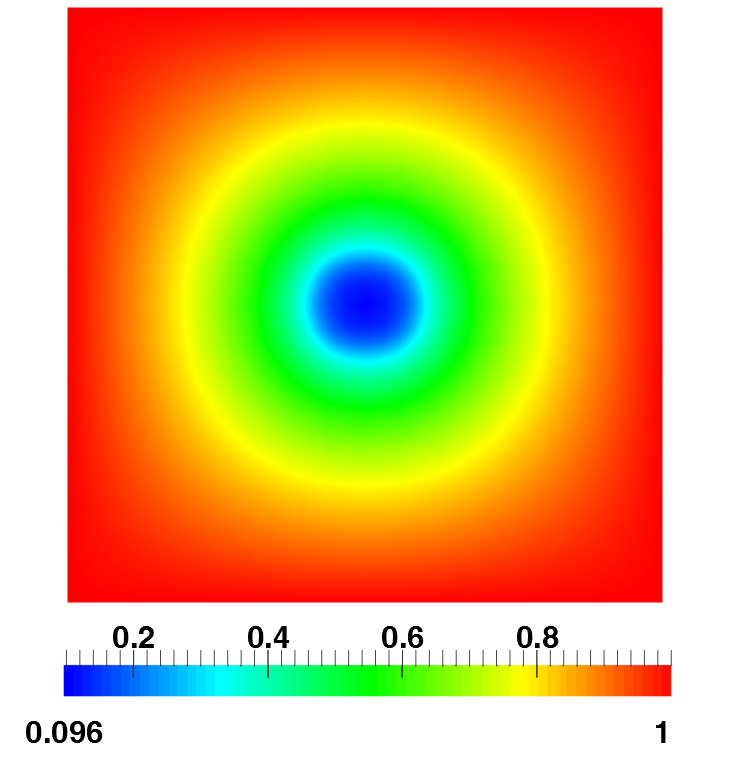}}\hspace{0mm}
\subfigure{\includegraphics[width=.27\textwidth,angle=0]{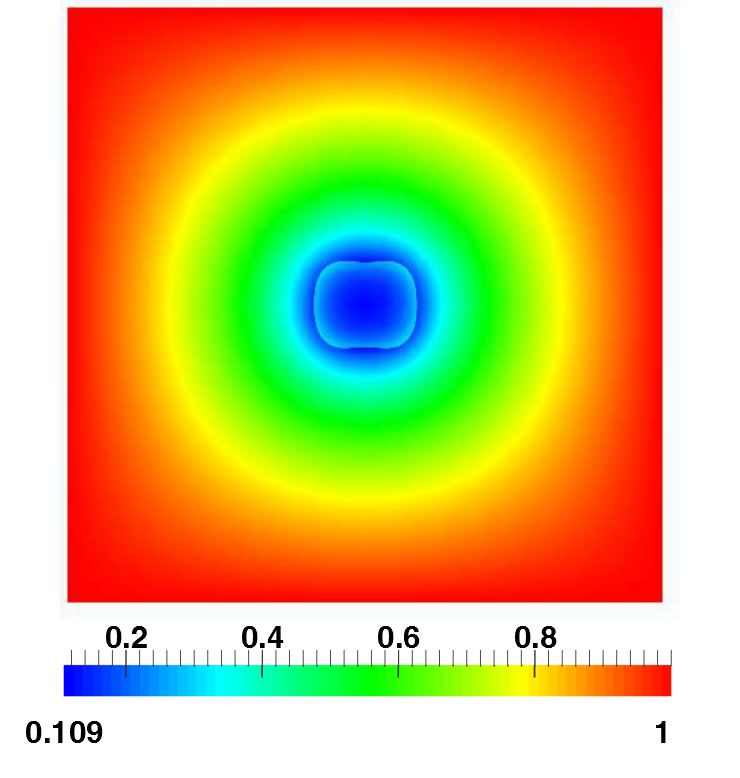}}
\subfigure{\includegraphics[width=.27\textwidth,angle=0]{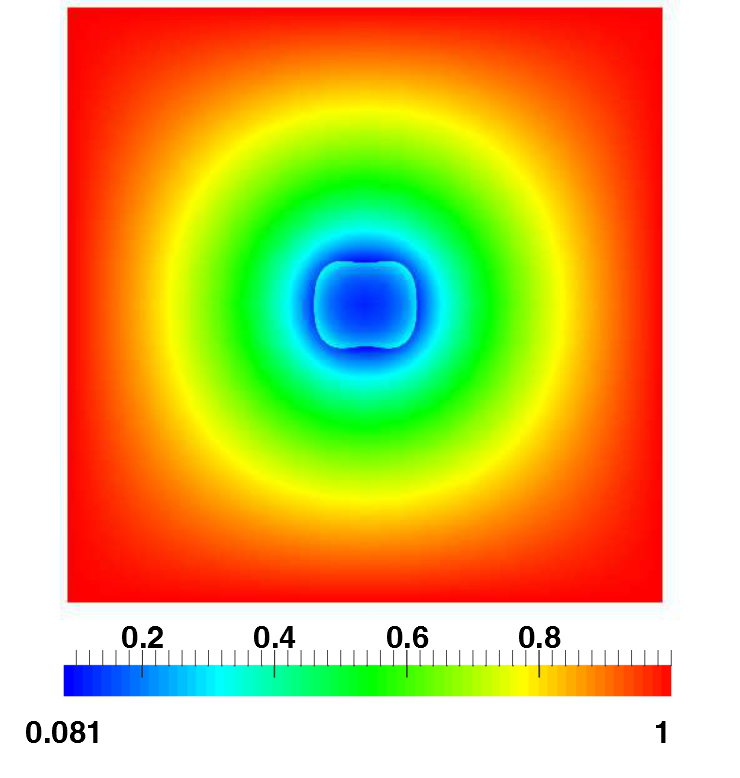}}
\caption{Solutions of (\ref{FE:Cristini}) with $\chi_{\varphi} = 5$, $\mathcal{P} = 0.1$, at $t = 4$, $\lambda = 0$ (left), $\lambda = 0.07$ (centre) 
and $\lambda = 0.09$ (right).}
\label{f:jump1}
\end{center}
\end{figure} 

\subsubsection*{Influence of the active transport parameter $\lambda$}
\begin{figure}[!h]
\centering %original width = 0.28\textwidth
\subfigure{\includegraphics[width = 0.32\textwidth]{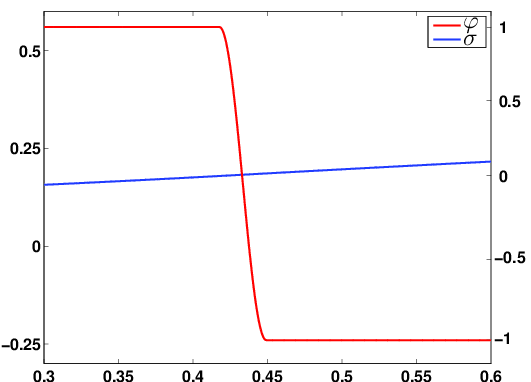}} \hspace{-1mm}
\subfigure{\includegraphics[width = 0.32\textwidth]{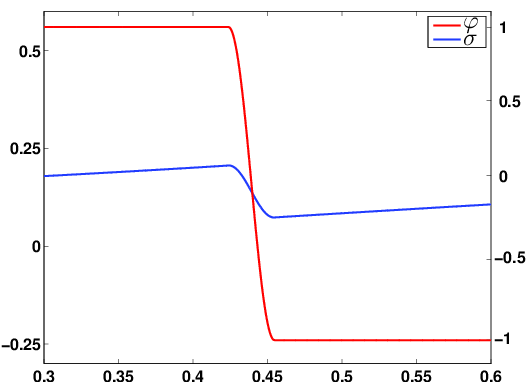}}  \hspace{-1mm}
\subfigure{\includegraphics[width = 0.32\textwidth]{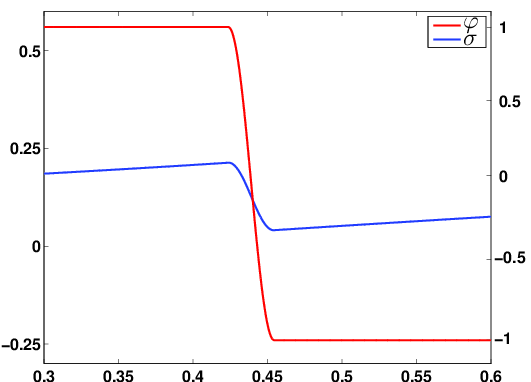}} 
\caption{Comparison of $\sigma$ across the interfacial region for (\ref{FE:Cristini}) at $t = 4$ with $\lambda = 0$ (left), $\lambda = 0.07$ (centre) and $\lambda = 0.09$ (right).}
\label{f:jump2}
\end{figure}

\begin{figure}[!h]
\centering
\subfigure{\includegraphics[width = 0.5\textwidth]{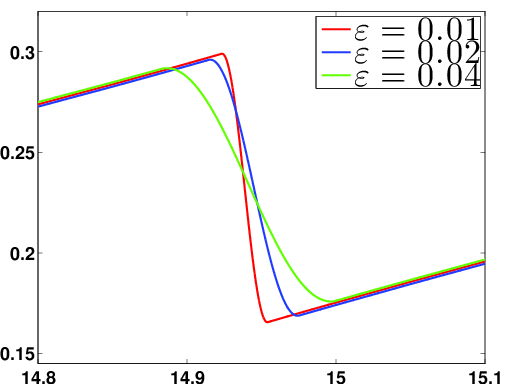}} 
\caption{Convergence of $\sigma$ as $\eps$ decreases for (\ref{FE:Cristini}) with $\lambda = 0.07$ at $t = 4$.}
\label{f:jump3}
\end{figure}

In Figures \ref{f:jump1} - \ref{f:jump3} we investigate the influence of $\lambda$. We set $\mathcal{P} = 0.1$ and $\chi_{\varphi} = 5$.  In Figure \ref{f:jump1} we show $\varphi$ (top row) and $\sigma$ (bottom row) at  $t=4$, with $\lambda = 0$ (left), $\lambda = 0.07$ (centre) and $\lambda = 0.09$ (right).  From this figure we see that when $\lambda = 0$ the variation of $\sigma$ across the interfacial region is smooth while taking $\lambda > 0$ leads to a drastic change in $\sigma$.  This change in $\sigma$ can be seen better in Figure \ref{f:jump2} where we show plots of $\varphi$ and $\sigma$ along a line that spans the interfacial region. The scales for $\sigma$ and $\varphi$ are shown on the left and right axes respectively. 

Here we see that the change in $\sigma$ across the interfacial region is more pronounced for larger values of $\lambda$.  In Figure \ref{f:jump3} we display the influence of $\eps$ on the change in $\sigma$ across the interfacial region, we 
set $\lambda = 0.07$ and plot $\sigma$ along a line that spans the interfacial region for $\eps = 0.04,~0.02,~0.01$.  From this figure we see the convergence of $\sigma$ as $\eps$ decreases.  In Figure \ref{f:jump3} the jump in $\sigma$ across the interfacial region for $\eps = 0.01$ is $0.1327 \approx 2\lambda$ which is consistent with the formal asymptotic analysis, recall (\ref{Cristini:jumpsigmalambda}).  

\subsection{Numerical computations with Darcy flow}
For positive constants $m_{0}$ and $K$, we now consider the model
\begin{subequations}
\begin{align}
\div \bm{v} & = \alpha \Gamma, \\
\bm{v} & = -K (\nabla p - (\mu + \chi_{\varphi} \sigma) \nabla \varphi), \\
\pd_{t} \varphi + \div (\varphi \bm{v}) & = m_{0} \Laplace \mu + \rho_{S} \Gamma, \\
\mu & = \frac{\beta}{\eps} \Psi'(\varphi) - \beta \eps \Laplace \varphi - \chi_{\varphi} \sigma, \\
0 & = \div (\mathcal{D}(\varphi) ( \nabla \sigma - \lambda \nabla \varphi)) - \frac{1}{2} \mathcal{C} \sigma (\varphi + 1),
\end{align}
\end{subequations}
where we recall that $\alpha := \frac{1}{\overline{\rho}_{2}} - \frac{1}{\overline{\rho}_{1}}$, $\rho_{S} := \frac{1}{\overline{\rho}_{2}} + \frac{1}{\overline{\rho}_{1}}$, $\Gamma = \frac{1}{2}(\mathcal{P} \sigma - \mathcal{A})(\varphi + 1)$, and $\mathcal{D}$ is defined in \eqref{Cristini:Specific}.  As additional boundary condition we prescribe 
\begin{align*}
p = 0 \text{ on } \pd \Omega,
\end{align*}
while we take homogeneous Neumann boundary conditions for $\varphi$ and $\mu$, and the Dirichlet boundary condition $\sigma = \sigma_{B} \in \R$ on $\pd \Omega$.  Recalling the finite element spaces $K_{h}$, $S_{h}$, $S_{h}^{B}$ and $S_{h}^{0}$ defined at the start of Section \ref{sec:Numerics}, for the double-obstacle potential \eqref{defn:DoubleObstacle}, we propose the following scheme for the above system: Find
\begin{align*}
\varphi_{h}^{n} \in K_{h}, \quad \mu_{h}^{n} \in S_{h}, \quad \sigma_{h}^{n} \in S_{h}^{B}, \quad p_{h}^{n} \in S_{h}^{0}
\end{align*}
such that for all $\eta_{h} \in S_{h}$, $\zeta_{h} \in K_{h}$ and $\chi_{h} \in S_{h}^{0}$,
\begin{subequations}
\begin{align}
\notag & \; \frac{1}{\tau}(\varphi_{h}^{n} - \varphi_{h}^{n-1}, \eta_{h})_{h} + m_{0}(\nabla \mu_{h}^{n}, \nabla \eta_{h}) \\
\notag & \; = \frac{\rho_{S}}{2}((\mathcal{P} \sigma_{h}^{n-1} - \mathcal{A})(\varphi_{h}^{n} + 1), \eta_{h})_{h} - \frac{\alpha}{2}(\varphi_{h}^{n-1}(\mathcal{P} \sigma_{h}^{n-1} - \mathcal{A})(\varphi_{h}^{n-1} + 1), \eta_{h})_{h} \\
 & \; + K ( \nabla p_{h}^{n-1} \cdot \nabla \varphi_{h}^{n-1} - (\mu_{h}^{n-1} + \chi_{\varphi} \sigma_{h}^{n-1}) \abs{\nabla \varphi_{h}^{n-1}}^{2}, \eta_{h})_{h}, \\
& \; \left ( \mu_{h}^{n} + \frac{\beta}{\eps} \varphi_{h}^{n-1} + \chi_{\varphi} \sigma_{h}^{n-1}, \zeta_{h} - \varphi_{h}^{n} \right )_{h} \leq \beta \eps (\nabla \varphi_{h}^{n}, \nabla (\zeta_{h} - \varphi_{h}^{n})),  \\
& \; (\mathcal{D}(\varphi_{h}^{n}) \nabla \sigma_{h}^{n}, \nabla \chi_{h})_{h} - \lambda (\mathcal{D}(\varphi_{h}^{n}) \nabla \varphi_{h}^{n}, \nabla \chi_{h})_{h} = -\frac{1}{2} \mathcal{C} (\sigma_{h}^{n}(\varphi_{h}^{n} + 1), \chi_{h})_{h},\\
& \;  (\nabla p_{h}^{n}, \nabla \chi_{h}) = ((\mu_{h}^{n} + \chi_{\varphi} \sigma_{h}^{n}) \nabla \varphi_{h}^{n}, \nabla \chi_{h})_{h} + \frac{\alpha}{2K} ((\mathcal{P} \sigma_{h}^{n} - \mathcal{A})(\varphi_{h}^{n} + 1), \chi_{h})_{h}.
\end{align}
\label{eq:hhh}
\end{subequations}
As initial condition for $p$ and $\mu$, we always choose $p_{h}^{0} = 0$ and $\mu_{h}^{0} = 0$.  We perform three different numerical simulations in which we vary the tumour and healthy cell densities.  The three cases are given as follows:
\begin{itemize}
\item (Case (1)) $\alpha = 0$ and $\rho_{S} = 2$ with $\overline{\rho}_{1} = \overline{\rho}_{2} = 1$ so that we solve for
\begin{align*}
\div \bm{v} = 0, \quad \pd_{t}\varphi + \nabla \varphi \cdot \bm{v} = m_{0} \Laplace \mu + (\mathcal{P} \sigma - \mathcal{A})(\varphi + 1);
\end{align*}
\item (Case (2)) $\alpha = \frac{2}{3}$ and $\rho_{S} = 2$ with $\overline{\rho}_{1} = \frac{3}{2}$, $\overline{\rho}_{2} = \frac{3}{4}$ so that we solve for
\begin{align*}
\div \bm{v} = \frac{1}{3}(\mathcal{P} \sigma - \mathcal{A})(\varphi + 1), \quad \pd_{t}\varphi + \div (\varphi \bm{v}) = m_{0} \Laplace \mu + (\mathcal{P} \sigma - \mathcal{A})(\varphi + 1);
\end{align*}
\item (Case (3)) $\alpha = -\frac{2}{3}$ and $\rho_{S} = 2$ with $\overline{\rho}_{1} = \frac{3}{4}$, $\overline{\rho}_{2} = \frac{3}{2}$ so that we solve for
\begin{align*}
\div \bm{v} = -\frac{1}{3}(\mathcal{P} \sigma - \mathcal{A})(\varphi + 1), \quad \pd_{t}\varphi + \div (\varphi \bm{v}) = m_{0} \Laplace \mu + (\mathcal{P} \sigma - \mathcal{A})(\varphi + 1).
\end{align*}
\end{itemize}
We always take $\sigma_{B} = 1$, $\beta = 0.1$, $\mathcal{P} = 0.1$, $\mathcal{A} = 0$, $\mathcal{C} = 1$, $\chi_{\varphi} = 10$, $\eps = 0.01$, $m_{0} = 1$ $K = 0.01$, $\lambda = 0.03$ and $D = 1$. 

\begin{figure}[h]
\begin{center}
\subfigure{\includegraphics[width=.28\textwidth,angle=0]{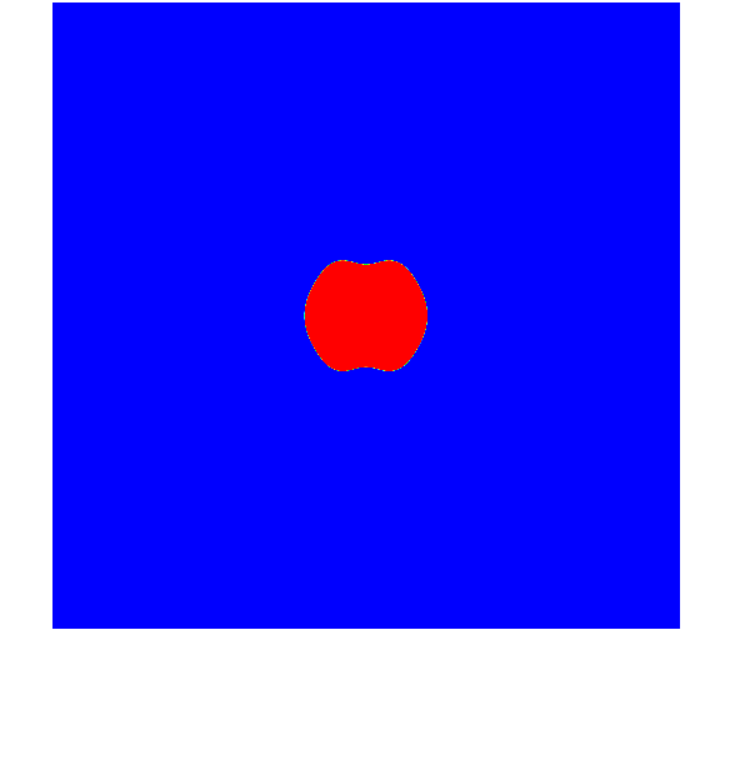}}\hspace{4mm}
\subfigure{\includegraphics[width=.28\textwidth,angle=0]{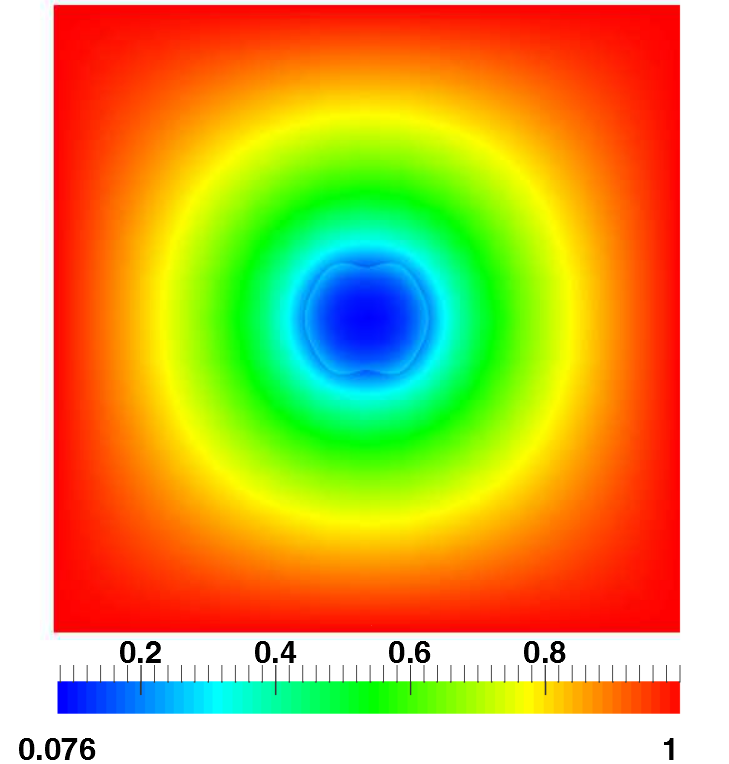}}\hspace{4mm}
\subfigure{\includegraphics[width=.28\textwidth,angle=0]{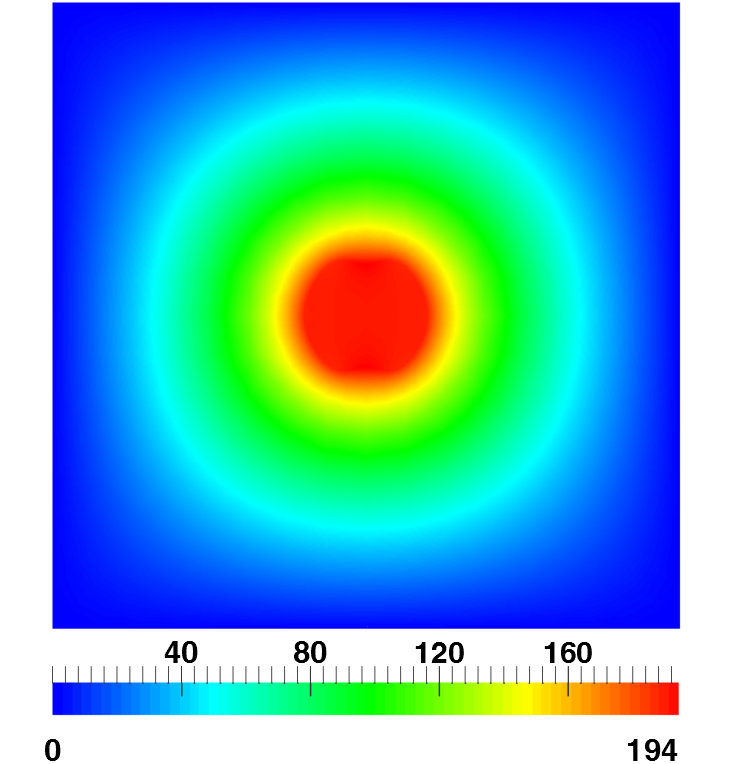}}
\caption{Solutions of (\ref{eq:hhh}) for case (1), $\varphi_{h}^{n}(x)$ (left), $\sigma_{h}^{n}(x)$ (centre) $p_{h}^{n}$ (right) at $t=1.5$.}
\label{f:vel1}
\end{center}
\end{figure}

\begin{figure}[h]
\begin{center}
\subfigure{\includegraphics[width=.28\textwidth,angle=0]{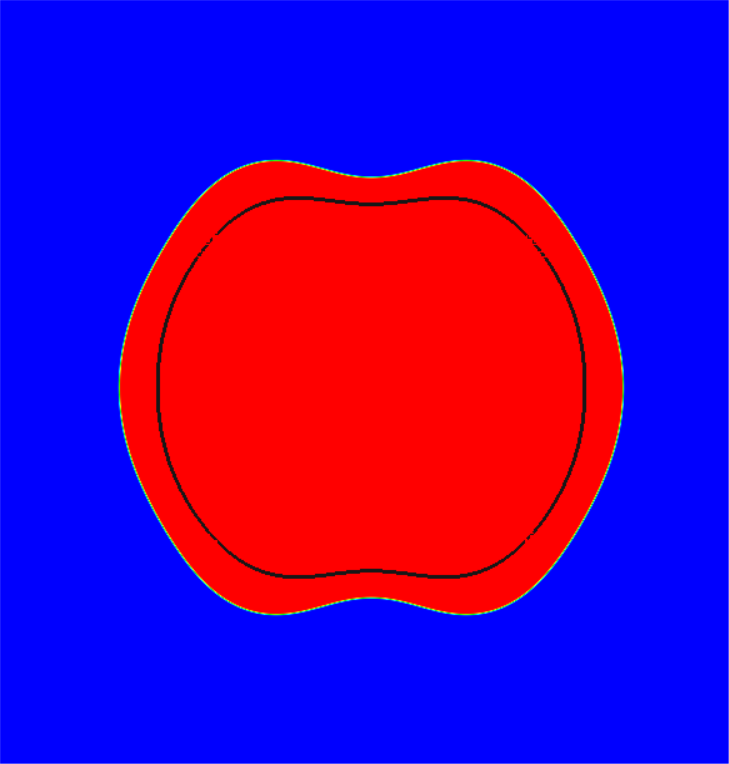}}\hspace{4mm}
\subfigure{\includegraphics[width=.28\textwidth,angle=0]{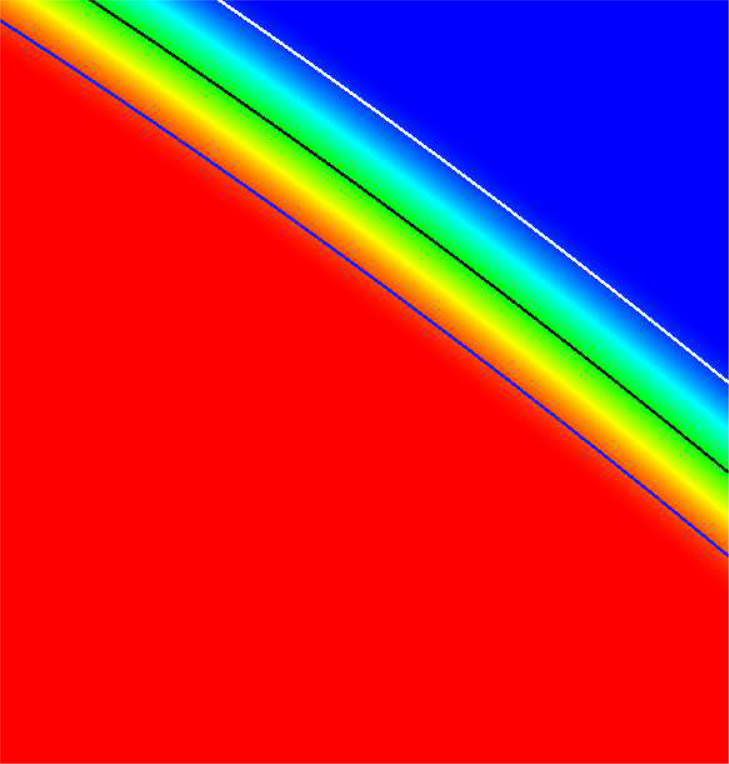}}%\hspace{4mm}
\caption{The left plot displays $\varphi_{h}^{n}(x)$ from (\ref{eq:hhh}) for case (2) at $t=1.5$ with the black line denoting the zero level line of $\varphi_{h}^{n}(x)$ from (\ref{FE:Cristini}).  The right plot displays the zero level lines of $\varphi_{h}^{n}(x)$ from (\ref{eq:hhh}) for case (1) (black line), case (2) (white line) and case (3) (blue line) at $t=1.5$.}
\label{f:vel2}
\end{center}
\end{figure}

In Figure \ref{f:vel1} we display the solutions of the Darcy flow model (\ref{eq:hhh}) for case (1) at $t = 1.5$; the left plot is of $\varphi$, the centre plot of $\sigma$ and the right plot is of $p$.  In the left plot of Figure \ref{f:vel2} we display a zoomed in plot of $\varphi$ at $t=1.5$ obtained from the Darcy flow model (\ref{eq:hhh}) for case (2) together with the zero level line of $\varphi_{h}^{n}(x)$ (depicted in black) with $K = \alpha = 0$ (which is equivalent to (\ref{FE:Cristini})).  In the right plot 
we show the influence of $\alpha$ on the position of the tumour, we display a zoomed in plot of the solution $\varphi$ at $t=1.5$; the black, white, and blue lines are the zero level lines of $\varphi$ for case (1), case (2) and case (3), respectively.  One observes that the model variant with Darcy flow enhances the growth velocities of the tumour.  In addition, the velocity is largest when the density of the tumour is smaller than the density of the healthy cells.

\section*{Acknowledgements}
The authors gratefully acknowledge the support of the Regensburger Universit\"{a}tsstiftung Hans Vielberth.  The fourth author is supported by the Engineering and Physical Sciences Research Council, UK grant (EP/J016780/1) and the Leverhulme Trust  Research Project Grant (RPG-2014-149).

\bibliographystyle{plain}
\bibliography{GLSSTumour}

\end{document}